\documentclass[11pt]{article}

\usepackage{mystyle}
\usepackage{contributors}
\usepackage{mymath}
\usepackage{jlcode}

\titleformat{\paragraph}
{\normalfont\small\bfseries}{\theparagraph}{1em}{}
\titlespacing*{\paragraph}
{0pt}{3.25ex plus 1ex minus .2ex}{1.5ex plus .2ex}

\title{Differentiable Programming for Differential \\ Equations: A Review}

\date{}

\begin{document}
\maketitle


\begin{abstract}
{\footnotesize
The differentiable programming paradigm is a cornerstone of modern scientific computing. 
It refers to numerical methods for computing the gradient of a numerical model's output. 
Many scientific models are based on differential equations, where differentiable programming plays a crucial role in calculating model sensitivities, inverting model parameters, and training hybrid models that combine differential equations with data-driven approaches.
Furthermore, recognizing the strong synergies between inverse methods and machine learning offers the opportunity to establish a coherent framework applicable to both fields.
Differentiating functions based on the numerical solution of differential equations is non-trivial. 
Numerous methods based on a wide variety of paradigms have been proposed in the literature, each with pros and cons specific to the type of problem investigated.
Here, we provide a comprehensive review of existing techniques to compute derivatives of numerical solutions of differential equations.
We first discuss the importance of gradients of solutions of differential equations in a variety of scientific domains.
Second, we lay out the mathematical foundations of the various approaches and compare them with each other. 
Third, we cover the computational considerations and explore the solutions available in modern scientific software.
Last but not least, we provide best-practices and recommendations for practitioners. 
We hope that this work accelerates the fusion of scientific models and data, and fosters a modern approach to scientific modelling.
}

\end{abstract}
\newpage


\tableofcontents

\mbox{}
\vfill
\begin{quote}
    \textit{This manuscript was conceived with the goal of shortening the gap between developers and practitioners of differentiable programming applied to modern scientific machine learning. 
    With the advent of new tools and new software, it is important to create pedagogical content that allows the broader community to understand and integrate these methods into their workflows. 
    We hope this encourages new people to be an active part of the ecosystem, by using and developing open-source tools. 
    This work was done under the premise of \textbf{open-science from scratch}, meaning all the contents of this work, both code and text, have been in the open from the beginning and that any interested person can contribute to the project. 
    }
\end{quote}
\newpage


\section{Introduction}
Models based on differential equations (DEs), including ordinary differential equations (ODEs), partial differential equations (PDEs), and Stochastic Differential Equations (SDEs), play a central role in describing the behavior of dynamical systems in applied, natural, and social sciences.
For instance, DEs are used for the modelling of the dynamics of the atmospheric and ocean circulation in climate science, for the modelling of ice or mantle flow in solid Earth geophysics, and for the modelling of the spatio-temporal dynamics of species abundance in ecology.  
For centuries, scientists have relied on theoretical and analytical methods to solve DEs. 
By allowing to approximate the solutions of large, nonlinear DE-based models, numerical methods and computers have lead to fundamental advances in the understanding and prediction of physical, biological, and social systems, among others \cite{Dahlquist_1985, hey2009, Rude:2018jv}.


Quantifying how much the output of a DE-based model changes with respect to its input parameters is fundamental to many scientific computing and machine learning applications, including optimization, sensitivity analysis, Bayesian inference, inverse methods, and uncertainty quantification, among many \cite{Razavi.2021}. 
Mathematically, quantifying this change involves evaluating the gradient of the model, i.e., calculating a vector whose components are the partial derivatives of the model evaluated at the model parameter values. 
In sensitivity analysis, gradients are crucial for comprehending the relationships between model inputs and outputs, assessing the influence of each parameter, and evaluating the robustness of model predictions. 
In optimization and inverse modelling, where the goal is to fit models to data and/or invert for unknown or uncertain parameters, gradient-based methods are more efficient at finding and converging to a minimum than gradient-free methods.
In Bayesian inference, gradient-based sampling strategies are often better at estimating the posterior distribution than gradient-free methods (e.g. \textcite{neal2011mcmc}).
Therefore, accurately determining model gradients is essential for robust model understanding and effective data assimilation that leverage strong physical priors while offering flexibility to adapt to observations. 
This is very appealing in fields such as computational physics, geophysics, and biology, to mention a few, where there is a broad literature on DE-based models.
The techniques used to compute these gradients fall within the framework of differentiable programming.

Differentiable programming (DP) refers to a programming paradigm that enables the end-to-end computation of gradients/sensitivities of a computer program with respect to input variables or parameters \cite{Shen_diff_modelling, Innes_Zygote, blondel2024elements}.
DP builds on top of computer science, mathematics, and statistics, among others \cite{blondel2024elements}.
Just as neural networks have been proven to be very flexible in learning nonlinear patterns from data using algorithmic differentiation, DE-based models require DP in order to explore the parameter space.
One central component behind the DP paradigm is automatic differentiation.
The set of tools known as automatic or algorithmic differentiation (AD) aims to compute derivatives of a model rendered on a computer by applying the chain rule to the sequence of unit operations that constitute a computer program \cite{Griewank:2008kh, Naumann.2011}. 
The premise is simple: every computer program is ultimately an algorithm described by a nested concatenation of elementary algebraic operations, such as addition and multiplication.
These operations are individually easy to differentiate, and their composition can be easily differentiated using the chain rule \cite{Giering_Kaminski_1998}.
During the last decades, reverse mode AD, also known as backpropagation, has enabled the fast growth of deep learning by efficiently computing gradients of neural networks \cite{griewank2012invented}.
Some authors have recently suggested DP as the bridge between modern machine learning and traditional scientific models \cite{Ramsundar_Krishnamurthy_Viswanathan_2021, Shen_diff_modelling, Gelbrecht-differential-programming-Earth, rackauckas2021generalized}. 
More broadly than AD, DP tools for DE-models further include forward sensitivity and adjoint methods that compute the gradient by relying on an auxiliary set of differential equations.

The development of DP for DE-based models has a long tradition across different scientific communities.
In statistics, gradients of the likelihood function of DE-based models enable inference on the model parameters \cite{ramsay2017dynamic}. 
In numerical analysis, sensitivities quantify how the solution of a differential equation fluctuates with respect to certain parameters. 
This is particularly useful in optimal control theory, where the goal is to find the optimal value of some control (e.g. the shape of a wing) that minimizes a given loss function \cite{Giles_Pierce_2000}. 
In recent years, there has been an increasing interest in designing machine learning workflows that include constraints in the form of DEs and are trained using gradient descent techniques.
This emerging sub-field is usually referred to as physics-based or physics-informed machine learning \cite{Karniadakis_Kevrekidis_Lu_Perdikaris_Wang_Yang_2021, thuerey2021pbdl, vadyala2022review}.

The need for model gradients is even more critical as the total number of parameters and the expressivity of the model increase, especially when dealing with highly non-linear processes.
In such circumstances, the curse of dimensionality renders gradient-free optimization and sampling methods computationally intractable. 
This is the case in inverse methods \cite{Tarantola:2007wu, Ghattas.2021} and in machine learning applications \cite{LeCun2015}, where highly parametrized regressor functions (e.g., neural networks) are used to approximate unknown non-linear function.
Furthermore, for stochastic forward models, the intractability of the likelihood function represents a major challenge for statistical inference.
The integration of DP has provided new tools for resolving complex simulation-based inference problems \cite{Cranmer_Brehmer_Louppe_2020}.

Computing gradients of functions, represented by DE-based simulation codes, with respect to their (high-dimensional) inputs is challenging due to the complexities in both the mathematical framework and the software implementation involved.
Except for a small set of particular cases, most DEs require numerical methods to approximate their solution, which means that solutions cannot be differentiated analytically. 
Furthermore, numerical solutions introduce approximation errors. 
These errors can be propagated to the computation of the gradient, leading to inaccurate or inconsistent gradient values.  
In addition to these complexities, the existence of a broad family of numerical methods, each one of them with different advantages depending on the DE \cite{hairer-solving-1, hairer-solving-2}, means that the tools developed to compute gradients need to be universal enough in order to be applied to all or at least to a large set of numerical solvers.

There exists a large family of methods to compute derivatives of DE-based models. 
The differences between methods to compute derivatives arise from their mathematical formulation, numerical stability, and their computational implementation. 
They can be roughly classified as continuous (differentiate-then-discretize) or discrete (discretize-then-differentiate) and forward or reverse. 
Different methods guarantee different levels of accuracy, have different computational complexity, and require different trade-offs between run time and memory usage. 
These properties further depend of the total number of parameters and size of the DE. 
Despite their independent success, integrating DP with DE-based models remains a major challenge in high-performance scientific computing \cite{Naumann.2011}, particularly in relation to research software engineering practices to ensure efficiency, scalability, sustainability and reproducibility.
\\ \\
This paper presents a comprehensive review of methods for calculating derivatives of functions of the numerical solution of differential equations, with a focus on efficiently computing gradients. 
We review differentiable programming methods for differential equations from three different perspectives: a domain science perspective (Section \ref{sectopn:motivation}), a mathematical perspective (Section \ref{section:methods}), and a computer science perspective (Section \ref{sec:computational-implementation}). 
In Section \ref{sectopn:motivation} we introduce some of the applications of DP for the modelling of complex systems in the natural and social sciences. 
In Section \ref{section:methods} we present a coherent mathematical framework to understand the theoretical differences between the DP methods.
In Section \ref{sec:computational-implementation} we show how these methods are computationally implemented and their numerical advantages and disadvantages.
For simplicity, all the methods introduced in Sections \ref{section:methods} and \ref{sec:computational-implementation} focus exclusively on first-order ODEs. 
How these methods generalize to other DE-based models, including PDEs and SDEs, is discussed in Section \ref{section:generalization}.
We conclude the paper with a series of recommendations in Section \ref{section:recomendations}.
By providing a common framework across all methods and applications, we hope to facilitate the development of scalable, practical, and efficient differentiable DE-based models.

\section{Scientific motivation: A domain science perspective}
\label{sectopn:motivation}
Mechanistic (or process-based) models play a central role in a wide range of scientific disciplines. 
They consist of mathematical descriptions of physical mechanisms that include the modelling of dependencies between components of the system under consideration. 
These mathematical representations often take the form of DEs. 
Together with the numerical methods to approximate their solutions, DEs have led to fundamental advances in the understanding and description of real-world systems.

DEs usually depend on inputs or parameters that change the obtained solutions. 
While direct or forward modelling usually refers to understanding how these parameters map into solutions or observations of the DE-based model, 
the overarching goal of inverse modelling is to find a set of optimal model parameters that minimizes an objective or cost function quantifying the misfit between observations and the simulated state.
This goal can be achieved via the construction the corresponding adjoint model that computes the gradient of the objective function with respect to all inputs \cite{Vadlamani.2020, Givoli_2021}.
Gradient-based optimization enables the inversion of optimal values of the unknown or uncertain inputs.
Depending on the nature of the inversion, we may distinguish between the following cases:
\begin{itemize}
    \item[$\blacktriangleright$] \textbf{Initial conditions.} Inverting for uncertain initial conditions, which, when integrated using the model, leads to an optimal match between the observations and the simulated state (or diagnostics); variants thereof are used for optimal forecasting.
    \item[$\blacktriangleright$] \textbf{Boundary conditions.} Inverting for uncertain surface (e.g., interface fluxes), bottom (e.g., bed properties), or lateral (e.g., open boundaries of a limited domain) boundaries, which, when used in the model, produce an optimal match of the observations. Variants thereof are used in tracer or boundary (air-sea) flux inversion problems, e.g., related to the global carbon cycle.
    \item[$\blacktriangleright$] \textbf{Model parameters.} Inverting for uncertain model parameters amounts to an optimal model calibration problem. As a \textit{learning of optimal parameters from data} problem, it is the closest to machine learning applications. Parametrization is a special case of parameter inversion, where a parametric function (e.g., a neural network) is used to approximate processes. 
\end{itemize}
Besides the use of sensitivity methods for optimization, inversion, estimation, or learning, gradients have also proven powerful tools for computing comprehensive sensitivities of quantities of interest; computing optimal perturbations (in initial or boundary conditions) that lead to maximum, non-normal amplification of specific norms of interest; and characterizing and quantifying uncertainties by way of second derivative (Hessian) information.
The availability of second derivatives further helps to improve the convergence rates of optimization algorithms.

In the following subsections, we present selected examples from a wide range of scientific communities where DP techniques have been used for the modelling of systems described using DEs. 
Rather than providing a exhaustive overview, we hope some of the applications enumerated here can give domain experts to a quick glance of how DP has been used in their field, as well as showcase the different historical scientific trajectories that these methods have had in each community. 

\subsection{Machine learning}

In recent years the use of machine learning methods has become more popular in many scientific domains (e.g. \textcite{rasp2018}, \textcite{pichler2023}, \textcite{meuwly2021machine}, \textcite{borowiec2022}, \textcite{lai2024machine_final}). 
By learning nonlinear patterns from large datasets at multiple levels of abstraction \cite{LeCun2015}, these methods are highly flexible with respect to inputs and outputs required, and can be exploited by many different domain-specific problems.
In contrast to purely statistical models, the process knowledge embedded in the structure of mechanistic models renders them more robust for predicting dynamics under different conditions.
The fields of mechanistic modelling and statistical modelling have mostly evolved independently due to several reasons \cite{zdeborova_understanding_2020}. 
On the one hand, domain scientists have often been reluctant to learning about machine learning methods, judging them as opaque black boxes, unreliable, and not respecting domain-established knowledge \cite{Coveney:2016eb}. 
On the other hand, the field of machine learning has mainly been developed around data-driven applications, without including any a priori physical knowledge. 
However, there has been an increasing interest in making mechanistic models more flexible, as well as introducing domain-specific or physical constraints and interpretability in machine learning models. 
This sub-field, usually known as physics-informed machine learning, refers to the collection of machine learning techniques that explicitly introduce biases to satisfy certain physical constraints. 
These biases can be forced by the design of algorithms that include symmetries, conservation laws, and constraints in the form of DEs \cite{Karniadakis_Kevrekidis_Lu_Perdikaris_Wang_Yang_2021}. 
It includes methods that numerically solve DEs, such as physics-informed neural networks \cite{PINNs_2019}, biology-informed neural networks \cite{Yazdani2020,Lagergren_Nardini_Baker_Simpson_Flores_2020}, NeuralPDEs \cite{Zubov_McCarthy_Ma_Calisto_Pagliarino_Azeglio_Bottero_Luján_Sulzer_Bharambe_et}, and mesh-free methods for solving high-dimensional PDEs \cite{boussange2023a}. 
Furthermore, there has been an increased interest in augmenting DE-based models by embedding a rich family of parametric functions (e.g., neural networks) inside the DE.  
This approach is known as universal differential equations \cite{rackauckas2020universal}, which also include the case of neural ordinary differential equations \cite{chen_neural_2019} and neural stochastic differential equations \cite{li2020scalable}.


\subsection{Computational physics and optimal design}

There is a long tradition of computational physics models based on adjoint methods and AD pipelines, where sensitivity methods have been used for optimal design and optimal control since the 1960s \cite{lions1971optimal}. 
These models are often based on PDEs and are applied in various fields to improve engineering designs or model parameters with respect to some objective function. 
The models can involve thousands of parameters, and they require efficient derivative calculations for the use of gradient-based optimizers such as quasi-Newton methods \cite{nocedal1999numerical}. 
Both discrete and continuous adjoints methods, which we will introduce in Sections \ref{section:discrete-adjoint} and \ref{section:continuous-adjoint}, respectively, have been used extensively, each having different benefits depending on the application.

\subsubsection{Computational fluid dynamics}

DP methods, including AD and adjoint methods, have been crucial in advancing computational fluid dynamics (CFD) applications \cite{KENWAY2019100542}. 
These techniques have been employed in optimizing the aerodynamics of aircraft for drag reduction, or for weight reduction in aircraft design, leading to significant fuel savings and enhanced performance \cite{jameson2003aerodynamic}. 
In aeroacoustic designs, adjoint methods can be used to minimize noise emissions (e.g. \textcite{FREUND201054}). 
The objective function in these applications typically relates to performance metrics or cost considerations, and a wide array of design parameters can be optimized.
Entire geometries can be parameterized for \emph{shape optimization}, enabling the refinement of complex structures like airfoils, which are critical for aerodynamic efficiency. 
Pironneau introduced fundamental methods for shape optimization in fluid mechanics \cite{Pironneau_1974}, and Jameson developed adjoint-based optimization methods that significantly improved aerodynamic designs \cite{Jameson_1988}. 
For comprehensive reviews of shape optimization using adjoint methods, we refer to \textcite{Giles_Pierce_2000} and \textcite{mohammadi2009applied}. 
Adjoints have also been used for topology optimization \cite{allaire2014shape}.

For aerospace applications, adjoint methods have been used to design supersonic aircraft, enhancing performance and reducing sonic boom impacts \cite{hu2010supersonic,fike2013multi}. 
Entire aircraft configurations have also been optimized using adjoint methods \cite{chen2016aerodynamic}. 
Beyond aerospace, other significant applications include optimizing ship hull designs to reduce drag and improve fuel efficiency \cite{kroger2018ships}, the aerodynamic shaping of cars to enhance speed and stability \cite{Othmer2014caraerodynamics}, and the design of wind turbines to maximize energy capture and structural resilience \cite{dhert2017aerodynamic}.

\subsubsection{Quantum physics}

Quantum optimal control has applications spanning a broad spectrum of quantum systems. 
Optimal control methods have been used to optimize pulse sequences, enabling the design of high-fidelity quantum gates, and the preparation of complex entangled quantum states. 
Typically, the objective is to maximize the fidelity to a target state or unitary operation, accompanied by additional constraints or costs specific to experimental demands. 
The predominant control algorithms are gradient-based optimization methods, and rely on the computation of derivatives for solutions of the ODEs modeling the time evolution of the quantum system. 
In cases where the analytical calculation of a gradient is impractical, numerical evaluation using AD becomes a viable alternative~\cite{jirari:2009, leung:2017, jirari2019quantum, schaefer:2020}. 
Specifically, AD streamlines the adjustment to diverse objectives or constraints, and its efficiency can be enhanced by employing custom derivative rules for the time propagation of quantum states as governed by solutions to the Schrödinger equation \cite{goerz:2022}. 
Moreover, sensitivity methods facilitate the design of feedback control schemes necessitating the differentiation of solutions to stochastic differential equations \cite{schaefer:2021}.

\subsubsection{Other applications}

Adjoint methods have also been applied successfully to a wide range of other computational physics problems. 
In particle physics, they enable precise parameter estimation and simulation improvements \cite{Dorigo.2022}.
In quantum chemistry, adjoint methods can be used to optimize molecular structures and reaction pathways (e.g. \textcite{Arrazola.2021}). 
The design of nanophotonic devices, such as photonic crystals and waveguides, has been significantly advanced through these techniques \cite{Molesky_Lin_Piggott_Jin_Vucković_Rodriguez_2018}. 
Electromagnetic applications, including the optimization of antenna designs and microwave circuits, benefit from the fine-tuning capabilities provided by adjoint methods \cite{Georgieva_Glavic_Bakr_Bandler_2002}. 
Stellarator coil design for nuclear fusion reactors is another important area, where adjoint methods contribute to optimizing magnetic confinement configurations \cite{McGreivy_stellarator_2021}.
Sensitivity analysis methods are also very popular in topological and structural design \cite{min1999optimal, van2005review}.

\subsection{Geosciences}

Many geoscientific phenomena are governed by conservation laws along with a set of empirical constitutive laws and subgrid-scale parameterization schemes. 
Together, they enable efficient description of the system's spatio-temporal evolution in terms of a set of PDEs.
An example is geophysical fluid dynamics \cite{Vallis.2016}, describing geophysical properties of many Earth system components, such as the atmosphere, ocean, land surface, and glaciers.
In such models, calibrating model parameters is extremely challenging, due to observational data being sparse in both space and time, heterogeneous, and noisy; and computational models involving high-dimensional parameter spaces, often on the order of $O(10^3) - O(10^8)$.
Moreover, many existing mechanistic models can only partially describe observations, with many detailed physical processes being ignored or poorly parameterized.

\subsubsection{Numerical weather prediction}
\label{section:meteorlogy}

Adjoint methods have played an important role in numerical weather prediction (NWP) to infer initial conditions that minimize the misfit between simulations and weather observations \cite{Lewis.1985,Talagrand.1987,Courtier.1987}, with the value of second-derivative information also being recognized \cite{Dimet.2002}. 
This led to the development of the four-dimensional variational (4D-Var) technique at the European Centre for Medium-Range Weather Forecasts (ECMWF) as one the most advanced data assimilation approaches \cite{Rabier.1992,Rabier:2000uu}.
Within the framework of transient non-normal amplification \cite{Farrell.1988,Farrell:1996jx}, derivative models have been used to infer patterns in initial conditions that over a finite time contribute to maximum uncertainty growth in forecasts \cite{Palmer:1994br,Buizza:1995in} and to infer the so-called Forecast Sensitivity-based Observation Impact (FSOI) \cite{Langland:2004jo}.
Except for early research applications \cite{Park.1996,Park.2000} and for experimental purposes \cite{Giering.2006}, AD has not been widely used in the development of adjoint models in NWP.
Instead, the adjoint code has been, for the most part, derived and implemented manually.

\subsubsection{Oceanography}

The recognition of the benefit of adjoint methods for use in data assimilation in the ocean coincided roughly with that in NWP \cite{Thacker:1988kp,Thacker:1989jf,Tziperman.1989,Tziperman:1992hg}. 
An important detail is that their work already differed from the 4D-Var problem of NWP (Section \ref{section:meteorlogy}) in that sensitivities were computed not only with respect to initial conditions but also surface boundary conditions.
Similar to the work on calculating singular vectors in the atmosphere, 
the question of El Ni\~no \cite{Moore.1997ah,Moore.1997}
and Atlantic Meridional Overturning (AMOC) \cite{Zanna:2012dw}
predictability invited model-based singular vector computations into ocean models. 
More recent data assimilation frameworks with fully hand-written adjoint codes include the NEMO model \cite{Weaver.2003,Vidard:2015kj} and the ROMS model \cite{Moore:2004fk,Moore:2011bc}.

Fully AD-based adjoint ocean models were developed beginning in the late 1990s in the context of state and parameter estimation \cite{Marotzke:1999, Heimbach.2005,Stammer.2002}.
Rigorous exploitation of AD enabled the extensions to vastly improved model numerics \cite{Forget.2015m9i} and coupling to other Earth system components.
Unlike NWP-type 4D-Var, it also enabled the extension to the problem of parameter calibration from observations \cite{Ferreira.2005,Stammer:2005dw,Liu:2012jd}. 
Arguably, this work heralded much of today's efforts in \textit{online learning} of parameterization schemes.
The desire to make AD for Earth system models written in Fortran (to date the vast majority) has spurned the development of AD tools with powerful reverse modes, both commercial \cite{Giering_Kaminski_1998,Giering.2006} and open-source \cite{Naumann:2006wq,Utke:2008ko,Hascoet.2013,Gaikwad.2025}.

\subsubsection{Climate science}

Similar goals driving the use of sensitivity information in NWP (optimal initial conditions for forecasting) or ocean science (state and parameter estimation) apply in the world of climate modeling.
The recognition that good initial conditions 
will lead to improved forecasts on subseasonal, seasonal, interannual, or even decadal time scales underlies major community efforts \cite{Meehl.2021}. 
However, there has been a lack of use of gradient information to achieve optimal initialization for coupled Earth system models. 
One conceptual challenge is the presence of multiple timescales in the coupled system and the utility of gradient information beyond many synoptic time scales in the atmosphere and ocean \cite{lea2000sensitivity,Lea:2002cv}.
Nevertheless, efforts are underway to enable adjoint-based parameter estimation of coupled atmosphere-ocean climate models, with AD again playing a crucial role 
\cite{Blessing.2014,Lyu.2018,Stammer:2018de}.
Additionally, recognizing the power of DP, efforts are also targeting the development of neural atmospheric general circulation models in \texttt{JAX}, which combine a differentiable dynamical core with neural operators as surrogate models of unresolved physics \cite{Kochkov.2024}.

\subsubsection{Glaciology}

Due to the difficulty of taking direct measurements of internal and basal rheological processes of glaciers and ice sheets, inverse methods based on adjoint models have been widely used to study them, following the pioneering work by \textcite{macayeal1992basal}. 
Since then, the adjoint method has been applied to many different studies investigating parameter and state estimation \cite{Vieli_Payne_2003, goldberg2013parameter}, ice volume sensitivity to basal, surface and initial conditions \cite{heimbach2009greenland}, inversion of initial conditions \cite{mosbeux2016comparison} or inversion of basal friction \cite{Petra.2012, morlighem2013inversion, Brinkerhoff_Johnson_2013}.
These studies either derived the adjoint with a manual implementation or combined AD with hand-written adjoint solvers. 
The use of AD has become increasingly widespread in glaciology, paving the way for more complex modelling frameworks \cite{hascoet2018source, Gaikwad.2023}. 
The use of second-derivative (Hessian) information has also been recognized as a powerful approach for conducting rigorous uncertainty quantification in the context of ice sheet parameter inversion \cite{Petra.2014,Isaac.2015}.
Recently, DP has also facilitated the development of hybrid frameworks, combining numerical methods with data-driven models by means of universal differential equations \cite{BolibarSapienza_UDEs}. 
Alternatively, some other approaches have dropped the use of numerical solvers in favor of physics-informed neural networks, exploring the inversion of rheological properties of glaciers \cite{wang2022discovering} and accelerating ice thickness inversions and simulations by leveraging GPUs \cite{Jouvet_Cordonnier_Kim_Lüthi_Vieli_Aschwanden_2021, jouvet2023inversion}. 
For a recent comprehensive review of data assimilation in ice sheet modeling, with emphasis on adjoint methods, see \textcite{Morlighem.2023}.


\subsection{Biology and ecology}
\label{section:biology}

DE-based models have been broadly used in biology and ecology to model 
neural firing \cite{hodgkin1952quantitative}, 
the dynamics of genes and alleles \cite{Page2002}, 
immune and disease processes \cite{colijn2006high}, 
and biomass and energy fluxes and transformation at ecosystem levels (e.g. \textcite{Weng2015}, \textcite{Schartau2017}), 
among others.
Parameters in DE-based models are often estimated through direct laboratory experiments (e.g. \textcite{hodgkin1952quantitative}) although this process is costly and difficult, and may result in simulations failing at capturing real biological dynamics \cite{Schartau2017, Watts2001}. 
Alternatively, inverse modeling methods have a long history of application, both for parameter estimation (e.g., \textcite{ramsay2007parameter}, \textcite{Schartau2017}, \textcite{ding2000h}, \textcite{fussmann2000crossing}) and model selection \cite{alsos2023,pantel2023} in DE-models.

When parameters are inferred along with their uncertainties, they can be interpreted to better understand the strengths and effects of the processes under consideration (e.g., \textcite{Pontarp2019} and \textcite{Curtsdotter2019}). 
However, in high-dimensional models parameters are often non-identifiable \cite{transtrum2011geometry}. 
Model selection, which does not require parameter identifiability, involves deriving candidate models that embed competing hypotheses about causal processes and computing the relative evidence for each model given the data to discriminate between hypotheses \cite{Johnson2004, alsos2023}.
Methods for inverse modelling with DE models must effectively handle the typically large number of parameters and the nonlinearities of biological models \cite{transtrum2011geometry, Gabor2015}. 
Sensitivity-based inference frameworks are strong candidates for this task.
Sensitivity methods also play a crucial role in assessing model fit quality, specifying system states, detecting stochastic noise \cite{hooker2009forcing, hooker2015goodness, liu2023specification}, and designing experiments to optimize parameter estimation precision \cite{bauer2000numerical}.

While inverse modeling based on AD is seldom used in biology, its potential has recently been emphasized \cite{frank2022, alsos2023}. 
Approaches involving AD have been recently proposed to accommodate the specific requirements of biological and ecological models where key processes are often not accurately represented (see \textcite{Boussange2024}, \textcite{Lagergren_Nardini_Baker_Simpson_Flores_2020}).
In particular, AD is bound to play a major role in hybrid modeling approaches that integrate data-driven parameterization of specific components of DE-based models (e.g. \textcite{ramsay1996principal},  \textcite{cao2008estimating},   \textcite{chen2017network},  \textcite{dai2022kernel}).

\subsection{Computational finance}
\label{section:comp_finance}

Computational finance deals with simulating the price of a range of financial assets (e.g., complex portfolios of structured derivatives) and their associated risks.
The price of financial assets is simulated with stochastic differential equations (SDEs) that capture the fluctuating nature of the asset's volatility as a geometric Brownian motion (e.g., \textcite{voit2005statistical}).
Conducting robust risk (i.e., uncertainty) analysis is essential in order to quantify the assets' exposure to market fluctuations.
Conventional risk analysis relies on sampling methods based on Monte Carlo methods which determine the sensitivity of the price to a large set of inputs.
The sensitivities of the asset price to its inputs are referred to as financial \textit{Greeks} (high-order Greeks represent higher derivatives).

The combination of PDE solvers with Monte Carlo methods render the comprehensive calculation of Greeks required for risk management a formidable computational task. 
A method called (forward) pathwise sensitivity calculation in the context of SDEs was proposed by \textcite{glasserman1999fast} to compute accurate estimates of price sensitivities (see also \textcite{glasserman2004monte}).
Using the example of a Libor market model, \textcite{Giles:2006vr} showed that the adjoint method (i.e., backward pathwise sensitivity) could reduce the computational cost of computing Greeks by orders of magnitude. 
They already refer to the use of AD to highlight the reduction in computational complexity of the reverse mode. 
The role of AD, or rather adjoint algorithmic differentiation (AAD), is exposed in more detail in \textcite{Capriotti:2010td}. 
We review the underlying application of AAD to SDEs in more detail in Section \ref{section:sde}.
Another important application is that of model parameter \textit{calibration} applied to financial asset price modeling.
Computational finance has since seen a flurry of work on the development of sophisticated algorithms involving AAD approaches in the derivatives industry, as recently reviewed by \textcite{Capriotti.2024}.

\section{Methods: A mathematical perspective}
\label{section:methods}
\begin{figure}[t]
    \centering
    \includegraphics[width=0.80\textwidth]{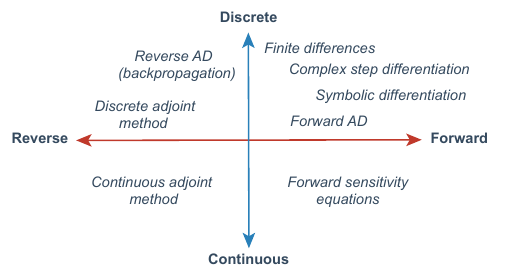}
    \caption{Schematic representation of the different methods available for differentiation involving differential equation solutions. These can be classified depending on whether they find the gradient by solving a new system of differential equations (\textit{continuous}) or whether they manipulate unit algebraic operations (\textit{discrete}). Additionally, these methods can be categorized depending on whether sensitivities are propagated from input to output (forward) or from output to input (reverse). There exist multiple intermediate methods to perform differentiation between the different axes of this figure. For example, elimination techniques can be used to efficiently evaluate derivatives (see Section \ref{sec:ad-further}) and continuous methods may rely on discrete methods.}
    \label{fig:scheme-all-methods}
\end{figure}

There is a large family of methods for computing gradients of functions based on solutions of DEs. 
Depending on the number of parameters and the characteristics of the DE (e.g., level of stiffness), they have different mathematical, numerical, and computational advantages.
These methods can be roughly classified as follows: 
\begin{itemize}
    \item[$\blacktriangleright$] \textit{Continuous} vs \textit{discrete}  methods
    \item[$\blacktriangleright$] \textit{Forward} vs \textit{reverse} methods
\end{itemize}
Figure \ref{fig:scheme-all-methods} displays a classification of some methods under this two-fold division. 

The \textit{continuous} vs \textit{discrete} distinction is one of mathematical and numerical nature. 
When solving for the gradient of a function of the solution of a differential equation, one needs to derive both a mathematical expression for the gradient (the differentiation step) and solve the differential equations using a numerical solver (the discretization step) \cite{bradley2013pde, Onken_Ruthotto_2020, FATODE2014, Sirkes_Tziperman_1997}. 
Depending on the order of these two operations, we refer to discrete methods (discretize-then-differentiate) or continuous methods (differentiate-then-discretize). 
In the case of \textit{discrete} methods, gradients are computed based on simple function evaluations of the solutions of the numerical solver (finite differences, complex step differentiation) or by manipulation of atomic operations inside the numerical solver (AD, symbolic differentiation, discrete adjoint method). 
It is worth noting that although both approaches are classified as discrete methods, their numerical properties are quite different.
In the case of \textit{continuous} methods, a new set of DEs is derived that allow the calculation of the desired gradient, namely the sensitivity (forward sensitivity equations) or the adjoint (continuous adjoint method) of the system.   
When comparing discrete to continuous methods, we are focusing, beyond computational efficiency, on the mathematical consistency of the method, that is, \textit{is the method estimating the right gradient?} 
Discrete methods compute the exact derivative of the numerical approximation to the loss function, but they do not necessarily yield to an approximation of the exact derivatives of the objective function \cite{Eberhard_Bischof_1996, Walther_2007}. 
On the other side, continuous methods lead to consistency error calculation of the sensitivity method, meaning that sensitivities are computed with a desired level of accuracy dictated by the discretization of the continuous inverse problem \cite{Keulen_Haftka_Kim_2005}. 

The distinction between \textit{forward} and \textit{reverse} regards whether sensitivities are computed for intermediate variables with respect to the input variable or parameter to differentiate (forward) or, on the contrary, we compute the sensitivity of the output variable with respect to each intermediate variable by defining a new adjoint variable (reverse). 
Mathematically speaking, this distinction translates to the fact that forward methods compute directional derivatives by mapping sequential mappings between tangent spaces, while reverse methods apply sequential mappings between co-tangent spaces from the direction of the output variable to the input variable (Section \ref{sec:vjp-jvp}).   
In all forward methods the DE is solved sequentially and simultaneously with the directional derivative during the forward pass of the numerical solver. 
On the contrary, reverse methods compute the gradient by solving a new problem that moves in the opposite direction as the original numerical solver.
In DE-based models, intermediate variables correspond to intermediate solutions of the DE.
In the case of ODEs and time-dependant PDEs, most numerical methods solve the DE by progressively moving forward in time, meaning that reverse methods solve for the gradient moving backwards in time. 
In other words, forward methods compute directional derivatives as they simultaneously solve the original DE, while reverse methods compute adjoints as they solve the problem from output to input.

As discussed in the following sections, forward methods are very efficient for problems with a small number of parameters we want to differentiate with respect to. 
Conversely, reverse methods, though more efficient for a large number of parameters, incur greater memory costs and computational overhead which need to be overcome using different performance tricks. 
With the exception of finite differences and complex step differentiation, the rest of the forward methods (i.e. forward AD, forward sensitivity equations, symbolic differentiation) compute the full sensitivity of the differential equation, which can be computationally expensive or intractable for large systems. 
Conversely, reverse methods are based on the computation of intermediate variables, known as the adjoint or dual variables, that cleverly avoid the unnecessary calculation of the full sensitivity at the expense of larger memory cost \cite{Givoli_2021}. 
Importantly for our discussion, other modes of differentiation via elimination techniques can achieve better performance, as we will discuss later in the manuscript.


The rest of this section is organized as follows. 
We first introduce some basic mathematical notions to facilitate the discussion of the DP methods (Section \ref{section:preliminaries}).
We then mathematically formalize each of the methods listed in Figure \ref{fig:scheme-all-methods}.
We finally discuss the mathematical foundations of these methods in \ref{section:compatison-math} with a comparison of some mathematical foundations of these methods. 

\subsection{Preliminaries}
\label{section:preliminaries}
Consider the first-order ODE given by
\begin{equation}
 \frac{du}{dt} = f(u, \theta, t)
 \label{eq:original_ODE}
\end{equation}
subject to the initial condition $u(t_0) = u_0$, where $u \in \mathbb{R}^n$ is the unknown solution vector of the ODE, $f: \R^n \times \R^p \times \R \mapsto \R^n$ is a function that depends on the state $u$, $\theta \in \mathbb R^p$ is a vector parameter, and $t \in [t_0, t_1]$ refers to time.
Here, $n$ denotes the size of the ODE and $p$ the number of parameters.
Except for a minority of functions $f(u,\theta, t)$ and initial conditions, solutions to Equation \eqref{eq:original_ODE} need to be computed using numerical solvers. 

\subsubsection{Numerical solvers for ODEs}
\label{section:intro-numerical-solvers}

Numerical solvers for the solution of ODEs or initial value problems can be classified as one-step methods, among which Runge-Kutta methods are the most widely used, and multi-step methods \cite{hairer-solving-1}.
Given an integer $s$, $s$-stage Runge-Kutta methods are defined by generalizing numerical integration quadrature rules as follows
\begin{align}
\begin{split}
    u^{m+1} 
    &= 
    u^m 
    + 
    \Delta t_m \sum_{i=1}^s b_i k_i \\
    k_i 
    &= 
    f \left(u^m + \sum_{j=1}^s a_{ij} k_j , \, \theta , \, t_m + c_i \Delta t_m \right) \qquad i=1,2, \ldots, s,
    \label{eq:Runge-Kutta-scheme}
\end{split}
\end{align}
where $u^{m} \approx u(t_m)$ approximates the solution at time $t_m$, $\Delta t_m = t_{m+1}-t_m$, and $a_{ij}$, $b_i$, and $c_j$ are scalar coefficients with $i,j=1, 2,\ldots, j$. 
A Runge-Kutta method is called explicit if $a_{ij}=0$ for $i \leq j$; diagonally implicit if $a_{ij}=0$ for $i < j$; and fully implicit otherwise. 
Different choices of the number of stages $s$ and coefficients give different orders of convergence of the numerical scheme \cite{Butcher_Wanner_1996, Butcher_2001}. 

In contrast, multi-step linear solvers are of the form 
\begin{equation}
    \sum_{i=0}^{k_1} \alpha_{mi} u^{m-i} 
    =
    \Delta t_m \sum_{j=0}^{k_2} \beta_{mj} f(u^{m-j}, \theta, t_{m-j})
\end{equation}
where $\alpha_{mi}$ and $\beta_{mj}$ are numerical coefficients, with $k_1$ and $k_2$ natural numbers \cite{hairer-solving-1}.
In most cases, including the Adams methods and backwards differentiation formulas (BDF), we have the coefficients $\alpha_{mi} = \alpha_i$ and $\beta_{mj}=\beta_j$, meaning that the coefficient do not depend on the iteration $m$. 
Notice that multi-step linear methods are linear in the function $f$, which is not the case in Runge-Kutta methods with intermediate evaluations \cite{ascher2008numerical}.
Explicit methods are characterized by $\beta_{m0} = 0$ and are easy to solve by direct iterative updates. 
For implicit methods, the usually non-linear equation 
\begin{equation}
    g_m(u^m; \theta) = u^m - \Delta t_m  \beta_{m0} f(u^m, \theta, t_m) - \bar{\alpha}_m = 0,
    \label{eq:solver-constriant-example}
\end{equation}
with $\bar{\alpha}_m$ a computed coefficient that includes the information of all past iterations, can be solved using predictor-corrector methods \cite{hairer-solving-1} or iteratively using Newton's method and preconditioned Krylov solvers at each nonlinear iteration \cite{SUNDIALS-hindmarsh2005sundials}.  

When choosing a numerical solver for differential equations, one crucial factor to consider is the stiffness of the equation.
Stiffness encompasses various definitions, reflecting its historical development and different types of instabilities \cite{Dahlquist_1985}.
Two definitions are noteworthy:
\begin{enumerate}
    \item[$ \blacktriangleright$] Stiff equations are equations for which explicit methods do not work and implicit methods work better \cite{hairer-solving-2}.
    \item[$ \blacktriangleright$] Stiff differential equations are characterized by dynamics with different time scales \cite{hairer-solving-2, kim_stiff_2021}, also characterized by the phenomena of increasing oscillations \cite{Dahlquist_1985}.
\end{enumerate} 
Stability properties can be achieved by different means, for example, by using implicit methods or stabilized explicit methods, such as Runge–Kutta–Chebyshev \cite{van1980internal, hairer-solving-2}. 
When using explicit methods, smaller timesteps may be required to guarantee stability. 


Numerical solvers usually estimate internally a scaled error computed as 
\begin{equation}
    \text{Err}_\text{scaled}^{m}
    =
    \sqrt{
    \frac{1}{n} \sum_{i=1}^n \left( \frac{\text{err}_i^{m}}{\mathfrak{abstol} + \mathfrak{reltol} \, \times \, M^m_i} \right)^2 },
    \label{eq:internal-norm-wrong}
\end{equation}
with $\mathfrak{abstol}$ and $\mathfrak{reltol}$ the adjustable absolute and relative solver tolerances, respectively, $M_i^m$ is the maximum expected value of the numerical solution, and $\text{err}^m$ is an estimation of the numerical error at step $m$ \cite{hairer-solving-1}. 
Estimations of the local error $\text{err}^{m}$ can be based on two approximation to the solution based on embedded Runge-Kutta pairs \cite{Ranocha_Dalcin_Parsani_Ketcheson_2022, hairer-solving-1}, or in theoretical upper bounds provided by the numerical solver. 
In the first case, common choices for these include $M^m_i = \max \{ u_i^{m}, \hat u_i^{m} \}$ and $\text{err}_i^{m} = u_i^{m} - \hat u_i^{m}$, with $u^m$ and $\hat u^m$ two different approximations for $u(t_m)$, but these can vary between solvers. 

Modern solvers include stepsize controllers that pick $\Delta t_m$ as large as possible to minimize the total number of steps while preventing large errors by keeping $\text{Err}^m_\text{scaled} \leq 1$. 
One of the most used methods to archive this is the proportional-integral controller (PIC) that updates the stepsize according to 
\begin{equation}
    \Delta t_{m} = \eta \, \Delta t_{m-1}
    \qquad 
    \eta = w_{m}^{\beta_1 / q} w_{m-1}^{\beta_2 / q} w_{m-2}^{\beta_3 / q},
    \label{eq:PIC}
\end{equation}
with $w_{m} = 1 / \text{Err}_\text{scaled}^{m}$ the inverse of the scaled error estimates; $\beta_1$, $\beta_2$, and $\beta_3$ numerical coefficients defined by the controller; and $q$ the order of the numerical solver \cite{hairer-solving-2, Ranocha_Dalcin_Parsani_Ketcheson_2022}. 
If the stepsize $\Delta t_{m}$ proposed in Equation \eqref{eq:PIC} to update from $u^{m}$ to $u^{m+1}$ does not satisfy $\text{Err}_\text{scaled}^{m+1} \leq 1$, a new smaller stepsize is proposed. 
When $\eta < 1$ (which is the case for simple controllers with $\beta_2 = \beta_3 = 0$), Equation \eqref{eq:PIC} can be used for the local update. 
It is also common to restrict $\eta \in [\eta_\text{min}, \eta_\text{max}]$ so the stepsize does not change abruptly \cite{hairer-solving-1}.

\subsubsection{What to differentiate?}

In most applications, the need for differentiating the solution of ODEs stems from the need to obtain the gradient of a function $L(\theta) = L(u(\cdot, \theta))$ with respect to the parameter $\theta$, where $L$ can denote:
\begin{itemize}
    \item[$ \blacktriangleright$] \textbf{A loss function or an empirical risk function}. This is usually a real-valued function that quantifies the level of agreement between the model prediction and observations. Examples of loss functions include the squared error
    \begin{equation}
         L(\theta) = \frac{1}{2} \left \| u(t_1; \theta) - u^{\text{target}}(t_1) \right \|_2^2,
         \label{eq:quadratic-loss-function}
    \end{equation}
    where $u^{\text{target}}(t_1)$ is the desired target observation at some later time $t_1$, and $\| \cdot \|_2$ is the Euclidean norm.
    More generally, we can evaluate the loss function at points of the time series for which we have observations, 
    \begin{equation}
        L(\theta) 
        = 
        \frac{1}{2} \sum_{i=1}^N 
        \, w_i 
        \left \| u(t_i; \theta) - u^{\text{target}}(t_i) \right \|_2^2,
        \label{eq:quadratic-loss-point}
    \end{equation}
    with $w_i$ some arbitrary non-negative weights.
    More generally, misfit functions used in optimal estimation and control problems are composite maps from the parameter space $\theta$ via the model's state space (in this case, the solution $u(t; \theta)$) to the observation space defined by a new variable $y(t) = H(u(t; \theta))$, where $H: \R^n \mapsto \R^o$ is a given function mapping the latent state to observational space \cite{1975-Bryson-Ho-optimal-control}. 
    In these cases, the loss function generalizes to 
    \begin{equation}
        L(\theta) 
        =
        \frac{1}{2} 
        \sum_{i=1}^N
        \, w_i 
        \left \| H(u(t_i; \theta)) - y^{\text{target}}(t_i) \right \|_2^2.
        \label{eq:loss-state-observation}
    \end{equation}
    We can also consider the continuous evaluated loss function of the form
    \begin{equation}
         L(u(\cdot; \theta)) = \int_{t_0}^{t_1} h( u(t;\theta), \theta)  dt, 
         \label{eq:integrated-loss-function}
    \end{equation}
    with $h$ being a function that quantifies the contribution of the error term at every time $t \in [t_0, t_1]$. 
    Defining a loss function where just the empirical error is penalized is known as trajectory matching \cite{ramsay2017dynamic}. 
    Other methods like gradient matching and generalized smoothing the loss depends on smooth approximations of the trajectory and their derivatives. 
    \item[$ \blacktriangleright$] \textbf{The likelihood function or posterior probability.} From a statistical and physical perspective, it is common to assume that observations correspond to noisy observations of the underlying dynamical system, $y_i = H(u(t_i; \theta)) + \varepsilon_i$, with $\varepsilon_i$ errors or residual that are independent of each other and of the trajectory $u(\cdot ; \theta)$ \cite{ramsay2017dynamic}.
    When $H$ is the identity, each $y_i$ corresponds to the noisy observation of the state $u(t_i; \theta)$.
    If $p(Y | t , \theta)$ is the probability distribution of $Y=(y_1, y_2, \ldots, y_N)$, 
    the maximum likelihood estimator (MLE) of $\theta$ is defined as 
    \begin{equation}
        \theta^\text{MLE} 
        = 
        \argmax{\theta} \,\, \ell (Y | \theta) 
        = 
        \prod_{i=1}^N p(y_i | \theta, t_i) .
    \end{equation}
    When $\varepsilon_i \sim N(0, \sigma_i^2 \I)$ is the isotropic multivariate normal distribution, the maximum likelihood principle is the same as minimizing $- \log \ell(Y | \theta)$ which coincides with the mean squared error of Equation \eqref{eq:loss-state-observation} \cite{hastie2009elements},
    \begin{equation}
        \theta^\text{MLE} 
        = 
        \argmin{\theta} \, \left \{ - \log \ell (Y | \theta) \right \}
        = 
        \argmin{\theta} \, \sum_{i=1}^N 
        \, \frac{1}{2\sigma_i^2} \,
        \left \| y_i - H(u(t_i; \theta)) \right \|_2^2 .
        \label{eq:MLE}
    \end{equation}
    A Bayesian formulation of equation \eqref{eq:MLE} would consist in deriving a point estimate $\theta^\text{MLE}$, the posterior mean of the maximum a posteriori (MAP), based on the posterior distribution for $\theta$ following Bayes theorem as $p(\theta | Y) = {p(Y | \theta) \, p (\theta)}/{p(Y)}$, where $p(\theta)$ is the prior distribution \cite{pml1Book}.
    In most realistic applications, the posterior distribution is approximated using sampling algorithms such as Markov chain Monte Carlo (MCMC) \cite{gelman2013bayesian} or Sequential Monte Carlo (SMC) \cite{delmoralSequentialMonteCarlo2006}. 
    However, being able to compute gradients of the joint distribution allows for the design of more efficient inference algorithms, such as Hamiltonian Monte Carlo \cite{neal2011mcmc, Betancourt_2017}, integrated nested Laplace approximations \cite{Rue_Martino_Chopin_2009}, or variational inference techniques \cite{blei2017variational}.
    \item[$ \blacktriangleright$] \textbf{A quantity of interest.} In some applications we are interested in quantifying how the solution of the differential equation changes as we vary the parameter values; or more generally when it returns the value of some variable that is a function of the solution of a differential equation. The latter corresponds to the case in design control theory, a popular approach in aerodynamics modelling where goals include maximizing the speed of an airplane or the lift of a wing given the solution of the flow equation for a given geometry profile \cite{Jameson_1988,Giles_Pierce_2000,Mohammadi:2004dg}. 
\end{itemize}
In the rest of the manuscript we will use letter $L$ to emphasize that in many cases this will be a loss function, but without loss of generality this includes the richer class of functions included in the previous examples. 


In the context of optimization, the goal is to find the parameter $\theta$ that is a minimizer of $L(\theta)$. 
There exists a broad literature of optimization methods based on gradients, including gradient descent and its many variants \cite{ruder2016overview-gradient-descent}.
Gradient-based methods tend to outperform gradient-free optimization schemes when $1 \ll p$, as they are not prone to the curse of dimensionality \cite{Schartau2017}. 
In the case of gradient descent, the parameter $\theta$ is updated based on the iterative procedure given by
\begin{equation}\label{eq:gradient-descent}
    \theta^{m+1} 
    = 
    \theta^m 
    - 
    \alpha_m 
    \frac{dL}{d\theta}(\theta^m),
\end{equation}
with $\alpha_m$ some choice of the stepsize and some initialization $\theta^0 \in \R^p$. 
A direct implementation of gradient descent following Equation \eqref{eq:gradient-descent} is prone to converge to a local minimum and slows down in a neighborhood of saddle points. 
To address these issues, variants of this scheme employing more advanced updating strategies have been proposed, including Newton-type methods \cite{second-order-optimization}, quasi-Newton methods, acceleration techniques \cite{JMLR:v22:20-207}, and natural gradient descent methods \cite{doi:10.1137/22M1477805}. 

\subsubsection{Sensitivity matrix}

In general, loss functions considered are of the form $L(\theta) = L(u(\cdot, \theta), \theta)$. 
Using the chain rule we can derive 
\begin{equation} 
 \frac{dL}{d\theta} = \frac{\partial L}{\partial u} \frac{\partial u}{\partial \theta} + \frac{\partial L}{\partial \theta}.
 \label{eq:dLdtheta_VJP}
\end{equation} 
The two partial derivatives involving the loss function on the right-hand side are usually easy to evaluate.
For example, for the loss function in Equation \eqref{eq:quadratic-loss-function}, these are simply given by 
\begin{equation}
    \frac{\partial L}{\partial u} = u - u^{\text{target}}(t_1)
    \qquad 
    \frac{\partial L}{\partial \theta} = 0.
    \label{eq:dLdu}
\end{equation}
In most applications, the empirical component of the loss function $L(\theta)$, that is, the part of the loss that is a function on the data, will depend on $\theta$ just through $u$, meaning $\frac{\partial L}{\partial \theta} = 0$. 
However, regularization terms added to the loss can directly depend on the parameter $\theta$, that is $\frac{\partial L}{\partial \theta} \neq 0$.
In both cases, the complicated term to compute is the matrix of derivatives $\frac{\partial u}{\partial \theta}$, usually referred to as the \textit{sensitivity} $s$, and represents how much the full solution $u$ varies as a function of the parameter $\theta$, 
\begin{equation}
 s 
 = 
 \frac{\partial u}{\partial \theta} 
 =
 \begin{bmatrix}
   \frac{\partial u_1}{\partial \theta_1} & \dots & \frac{\partial u_1}{\partial \theta_p} \\
   \vdots & \ddots & \vdots \\
   \frac{\partial u_n}{\partial \theta_1} & \dots & \frac{\partial u_n}{\partial \theta_p}
 \end{bmatrix}
 \in \mathbb R^{n \times p}.
 \label{eq:sensitivity-definition}
\end{equation}
The sensitivity $s = s(t; \theta)$ defined in Equation \eqref{eq:sensitivity-definition} is a \textit{Jacobian}, that is, a matrix of first derivatives of a vector-valued function. 

Notice that the product $s v$, with $v \in \R^p$, is the directional derivative of the function $u(t; \theta)$ in the direction $v$, that is 
\begin{equation}
    s v
    = 
    \frac{\partial u}{\partial \theta} v 
    = 
    \lim_{h \rightarrow 0} \frac{u(t; \theta + h v) - u(t; \theta)}{h},
    \label{eq:directional-derivative}
\end{equation}
which represents how much the function $u(t; \theta)$ changes when we perturb $\theta$ in the direction of $v$. 


\subsection{Finite differences}
\label{section:finite-diffences}
Finite differences are arguably the simplest scheme to obtain the derivative of a function. 
In the case of the function $L : \R^p \mapsto \R$, a first-order Taylor expansion yields to the following expression for the directional derivative
\begin{equation}
 \frac{dL}{d\theta_i} (\theta) = \frac{L(\theta + \varepsilon e_i ) - L(\theta)}{\varepsilon} + \mathcal O (\varepsilon),
 \label{eq:finite_diff}
\end{equation}
with $e_i$ the $i$-th canonical vector and $\varepsilon$ the stepsize. 
Even better, the centered difference scheme leads to
\begin{equation}
 \frac{dL}{d\theta_i} (\theta) 
 =
 \frac{L(\theta + \varepsilon e_i ) - L(\theta - \varepsilon e_i)}{2\varepsilon}
 + \mathcal O (\varepsilon^2).
 \label{eq:finite_diff2}
\end{equation}
While Equation \eqref{eq:finite_diff} gives the derivative to an error of magnitude $\mathcal O (\varepsilon)$, the centered differences schemes improves the accuracy to $\mathcal O (\varepsilon^2)$ \cite{ascher2008-numerical-methods}. 
Further finite difference stencils of higher order exist in the literature \cite{Fornberg1988}. 
 
Finite difference scheme are subject to a number of issues, related to the parameter vector dimension and rounding errors.
Firstly, calculating directional derivatives requires at least one extra function evaluation per parameter dimension.
For the centered differences approach in Equation \eqref{eq:finite_diff2}, this requires a total of $2p$ function evaluations which demands solving the DE each time for a new set of parameters.
Second, finite differences involve the subtraction of two closely valued numbers, which can lead to floating point cancellation errors when the step size $\varepsilon$ is small \cite{Goldberg_1991_floatingpoint}. 
While small values of $\varepsilon$ lead to cancellation errors, large values of the stepsize give inaccurate estimations of the derivative. 
Furthermore, numerical solutions of DEs have errors that are typically larger than machine precision, which leads to inaccurate estimations of the gradient when $\varepsilon$ is too small (see also Section \ref{section:software-finite-differences}).
Finding the optimal value of $\varepsilon$ that balances these two effects is sometimes known as the \textit{stepsize dilemma}, for which algorithms based on prior knowledge of the function to be differentiated or algorithms based on heuristic rules have been introduced \cite{mathur2012stepsize-finitediff, BARTON_1992_finite_diff, SUNDIALS-hindmarsh2005sundials}. 

Despite these caveats, finite differences can prove useful in specific contexts, such as computing Jacobian-vector products (JVPs). 
Given a Jacobian matrix $J = \frac{\partial f}{\partial u}$ (or the sensitivity $s = \frac{\partial u}{\partial \theta}$) and a vector $v$, the product $Jv$ corresponding to the directional derivative and can be approximated as 
\begin{equation}
    Jv \approx \frac{f(u + \varepsilon v, \theta, t) - f(u, \theta, t)}{\varepsilon}.
\end{equation}
This approach is used in numerical solvers based on Krylov methods, where linear systems are solved by iteratively solving matrix-vectors products \cite{Ipsen_Meyer_1998}.
Furthermore, finite differences is commonly used as ground-truth in unit testing for more complex differentiation methods.

\subsection{Automatic differentiation}
\label{section:automatic-differentiation}
Automatic differentiation (AD) is a technique that generates new code representing derivatives of a given computer program defined by some evaluation procedure. 
Examples are code representing the tangent linear or adjoint operator of the original parent code. 
The names \textit{algorithmic} and \textit{computational} differentiation are also used in the literature, emphasizing the algorithmic rather than automatic nature of AD \cite{Griewank:2008kh, Naumann.2011, Margossian.2019}. 
The AD community has established a vibrant research landscape over more than three decades, as documented in a series of proceedings accompanying the International Conference on Algorithmic Differentiation held roughly every four years since 1991, with the earliest proceedings published in 1991 \cite{Griewank1991ADo} and the latest in 2024 covering a range of fundamental contributions on the algorithmic theory underlying AD (a full list of proceedings is available from \url{https://autodiff.org}). 
Acknowledging the extensive literature on AD, this section gives a condensed review, highlighting salient features that are relevant in the context of differentiable programming applied to DE-based models.

Any computer program implementing a given function can be reduced to a sequence of simple algebraic operations that have straightforward derivative expressions, based upon elementary rules of differentiation \cite{juedes1991taxonomy}.
The derivatives of the outputs of the computer program (dependent variables) with respect to their inputs (independent variables) are then combined using the chain rule.
One advantage of AD systems is their capacity to differentiate complex programs that include control flow, such as branching, loops or recursions. 

AD falls under the category of discrete methods.
Depending on whether the concatenation of the elementary derivatives is done as the program is executed (from input to output) or in a later instance where we trace-back the calculation from the end (from output to input), we refer to \textit{forward} or \textit{reverse} mode AD, respectively.
Neither forward nor reverse mode is more efficient in all cases \cite{Griewank_1989}, as we will discuss in Section \ref{sec:vjp-jvp}.

\subsubsection{Forward mode}

Forward mode AD can be implemented in different ways depending on the data structures we use when representing a computer program. 
Examples of these data structures include dual numbers and computational graphs \cite{Baydin_Pearlmutter_Radul_Siskind_2015}. 
These representations are mathematically equivalent and lead to the same implementation except for details in the compiler optimizations with respect to floating point ordering.

\paragraph{Dual numbers}
\label{section:dual-numbers}

Dual numbers extend the definition of a numerical variable that takes a certain value to also carry information about its derivative with respect to a certain parameter \cite{clifford1871dualnumbers}. 
We define a dual number based on two variables: a \textit{value} coordinate $x_1$ that carries the value of the variable and a \textit{derivative} (also known as partial or tangent) coordinate $x_2$ with the value of the derivative $\frac{\partial x_1}{\partial \theta}$. 
Just as complex numbers, we can represent dual numbers as an ordered pair $(x_1, x_2)$, sometimes known as Argand pair, or in the rectangular form 
\begin{equation}
 x_\epsilon = x_1 + \epsilon \, x_2,
\end{equation}
where $\epsilon$ is an abstract number called a perturbation or tangent, with the properties $\epsilon^2 = 0$ and $\epsilon \neq 0$.
This last representation is quite convenient since it naturally allows us to extend algebraic operations, like addition and multiplication, to dual numbers \cite{Karczmarczuk2001}. 
For example, given two dual numbers $x_\epsilon = x_1 + \epsilon x_2$ and $y_\epsilon = y_1 + \epsilon y_2$, it is easy to derive, using the fact $\epsilon^2=0$, that
\begin{equation}
 x_\epsilon + y_\epsilon = (x_1 + y_1) + \epsilon \, (x_2 + y_2)
 \qquad
 x_\epsilon y_\epsilon = x_1 y_1 + \epsilon \, (x_1 y_2 + x_2 y_1) .
\end{equation}
From these last examples, we can see that the derivative component of the dual number carries the information of the derivatives when combining operations (e.g., when the dual variables $x_2$ and $y_2$ carry the value of the derivative of $x_1$ and $x_2$ with respect to a parameter $\theta$, respectively). 

Intuitively, we can think of $\epsilon$ as being a differential in the Taylor series expansion, as evident in how the output of any scalar functions is extended to a dual number output:
\begin{align}
\begin{split}
    f(x_1 + \epsilon x_2)
    &= 
    f(x_1)
    + 
    \epsilon \, x_2 \,  f'(x_1)
    + 
    \epsilon^2 \cdot ( \ldots )\\
    &= 
    f(x_1)
    + 
    \epsilon \, x_2 \,  f'(x_1).
\end{split}
\label{eq:dual-number-function}
\end{align}
When computing first order derivatives, we can ignore everything of order $\epsilon^2$ or larger, which is represented in the condition $\epsilon^2 = 0$.
This implies that we can use dual numbers to implement forward AD through a numerical algorithm. 
In Section \ref{section:software-Forward-AD} we will explore how this is implemented. 
These ideas generalize to higher-order derivatives using the Taylor expansion of a variable to propagate derivatives, a method known as Taylor AD \cite{Griewank_Juedes_Utke_1996, Griewank:2008kh}.

Multidimensional dual numbers generalize dual numbers to include a different dual variable $\epsilon_i$ for each variable we want to differentiate with respect to \cite{Neuenhofen_2018, RevelsLubinPapamarkou2016}.
A multidimensional dual number is then defined as $x_\epsilon = x + \sum_{i=1}^p x_i \epsilon_i$, with the property that $\epsilon_i \epsilon_j = 0$ for all pairs $i$ and $j$.
Another extension of dual numbers that should not be confused with multidimensional dual numbers is hyper-dual numbers, which allow the computation of higher-order derivatives of a function \cite{fike2013multi}.

\paragraph{Computational graph}

A useful way of representing a computer program is via a computational graph with intermediate variables that relate the input and output variables. 
Most scalar functions of interest can be represented as a directed acyclic graph (DAG) with nodes associated to variables and edges to atomic operations \cite{Griewank:2008kh, Griewank_1989}, known as Kantorovich graph \cite{kantorovich1957mathematical} or its linearized representation via a Wengert trace/tape \cite{Wengert_1964, Griewank:2008kh}. 
We can define $v_{-p+1}, v_{-p+2}, \ldots, v_0 = \theta_1, \theta_2, \ldots, \theta_p$ the input set of variables; $v_{1}, \ldots, v_{m-1}$ the set of all the intermediate variables; and $v_m = L(\theta)$ the final output of a computer program. 
This can be done in such a way that the order is strict, meaning that each variable $v_i$ is computed just as a function of the previous variables $v_j$ with $j < i$. 
Once the graph is constructed, we can compute the derivative of every node with respect to the other, a quantity known as the tangent, using the Bauer formula \cite{Bauer_1974}:
\begin{equation}
    \frac{\partial v_j}{\partial v_i}
    = 
    \sum_{\substack{ \text{paths }w_0 \rightarrow w_1 \rightarrow \ldots \rightarrow w_K \\
                    \text{with } w_0=v_i, w_K = v_j}}
    \prod_{k=0}^{K-1} \frac{\partial w_{k+1}}{\partial w_{k}},
    \label{eq:bauer}
\end{equation}
where the sum is calculated with respect to all the directed paths in the graph connecting the input and target node.
Instead of evaluating the last expression for all possible paths, a simplification is to increasingly evaluate $j=-p+1, -p+1, \ldots, m$ using the recursion 
\begin{equation}
    \frac{\partial v_j}{\partial v_i}
    = 
    \sum_\text{$w$\text{ such that} $w \rightarrow v_j$}
    \frac{\partial v_j}{\partial w}
    \frac{\partial w}{\partial v_i}. 
    \label{eq:AD-graph-recursion}
\end{equation}
Since every variable node $w$ such that $w \rightarrow v_j$ is an edge of the computational graph has an index less than $j$, we can iterate this procedure as we run the computer program and solve for both the function and its derivative.
This is possible because in forward mode the term $\frac{\partial w}{\partial v_i}$ has been computed in a previous iteration, while $\frac{\partial v_j}{\partial w}$ can be evaluated at the same time the node $v_j$ is computed based on only the value of the parent variable nodes. 
The only requirement for differentiation is being able to compute the derivative/tangent of each edge/primitive and combine these using the recursion defined in Equation \eqref{eq:AD-graph-recursion}.

\subsubsection{Reverse mode}

Reverse mode AD is also known as the adjoint, or cotangent linear mode, or backpropagation in the field of machine learning. 
The reverse mode of AD has been introduced in different contexts \cite{griewank2012invented} and materializes the observation made by Phil Wolfe that if the chain rule is implemented in reverse mode, then the ratio between the computational cost of the gradient of a function and the function itself can be bounded by a constant that does not depend on the number of parameters to differentiate \cite{Griewank_1989, Wolfe_1982}, a point known as the \textit{cheap gradient principle} \cite{griewank2012invented}.  
Given a DAG of operations defined by a Wengert list, we can compute gradients of any given function in the same fashion as Equation \eqref{eq:AD-graph-recursion} but in decreasing order $j=m, m-1, \ldots, -p+1$ as
\begin{equation}
    \bar v_i 
    = 
    \frac{\partial \ell}{\partial v_i}
    = 
    \sum_\text{$w$\text{ such that} $v_i \rightarrow w$}
    \bar{w} \, \frac{\partial w}{\partial v_i}.
    \label{eq:reverse-mode-ad-definition}
\end{equation}
In this context, the notation $\bar{w} = \frac{\partial L}{\partial w}$ is introduced to signify the partial derivative of the output variable, here associated to the loss function, with respect to input and intermediate variables. 
This derivative is often referred to as the adjoint, dual, or cotangent, and its connection with the discrete adjoint method will be made more explicitly in Section \ref{section:comparison-discrete-adjoint-AD}. 

Since in reverse-mode AD the values of $\bar w$ are being updated in reverse order, in general
we need to know the state value of all the argument variables $v$ of $w$ in order to evaluate the terms $\frac{\partial w}{\partial v}$.
These state values (required variables) need to be either stored in memory during the evaluation of the function or recomputed on the fly in order to be able to evaluate the derivative. 
Checkpointing schemes exist to limit and balance the amount of storing versus recomputation (see section \ref{section:checkpointing}).

\subsubsection{AD connection with JVPs and VJPs}
\label{sec:vjp-jvp}
Forward and reverse AD is based on the sequential evaluation of Jacobian-vector products (JVPs) and vector-Jacobian products (VJPs), respectively. 
Let us consider for example the case of a loss function $L : \mathbb R^p \mapsto \mathbb R$ taking a total of $p$ arguments as inputs that is computed using the evaluation procedure $L(\theta) = \ell \circ g_{k} \circ \ldots \circ g_2 \circ g_1(\theta)$, with $\ell : \mathbb R^{d_k} \mapsto \mathbb R$ the final evaluation of the loss function after we apply in order a sequence of intermediate functions $g_i : \mathbb R^{d_{i-1}} \mapsto \mathbb R^{d_i}$, where we define $d_0 = p$ for simplicity. 
If we perturb the parameter $\theta \rightarrow \theta + \delta \theta$, this will produce a perturbation $L (\theta) \rightarrow L(\theta) + \delta L$ in the loss function that can be computed at first order in $\delta \theta$ using the chain rule as: 
\begin{equation}
     \delta L = \nabla_\theta L \cdot \delta \theta = \nabla \ell \cdot Dg_{k} \cdot Dg_{k-1} \cdot \ldots \cdot Dg_2 \cdot Dg_1 \cdot \delta \theta , 
    \label{eq:deltaL}
\end{equation}
with $Dg_i$ the Jacobian of each intermediate function evaluated at the intermediate values $g_{i-1} \circ g_{i-2} \circ \ldots \circ g_i (\theta)$ \cite{Giering_Kaminski_1998}.

\begin{figure}[p]
    \centering
    \includegraphics[width=0.85\textwidth]{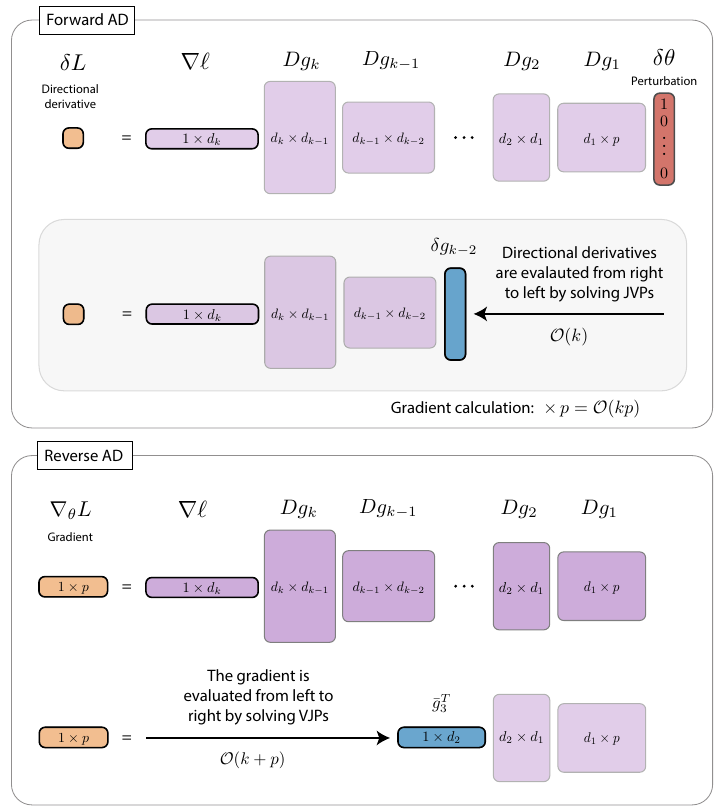}
    \caption{Comparison between forward and reverse AD. Changing the order of Jacobian and vector multiplications changes the total number of floating-point operations, which leads to different computational complexities between forward and reverse mode. When computing directional derivatives with forward AD, there is a total of $\mathcal O (k)$ JVPs that need to be computed, which considering we need to repeat this procedure $p$ times gives a total complexity of $\mathcal O (kp)$. This is the opposite of what happens when we carry the VJPs from the left side of the expression, where the matrix of size $d_1 \times p$ has no effect in the intermediate calculations, making all the intermediate calculations $\mathcal O (1)$ with respect to $p$ and a total complexity of $\mathcal O (k + p)$. }
    \label{fig:vjp-jvp}
\end{figure}

In forward AD, we can compute $\delta L$ from Equation \eqref{eq:deltaL} by defining the intermediate perturbation $\delta g_j$ as the sequential evaluation of the JVP given by the map between tangent spaces ${\delta x \mapsto Dg_j (x) \cdot \delta x}$ \cite{Griewank:2008kh}:
\begin{align}
    \delta g_0 &= \delta \theta \\
    \delta g_j &= D g_j \cdot \delta g_{j-1} \qquad j = 1, 2, \ldots, k \\
    \delta L &= \nabla \ell \cdot \delta g_{k}.
\end{align}
For $\| \delta \theta \|_2 = 1$, this procedure will return $\delta L$ as the value of the directional derivative of $L$ evaluated at $\theta$ in the direction $\delta \theta$ (see Equation \eqref{eq:directional-derivative}). 
In order to compute the full gradient $\nabla L \in \R^p$, we need to perform this operation $O(p)$ times, which requires a total of $p \, (d_2 d_1 + d_3 d_2 + \ldots + d_k d_{k-1} + d_k )= \mathcal O (kp)$ operations.

In the case of reverse AD, we observe that $\nabla \ell \in \mathbb R^{d_k}$ is a vector, so we can instead compute $\delta L$ for all possible perturbations $\delta \theta$ by solving the multiplication involved in Equation \eqref{eq:deltaL} starting from the left-hand side. 
This is carried by the sequential definition of intermediate variables $\bar g_j$ computed as VJPs that map between co-tangent (or normal spaces) $\bar y \mapsto \bar y^T \cdot Dg_j$:
\begin{align}
    \bar g_{k}^T &= \nabla \ell \\
    \bar g_{j-1}^T &= \bar g_{j}^T \cdot Dg_j \qquad j = k, k-1, \ldots, 1 \\
    \nabla L &= \bar g_0.
\end{align}
Since this procedure needs to be evaluated just once to evaluate $\nabla L$, we conclude that reverse AD requires a total of $ d_k d_{k-1} + d_{k-1} d_{k-2} + \ldots + d_2 d_1 + d_1 p = \mathcal O (k+p)$ operations. 

The reverse mode will in general be faster when $1 \ll p$. 
This example is illustrated in Figure \ref{fig:vjp-jvp}. 
In the general case of a function $L : \R^p \mapsto \R^q$ with multiple outputs and a total of $k$ intermediate functions, the cost of forward AD is $\mathcal O (pk + q)$ and the cost of reverse is $\mathcal O (p + kq)$.
When the function to differentiate has a larger input space than output ($q \ll p$), AD in reverse mode is more efficient as it propagates the chain rule by computing VJPs.
For this reason, reverse AD is often preferred in both modern machine learning and inverse methods.
However, notice that reverse mode AD requires saving intermediate variables through the forward run in order to run backwards afterwards \cite{Bennett_1973}, leading to performance overhead that makes forward AD more efficient when $p \lesssim q$ \cite{Griewank_1989, Margossian.2019, Baydin_Pearlmutter_Radul_Siskind_2015}. 

In practice, most AD systems are reduced to the computation of only directional derivatives (JVPs) or gradients (VJPs) \cite{Griewank:2008kh}.
Full Jacobians $J \in \R^{n \times p}$ (e.g., the sensitivity $s = \frac{\partial u}{\partial \theta} \in \R^{n \times p}$) can be fully reconstructed by the independent computation of the $p$ columns of $J$ via the JVPs $J e_i$, with $e_i \in \R^p$ the canonical vectors; or by the calculation of the $m$ rows of $J$ via the VJPs $e_j^T J$, with $e_j \in \R^n$.
In other words, forward AD computes Jacobians column-by-column while reverse AD does it row-by-row.  

\subsubsection{Further remarks}
\label{sec:ad-further}
In this section, we will briefly comment on some further implementation aspects of AD that are of particular importance when working with differential equations and numerical solvers. 

~\\
\noindent
\textit{Sparsity patterns.} 
The sparsity structure of the Jacobian can be exploited with the combination of forward and reverse AD.
When the sparsity pattern is known, colored AD efficiently chunk the calculation of the Jacobian into multiple JVPs or VJPs \cite{gebremedhin2005color}. 
This results in a smaller number of evaluations of JVPs/VJPs compared to the one required to compute all entries of a dense Jacobian (e.g., \textcite{pal2024nonlinearsolve}). 
An example of this is given by the arrowhead matrix $J_{\text{arrowhead}} \in \R^{n \times n}$ defined as:
\begin{equation}
    J_{\text{arrowhead}} = \begin{bmatrix}
        \bullet & \bullet & \bullet & \bullet & \cdots & \bullet  & \bullet   \\
        \bullet & \bullet & 0  &  0    &   & 0        & 0        \\
        \bullet & 0 & \bullet & 0 & & 0        & 0        \\
        \bullet & 0 & 0 & \bullet & & 0 & 0 \\
        \vdots & & & & \ddots & &  \\
        \bullet & 0 & 0 & 0   &   & \bullet   & 0        \\
        \bullet & 0 & 0        & 0    &    & 0        & \bullet
    \end{bmatrix},
\end{equation}
where $\bullet$ indicate non trivial zero entries of the Jacobian.
In this case, both forward and reverse AD will have to perform $n$ VJPs and JVPs, respectively, and there is no computational benefit of using colored AD. 
Instead, a combination of forward and reverse AD can be used to color the Jacobian with two forward and one reverse AD evaluation (namely, the two JVPs given by $J_{\text{arrowhead}} e_1$ and $J_{\text{arrowhead}} 1_n$, with $1_n \in \R^n$ the vector with all ones, and one JVP given by $e_1^T J_{\text{arrowhead}}$ are enough to solve for all the non-zero entries of the Jacobian). 
It is important to highlight that sparsity patters may play an important role in differentiation of the solver even in cases with dense sensitivities or Jacobians. 
This is due to the fact that intermediate Jacobians can still include known sparsity patterns that can be exploded (e.g., in solving the nonlinear problem involved at each different step of the solver).

~\\
\noindent
\textit{Beyond forward and reverse.} 
Forward and reverse AD are not the only ways of computing derivatives using the chain rule, and instead they both correspond to particular choices of how to evaluate the Bauer formula (Equation \eqref{eq:bauer}) \cite{Griewank:2008kh}. 
Efficient Jacobian accumulation refers to the problem of finding the optimal way of evaluating the Bauer formula so that it minimizes the total number of operations \cite{naumann2000optimized}. 
It can be shown that such a problem is NP-complete \cite{Naumann_2008}, and different heuristics have been introduced in the literature \cite{Naumann2023}. 

~\\
\noindent
\textit{AD of non-smooth function.} 
It is often unavoidable (or the result of dealing with legacy numerical schemes) that computational models contain non-smooth or non-differentiable functions. 
Recognizing and dealing with such cases in the context of AD goes back at least as far as \textcite{griewank1995automatic} who considered functions such as absolute value, minimum, maximum, and Euclidean norms. 
General approaches have been the introduction of generalized derivatives, successively piecewise linearizations, linear interpolating functions, sorting functions, or threshold functions (e.g, \textcite{clarke2008nonsmooth}, \textcite{khan2014generalized}, \textcite{fiege2018algorithmic}, \textcite{bethke2024semismooth}, \textcite{Griewank_2013, Griewank_Rojas_2019}, \textcite{Bolte_Boustany_Pauwels_Pesquet-Popescu_2022}). 
In the context of machine learning, such non-smooth functional forms may appear especially in activation functions. A
separate, fast-evolving ML literature on non-differentiable functions exists, a comprehensive review of which is well beyond the scope of this paper (but see, e.g., \textcite{krishna2024gd} for a recent discussion). 

~\\
\noindent
\textit{Parallel AD.} 
In order to harness high performance computing, many large-scale numerical simulations rely on parallel algorithms. 
Message Passing Interface (MPI) has established itself as a ubiquitous standard library that underlies most scalable application codes. 
Differentiating such codes has been considered early on for forward-mode AD (e.g., \textcite{hovland1998automatic}). 
The added complexity of direct differentiation in reverse of MPI libraries, such as reversal of \texttt{waitall} calls, has been discussed in in the context of operator overloading \cite{Bischof.2008} and source transformation \cite{Utke.2009}. 
Extensive research has since been conducted on this subject, including recent work by \textcite{Towara.2015} who present an adjoint MPI library in the context of a CFD code,
\textcite{Hueckelheim.2022} on AD of parallel loops in OpenMP, and \textcite{Moses.2022f0d} on compiler augmentations to support scalable AD-generated reverse-mode simulations.

\subsection{Complex step differentiation}
\label{section:comple-step-differentiation}
An alternative to finite differences that avoids subtractive cancellation errors is based on complex variable analysis. 
The first proposals originated in 1967 using the Cauchy integral theorem involving the numerical evaluation of a complex-valued integral \cite{Lyness_1967, Lyness_Moler_1967}.
A newer approach recently emerged that uses the complex generalization of a real function to evaluate its derivatives \cite{Squire_Trapp_1998_complex_diff, Martins_Sturdza_Alonso_2003_complex_differentiation}. 
Assuming that the function $L(\theta)$ admits a holomorphic extension (that is, it can be extended to a complex-valued function that is analytical and differentiable \cite{stein2010complex}), the Cauchy-Riemann conditions can be used to evaluate the derivative with respect to one single scalar parameter $\theta \in \R$ as
\begin{equation}
    \frac{dL}{d\theta} = \lim_{\varepsilon \rightarrow 0} \frac{\text{Im}(L(\theta + i \varepsilon))}{\varepsilon},
\end{equation}
where $i$ is the imaginary unit satisfying $i^2 = -1$. 
The order of this approximation can be found using the Taylor expansion:
\begin{equation}
    L(\theta + i \varepsilon)
    = 
    L(\theta) + i \varepsilon \frac{dL}{d\theta} 
    - 
    \frac 1 2  \varepsilon^2
    \frac{d^2 L}{d\theta^2}
    + 
    \mathcal O (\varepsilon^3).
\end{equation}
Computing the imaginary part $\text{Im}(L(\theta + i \varepsilon))$ leads to
\begin{equation}
    \frac{dL}{d\theta} 
    = 
    \frac{\text{Im}(L(\theta + i \varepsilon))}{\varepsilon}
    + 
    \mathcal{O} (\varepsilon^2).
    \label{eq:complex-step-definition}
\end{equation}
The method of \textit{complex step differentiation} consists then in estimating the gradient as $\text{Im}(L(\theta + i \varepsilon)) / \varepsilon$ for a small value of $\varepsilon$. 
Besides the advantage of being a method with precision $\mathcal{O}(\varepsilon^2)$, the complex step method avoids subtracting cancellation error and then the value of $\varepsilon$ can be reduced to almost machine precision error without affecting the calculation of the derivative. 
However, a major limitation of this method is that it only applicable to locally complex analytical functions \cite{Martins_Sturdza_Alonso_2003_complex_differentiation} and does not outperform AD (see Sections \ref{section:direct-methods} and \ref{section:recomendations}). 
One additional limitation is that it requires the evaluation of mathematical functions with small complex values, e.g., operations such as $\sin(1 + 10^{-16} i)$, which are not necessarily always computable to high accuracy with modern math libraries.
Extension to higher order derivatives can be obtained by introducing multicomplex variables \cite{Lantoine_Russell_Dargent_2012}.

\subsection{Symbolic differentiation}
\label{section:symbolic}
In symbolic differentiation, functions are represented algebraically instead of algorithmically, which is why many symbolic differentiation tools are included inside computer algebra systems (CAS) \cite{Symbolics_jl_2022}. 
Instead of numerically evaluating the final value of a derivative, symbolic systems assign variable names, expressions, operations, and literals to \textit{algebraic} objects. 
For example, the relation $y = x^2$ is interpreted as an expression with two variables, $x$ and $y$, and the symbolic system generates the derivative $y' = 2 \times x$ with $2$ a numeric literal, $\times$ a binary operation, and $x$ the same variable assignment as in the original expression.

The general issue with symbolic differentiation is \textit{expression swell}, i.e. the fact that the size of a derivative expression can be much larger than the original expression \cite{Baydin_Pearlmutter_Radul_Siskind_2015}. 
One way to visualize this swell is to note that the product rule grows an expression of $f(x)g(x)$ into two expressions, namely $\frac{d}{dx}(f(x)g(x)) = \frac{df}{dx}g(x) + f(x)\frac{dg}{dx}$, and thus the composition of many functions leads to a large derivative expression. 
AD avoids expression swell by instead numerically calculating the derivative of a given expression at some fixed value, never representing the general derivative but only at the values obtained by the forward pass. 
This eager evaluation of the derivative around a given value forces the intermediate computation into the JVPs or VJPs form as a way to continually pass forward/reverse the current state. 
Meanwhile, symbolic differentiation can represent the complete derivative expression and thus avoid being forced into a given computation order, but at the memory cost of having to represent larger expressions.

However, it is important to acknowledge the close relationship between AD and symbolic differentiation.
AD uses symbolic differentiation in its definition of primitives which are then chained together in a specific way to form VJPs and vector products. 
Forward AD can be expressed as a form of symbolic differentiation with a specific choice of common subexpression elimination, i.e. forward AD can be expressed as a symbolic differentiation with a specific choice of how to accumulate the intermediate calculations so that expression growth can be avoided \cite{juedes1991taxonomy, Elliott_2018, Laue2020, Dürrbaum_Klier_Hahn_2002}.
However, general symbolic differentiation can have many other choices for the differentiation order, and does not in general require computation using the JVPs or VJPs \cite{Baydin_Pearlmutter_Radul_Siskind_2015}. 
This is apparent for example when computing sparse Jacobians, where generally symbolic differentiation computes entries element-by-element while forward AD computes the matrix column-by-column and reverse AD computes row-by-row (see Section \ref{sec:ad-further}).

\subsection{Forward sensitivity equations}
\label{section:sensitivity-equation}
An easy way to derive an expression for the sensitivity $s$ defined in Equation \eqref{eq:sensitivity-definition} is by deriving the forward sensitivity equations \cite{ramsay2017dynamic}, a method also referred to as continuous local sensitivity analysis (CSA). 
If we consider the original ODE given by Equation \eqref{eq:original_ODE} and we differentiate with respect to $\theta$, we then obtain
\begin{equation}
    \frac{d}{d\theta} \left( \frac{du}{dt}  - f(u(t; \theta), \theta, t) \right) = 0.
\end{equation}
Assuming that a unique solution exists and both $\frac{\partial f}{\partial u}$ and $\frac{\partial f}{\partial \theta}$ are continuous in the neighborhood of the solution, or under the guarantee of interchangeability of the derivatives \cite{gronwall1919note}, for example by assuming that both $\frac{du}{dt}$ and $\frac{du}{d\theta}$ are differentiable \cite{math8111947}, we can derive
\begin{equation}
 \frac{d}{d\theta} \frac{du}{dt} 
 =
 \frac{d}{d\theta} f(u(t; \theta), \theta, t)
 = 
 \frac{\partial f}{\partial \theta}
 + 
 \frac{\partial f}{\partial u} \frac{\partial u}{\partial \theta}.
\end{equation}
Identifying the sensitivity matrix $s(t)$ now as a function of time, we obtain the \textit{sensitivity differential equation} 
\begin{equation}
 \frac{ds}{dt} = \frac{\partial f}{\partial u} s + \frac{\partial f}{\partial \theta}.
 \label{eq:sensitivity_equations}
\end{equation}
The initial condition is simply given by $s(t_0) = \frac{du_0}{d\theta}$, which is zero unless the initial condition explicitly depends on the parameter $\theta$.
Both the original ODE of size $n$ and the forward sensitivity equation of size $np$ are solved simultaneously, which is necessary since the forward sensitivity DE directly depends on the value of $u(t; \theta)$.  
This implies that as we solve the ODE, we can ensure the same level of numerical precision for the two of them inside the numerical solver.

In contrast to the methods previously introduced, the forward sensitivity equations find the derivative by solving a new set of continuous differential equations.
Notice also that the obtained sensitivity $s(t)$ can be evaluated at any given time $t$. 
This method can be labeled as forward, since we solve both $u(t; \theta)$ and $s(t)$ as we solve the DE forward in time, without the need of backtracking any operation though the solver.
By solving the forward sensitivity equation and the original ODE for $u(t; \theta)$ simultaneously, we ensure that by the end of the forward step we have calculated both $u(t; \theta)$ and $s(t)$. 


\subsection{Discrete adjoint method}
\label{section:discrete-adjoint}
Also known as the adjoint state method, the discrete adjoint method is another example of a discrete method that aims to find the gradient by solving an alternative system of linear equations, known as the \textit{adjoint equations}, simultaneously with the original system of equations defined by the numerical solver. 
These methods are extremely popular in optimal control theory in fluid dynamics, for example for the design of geometries for vehicles and airplanes that optimize performance \cite{Elliott_Peraire_1996, Giles_Pierce_2000} or in ocean state estimation \cite{Wunsch.2008}.

The idea of the adjoint method is to treat the DE as a constraint in an optimization problem and then differentiate an objective function subject to that constraint. 
Mathematically speaking, this can be treated both from a duality or Lagrangian perspective \cite{Giles_Pierce_2000}.
In agreement with other authors, we prefer to derive the equation using the former as it gives better insights to how the method works and allows generalization to other user cases \cite{Givoli_2021}.
We will introduce the adjoint method for computing first-order derivatives (gradient), but generalizations to higher-order derivatives (e.g., Hessian) has also been explored in the literature (e.g., \textcite{Papadimitriou_Giannakoglou_2008}, \textcite{pacaud2022batched}).

\subsubsection{Adjoint state equations}

The derivation of the discrete adjoint equations is carried out once the numerical scheme for solving Equation \eqref{eq:original_ODE} has been specified.  
Given a discrete sequence of timesteps $t_0, t_1, \ldots, t_M$, we aim to find approximate numerical solutions $u^m \approx u(t_m; \theta)$. 
Any numerical solver, including the ones discussed in Section \ref{section:intro-numerical-solvers}, can be understood as solving the (in general nonlinear) system of equations defined by $G(U; \theta) = 0$, where $U$ is the super-vector $U = (u_1, u_2, \ldots, u_M) \in \R^{nM}$, and we have combined the system of equations defined by the iterative solver as $G(U; \theta) = (g_1(u^1; \theta), \ldots, g_M(u^M; \theta)) = 0$ (see Equation \eqref{eq:solver-constriant-example}).

\begin{figure}[t]
    \centering
    \includegraphics[width=0.95\textwidth]{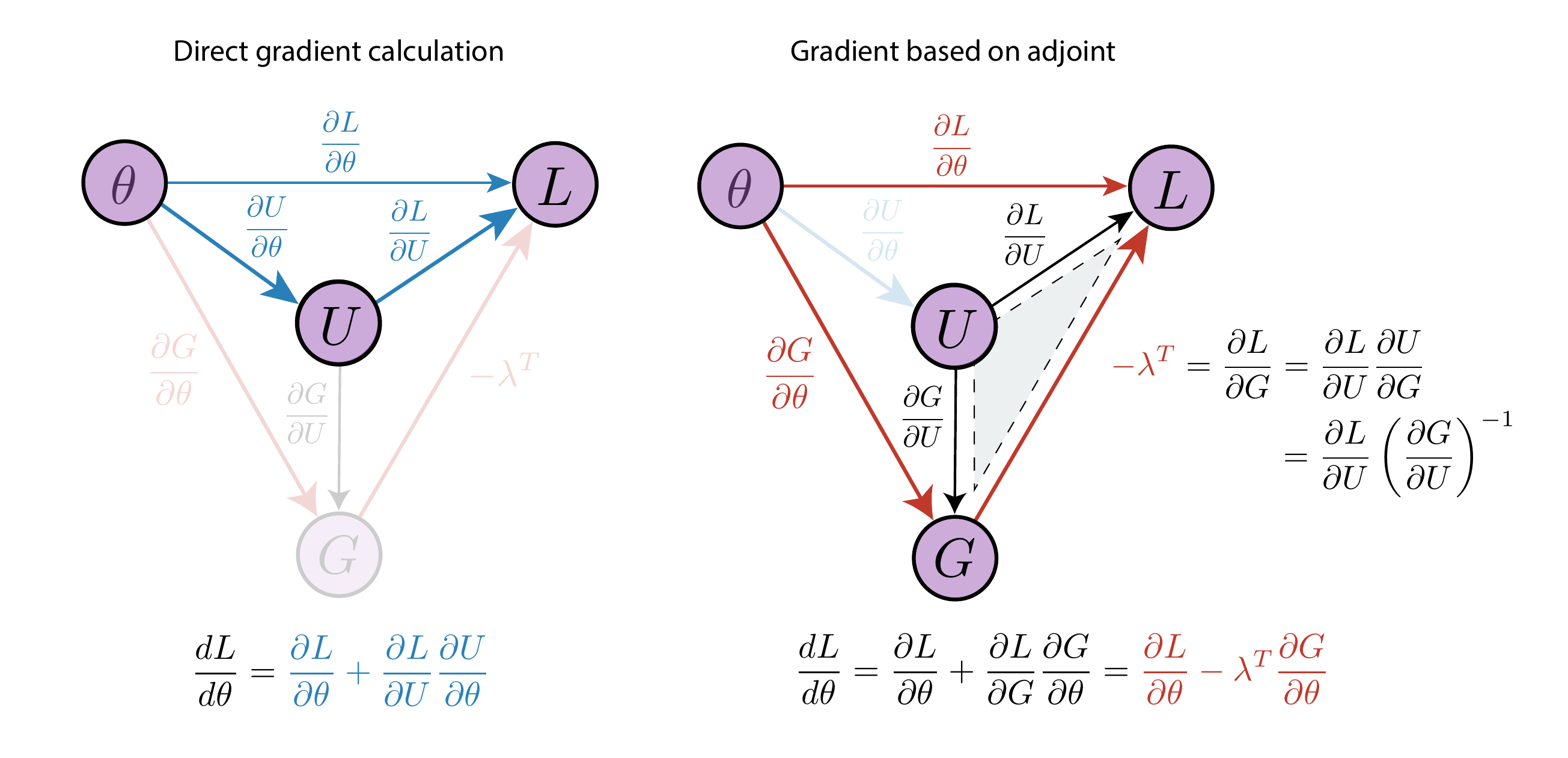}
    \caption{Diagram showing how gradients are computed using discrete adjoints. On the left, we see how gradients will be computed if we use the chain rule applied to the directed triangle defined by the variables $\theta$, $U$, and $L$ (blue arrows). However, we can define the intermediate vector variable $G = G(U; \theta)$, which satisfies $G  = 0$ as long as the discrete system of differential equations are satisfied, and apply the chain rule instead to the triangle defined by $\theta$, $G$, and $L$ (red arrows). In the red diagram, the calculation of $\frac{\partial L}{\partial G}$ is done by pivoting in $U$ as shown in the right diagram (shaded area). Notice that the use of adjoints avoids the calculation of the sensitivity $\frac{\partial U}{\partial \theta}$. The adjoint is defined as the partial derivative $\lambda^T = - \frac{\partial L}{\partial G}$ representing changes in the loss function due to variations in the discrete equation $G(U; \theta) = 0$. 
    }
    \label{fig:discrete-adjoint}
\end{figure}

We are interested in differentiating an objective or loss function $L(\theta) = L(U(\theta), \theta)$ with respect to the parameter $\theta$. 
Since here $U$ is the discrete set of evaluations of the solver, examples of loss functions now include 
\begin{equation}
    L(U, \theta) 
    = 
    \frac{1}{2} \sum_{m=1}^M w_m \left \| u^m - u_m^\text{obs} \right \|_2^2, 
\end{equation}
with $u_m^\text{obs}$ the observed time-series, and $w_m \geq 0$ some arbitrary weights (potentially, many of them equal to zero). 
Similarly to Equation \eqref{eq:dLdtheta_VJP} we have 
\begin{equation}
    \frac{dL}{d\theta} 
    = 
    \frac{\partial L}{\partial \theta} 
    + 
    \frac{\partial L}{\partial U} \frac{\partial U}{\partial \theta}.
    \label{eq:dhdtheta0}
\end{equation}
By differentiating the constraint $G(U; \theta) = 0$, we obtain
\begin{equation}
    \frac{dG}{d\theta} 
    = 
    \frac{\partial G}{\partial \theta} 
    + 
    \frac{\partial G}{\partial U} \frac{\partial U}{\partial \theta}
    =
    0,
\end{equation}
which is equivalent to 
\begin{equation}
    \frac{\partial U}{\partial \theta} 
    = 
    - \left( \frac{\partial G}{\partial U} \right)^{-1} \frac{\partial G}{\partial \theta}.
    \label{eq:adjoint-inversion-implicit-theorem}
\end{equation}
Replacing this last expression into Equation \eqref{eq:dhdtheta0}, we obtain
\begin{equation}
    \frac{dL}{d\theta} 
    =
    \frac{\partial L}{\partial \theta} 
    - 
    \frac{\partial L}{\partial U}
    \left( \frac{\partial G}{\partial U} \right)^{-1} 
    \frac{\partial G}{\partial \theta}.
    \label{eq:dhdtheta}
\end{equation}
The important trick used in the discrete adjoint method is the rearrangement of the multiplicative terms involved in equation \eqref{eq:dhdtheta}. 
Computing the full Jacobian/sensitivity $\partial U / \partial \theta$ will be computationally expensive and involves the product of two matrices (Equation \eqref{eq:adjoint-inversion-implicit-theorem}). 
However, we are not interested in the calculation of the Jacobian, but instead in the VJP given by $\frac{\partial L}{\partial U} \frac{\partial U}{\partial \theta}$. 
By rearranging these terms and relying on the intermediate variable $G(U; \theta)$, we can make the same computation more efficient. 
This leads to the definition of the adjoint $\lambda \in \R^{nM}$ as the solution of the linear system of equations 
\begin{equation}
    \left( \frac{\partial G}{\partial U}\right)^T \lambda 
    =  
    \left( \frac{\partial L}{\partial U} \right)^T,
    \label{eq:adjoint-state-equation}
\end{equation}
or equivalently,
\begin{equation}
    \lambda^T = \frac{\partial L}{\partial U} \left( \frac{\partial G}{\partial U} \right)^{-1}.
    \label{eq:def_adjoint}
\end{equation}
Replacing Equation \eqref{eq:def_adjoint} into \eqref{eq:dhdtheta} yields
\begin{equation}
    \frac{dL}{d\theta} 
    =
    \frac{\partial L}{\partial \theta} 
    - 
    \lambda^T \frac{\partial G}{\partial \theta}.
    \label{eq:gradient-adjoint-state-method}
\end{equation}
These ideas are summarized in the diagram in Figure \ref{fig:discrete-adjoint}, where we can also see an interpretation of the adjoint as being equivalent to $\lambda^T = - \frac{\partial L}{\partial G}$. 

Notice that the algebraic equation of the adjoint $\lambda$ in Equation \eqref{eq:adjoint-state-equation} is a linear system of equations, even when the original system $G(U; \theta)=0$ is not necessarily linear in $U$.
This means that while solving the original ODE may require multiple iterations in order to solve the non-linear system $G(U; \theta) = 0$ (e.g., by using Krylov methods), the backwards step to compute the adjoint is one single linear system of equations. 

\subsubsection{Simple linear system}

To gain further intuition about the discrete adjoint method, let us consider the simple case of the explicit linear one-step numerical solver, where at every step we need to solve the equation $u^{m+1} = g_m (u^m; \theta) = A_m (\theta) \, u^m + b_m(\theta)$, where $A_m(\theta) \in \R^{n \times n}$ and $b_m(\theta) \in \R^n$ are defined by the numerical solver \cite{Johnson}. 
This condition can be written in a more compact manner as $G(U; \theta)=A(\theta) U - b(\theta) = 0$, that is 
\begin{equation}
    A(\theta) U 
    = 
    \begin{bmatrix}
        \I_{n} & 0 &   &  & \\
        -A_1 & \I_{n} & 0 &  &  \\
          & -A_2 & \I_{n} & 0 &  \\
         &  &   & \ddots &   \\
         &  &  & -A_{M-1} & \I_{n}
    \end{bmatrix}
    \begin{bmatrix}
        u^1 \\
        u^2 \\
        u^3 \\
        \vdots \\
        u^M
    \end{bmatrix}
    = 
    \begin{bmatrix}
        A_0 u_0 + b_0 \\
        b_1 \\
        b_2 \\
        \vdots \\
        b_{M-1}
    \end{bmatrix}
    = 
    b(\theta), 
\end{equation}
with $\I_{n}$ the identity matrix of size $n \times n$.
Notice that in most cases, the matrix $A(\theta)$ is quite large and mostly sparse. 
While this representation of the discrete differential equation is  convenient for mathematical manipulations, when solving the system we rely on iterative solvers that save memory and computation. 

For the linear system of discrete equations $G(U; \theta)=A(\theta) U - b(\theta)=0$, we have 
\begin{equation}
    \frac{\partial G}{\partial \theta} 
    = 
    \frac{\partial A }{\partial \theta} U - \frac{\partial b}{\partial \theta},
\end{equation}
so the desired gradient in Equation \eqref{eq:gradient-adjoint-state-method} can be computed as 
\begin{equation}
    \frac{dL}{d\theta} 
    = 
    \frac{\partial L}{\partial \theta} 
    - 
    \lambda^T \left( \frac{\partial A }{\partial \theta} U - \frac{\partial b}{\partial \theta} \right),
    \label{eq:dhdtheta_linear}
\end{equation}
with $\lambda$ the discrete adjoint obtained by solving the linear system in Equation \eqref{eq:adjoint-state-equation},
\begin{equation}
    A(\theta)^T \lambda 
    =
    \begin{bmatrix}
        \I_{n} & -A_1^T &   &  & \\
        0 & \I_{n} & -A_2^T &  &  \\
          & 0 & \I_{n} & -A_3^T &  \\
         &  &   & \ddots & -A_{M-1}^T  \\
         &  &  & 0 & \I_{n}
    \end{bmatrix}
    \begin{bmatrix}
        \lambda_1 \\
        \lambda_2 \\
        \lambda_3 \\
        \vdots \\
        \lambda_M
    \end{bmatrix}
    = 
    \begin{bmatrix}
        w_1 (u^1 - u_1^\text{obs}) \\
        w_2 (u^2 - u_2^\text{obs}) \\
        w_3 (u^3 - u_3^\text{obs}) \\
        \vdots \\
        w_M (u^M - u_M^\text{obs})     
    \end{bmatrix}
    = 
    \frac{\partial L}{\partial U}^T.
    \label{eq:linea-adjoint-state-equation}
\end{equation}
This is a linear system of equations with the same size as the original $A(\theta) U = b(\theta)$, but involving the adjoint matrix $A^T$. 
Computationally this also means that if we can solve the original system of discretized equations then we can also solve the adjoint at the same computational cost (e.g., by using the LU factorization of $A(\theta)$). 
Another more natural way of finding the adjoints $\lambda_m$ is by noticing that the system of equations \eqref{eq:linea-adjoint-state-equation} is equivalent to the final value problem 
\begin{equation}
    \lambda_{m} = A_{m}^T \lambda_{m+1} + w_m(u^m - u_m^\text{obs}),
    \label{eq:adjoint-discrete-linear-example}
\end{equation}
with final condition $\lambda_M$. 
This means that we can efficiently compute the adjoint equation in reverse mode, starting from the final state $\lambda_M$ and computing the values of $\lambda_m$ in decreasing index order. 
Unless the loss function $L$ is linear in $U$, this procedure requires knowledge of the value of $u^m$ (or some equivalent form of it) at any given timestep $t_m$. 

\subsection{Continuous adjoint method}
\label{section:continuous-adjoint}
The continuous adjoint method, also known as continuous adjoint sensitivity analysis (CASA), operates by defining a convenient set of new DEs for the adjoint variable and using this to compute the gradient in a more efficient manner. 
The continuous adjoint method follows the same logic as the discrete adjoint method, but where the discretization of the DE does not happen until the very last step, when the solutions need to be solved numerically. 

Consider an integrated loss function defined in Equation \eqref{eq:integrated-loss-function} of the form 
\begin{equation}
    L(u; \theta) = \inttime h(u(t;\theta), \theta) dt,
\end{equation}
and its derivative with respect to the parameter $\theta$ given by the following integral involving the sensitivity matrix $s(t)$:
\begin{equation}
    \frac{dL}{d\theta}
    = 
    \inttime \left( \frac{\partial h}{\partial \theta} + \frac{\partial h}{\partial u} s(t) \right) dt.
    \label{eq:casa-loss}
\end{equation}
Just as in the case of the discrete adjoint method, the complicated term to evaluate in the last expression is the sensitivity $s(t)$.
Again, the trick is to evaluate the VJP $\frac{\partial h}{\partial u} s(t)$ by first defining an intermediate adjoint variable. 
The continuous adjoint equation is now obtained by finding the dual/adjoint equation associated to the forward sensitivity equation using the weak formulation of Equation \eqref{eq:sensitivity_equations} \cite{brezis2011functional}. 
The adjoint equation is obtained by writing the forward sensitivity equation in the form 
\begin{equation}
    \inttime \lambda(t)^T \left( \frac{ds}{dt} - \frac{\partial f}{\partial u} \, s - \frac{\partial f}{\partial \theta}  \right) dt 
    = 
    0,
    \label{eq:integrated-sensitivity-equation}
\end{equation}
where this equation must be satisfied for every suitable function $\lambda: [t_0, t_1] \mapsto \R^n$ in order for Equation \eqref{eq:sensitivity_equations} to be true. 
The next step is to get rid of the time derivative applied to the sensitivity $s(t)$ using integration by parts: 
\begin{equation}
    \inttime \lambda(t)^T \frac{ds}{dt} dt
    = 
    \lambda(t_1)^T s(t_1) - \lambda(t_0)^T s(t_0)
    -
    \inttime \frac{d\lambda^T}{dt} s(t)\, dt.
\end{equation}
Replacing this last expression into Equation \eqref{eq:integrated-sensitivity-equation} we obtain 
\begin{equation}
    \inttime \left( - \frac{d\lambda^T}{dt} -  \lambda(t)^T \frac{\partial f}{\partial u} \right) s(t) dt
    =
    \inttime \lambda(t)^T \frac{\partial f}{\partial \theta} dt 
    - 
    \lambda(t_1)^T s(t_1)
    + 
    \lambda(t_0)^T s(t_0).
    \label{eq:casa-semiadjoint}
\end{equation}
At first glance, there is nothing particularly interesting about this last equation. 
However, both Equations \eqref{eq:casa-loss} and \eqref{eq:casa-semiadjoint} involve $s(t)$ in a VJP. 
Since Equation \eqref{eq:casa-semiadjoint} must hold for every function $\lambda(t)$, we can pick $\lambda(t)$ to make the terms involving $s(t)$ in Equations \eqref{eq:casa-loss} and \eqref{eq:casa-semiadjoint} to perfectly coincide. 
This is done by defining the adjoint $\lambda(t)$ to be the solution of the new system of differential equations
\begin{equation}
    \frac{d\lambda}{dt} 
    = 
    - 
    \frac{\partial f}{\partial u}^T \lambda  
    - 
    \frac{\partial h^T}{\partial u} 
    \qquad \quad \lambda(t_1) = 0. 
    \label{eq:casa-adjoint-equation}
\end{equation}
Notice that the adjoint equation is defined with the final condition at $t_1$, meaning that it needs to be solved backwards in time from $t_1$ to $t_0$. 
The definition of the adjoint $\lambda(t)$ as the solution of this last ODE simplifies Equation \eqref{eq:casa-semiadjoint} to
\begin{equation}
    \inttime \frac{\partial h}{\partial u} s(t) dt
    = 
    \lambda(t_0)^T s(t_0)
    + 
    \inttime \lambda(t)^T \frac{\partial f}{\partial \theta} dt.
\end{equation}
Finally, replacing this inside the expression for the gradient of the loss function we have 
\begin{equation}
    \frac{dL}{d\theta}
    = 
    \lambda(t_0)^T s(t_0)
    + 
    \inttime
    \left( \frac{\partial h}{\partial \theta} + \lambda^T \frac{\partial f}{\partial \theta} \right) dt.
    \label{eq:casa-final-loss-gradient}
\end{equation}
The full algorithm to compute the full gradient $\frac{dL}{d\theta}$ can be described as follows: (i) Solve the original ODE given by $\frac{du}{dt} = f(u, t, \theta)$, (ii) Solve the reverse adjoint ODE given by Equation \eqref{eq:casa-adjoint-equation}, (iii) Compute the gradient using Equation \eqref{eq:casa-final-loss-gradient}.

\subsection{Mathematical comparison of the methods}
\label{section:compatison-math}
In Sections \ref{section:finite-diffences}-\ref{section:continuous-adjoint} we focused on merely introducing each one of the sensitivity methods classified in Figure \ref{fig:scheme-all-methods} as separate methods, postponing the discussion about their points in common. 
In this section, we compare one-to-one these methods and highlight differences and parallels between them.

\subsubsection{Forward AD and complex step differentiation}

\begin{figure}[t]
    \centering
    \includegraphics[width=0.75\textwidth]{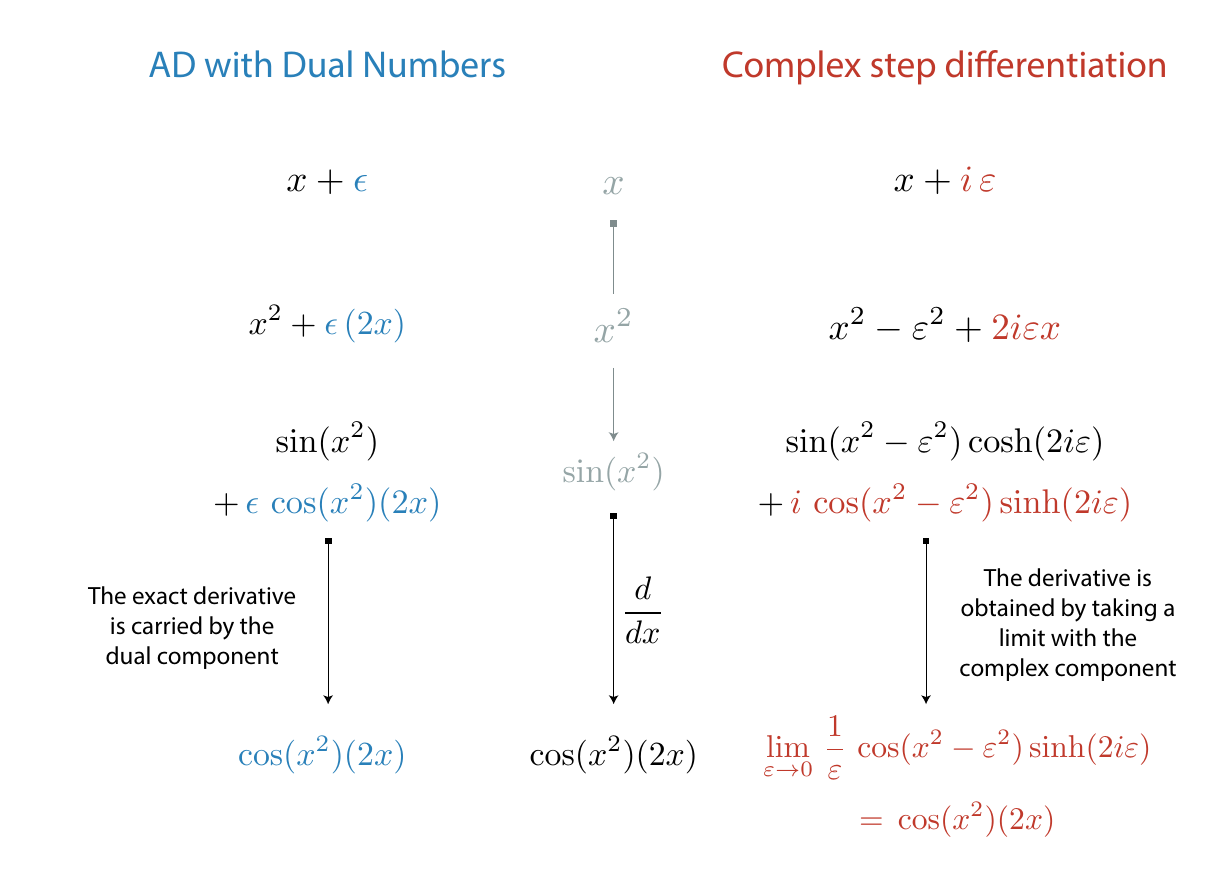}
    \caption{Comparison between AD implemented with dual numbers and complex step differentiation. For the simple case of the function $f(x) = \sin(x^2)$, we can see how each operation is carried in the forward step by the dual component (blue) and the complex component (red). Whereas AD gives the exact gradient at the end of the forward run, in the case of the complex step method we need to take the limit in the imaginary part. }
    \label{fig:complex-step-AD}
\end{figure}

Both AD based on dual numbers and complex-step differentiation introduce an abstract unit, $\epsilon$ and $i$, respectively, associated with the imaginary part of the dual variable that carries forward the numerical value of the gradient.
Although these methods seem similar, AD gives the exact gradient value, whereas complex step differentiation relies on numerical approximations that are valid only when the stepsize $\varepsilon$ is small. 
In Figure \ref{fig:complex-step-AD} we show how the calculation of the gradient of the function $\sin (x^2)$ is performed by these two methods.
Whereas the second component of the dual number has the exact derivative of the function $\sin(x^2)$ with respect to $x$, it is not until we take $\varepsilon \rightarrow 0$ that we obtain the derivative in the imaginary component for the complex step method.
The dependence of the complex step differentiation method on the step size gives it a closer resemblance to finite difference methods than to AD using dual numbers.
Further, notice the complex step method involves more terms in the calculation, a consequence of the fact that second order terms of the form $i^2 = -1$ are transferred to the real part of the complex number, while for dual numbers the terms associated to $\epsilon^2 = 0$ vanish \cite{martins2001connection}. 

\subsubsection{Discrete adjoints and reverse AD}
\label{section:comparison-discrete-adjoint-AD}

Both discrete adjoint methods and reverse AD are classified as discrete and reverse methods (see Figure \ref{fig:scheme-all-methods}). 
Furthermore, both methods introduce an intermediate adjoint associated with the partial derivative of the loss function (output variable) with respect to intermediate variables of the forward computation.
In the case of reverse AD, this adjoint is defined with the notation $\bar w$ (Equation \eqref{eq:reverse-mode-ad-definition}), while in the discrete adjoint method this correspond to each one of the variables $\lambda_1, \lambda_2, \ldots, \lambda_M$ (Equation \eqref{eq:linea-adjoint-state-equation}).
In this section we show that both methods are mathematically equivalent \cite{Zhu_Xu_Darve_Beroza_2021, li2020coupled}, but naive implementations using reverse AD can result in sub-optimal performance compared to that obtained by directly employing the discrete adjoint method \cite{Alexe_Sandu_2009}. 

In order to have a better idea of how this works in the case of a numerical solver, let us consider again the case of a one-step explicit method, not necessarily linear, where the updates $u^m$ satisfy the equation $u^{m+1} = g_{m}(u^{m}; \theta)$.
Following the same schematics as in Figure \ref{fig:discrete-adjoint}, we represent the computational graph of the numerical method using the intermediate variables $g_1, g_2, \ldots, g_{M-1}$.
The dual/adjoint variables defined in reverse AD in this computational graph are given by 
\begin{equation}
    \bar g_m^T
    = 
    (\bar u^{m+1})^T \frac{\partial u^{m+1}}{\partial g_m} 
    = 
    (\bar g_{m+1})^T \frac{\partial g_{m+1}}{\partial u^{m+1}}
    + 
    \left(\frac{\partial L}{\partial u^{m+1}}\right)^T.
\end{equation}
The updates of $\bar g_m$ then mathematically coincide with the updates in reverse mode of the adjoint variable $\lambda_m$ (see Equation \eqref{eq:adjoint-discrete-linear-example}) mapping between tangent spaces (see Section \ref{sec:vjp-jvp}).

\begin{figure}[t]
    \centering
    \includegraphics[width=0.65\textwidth]{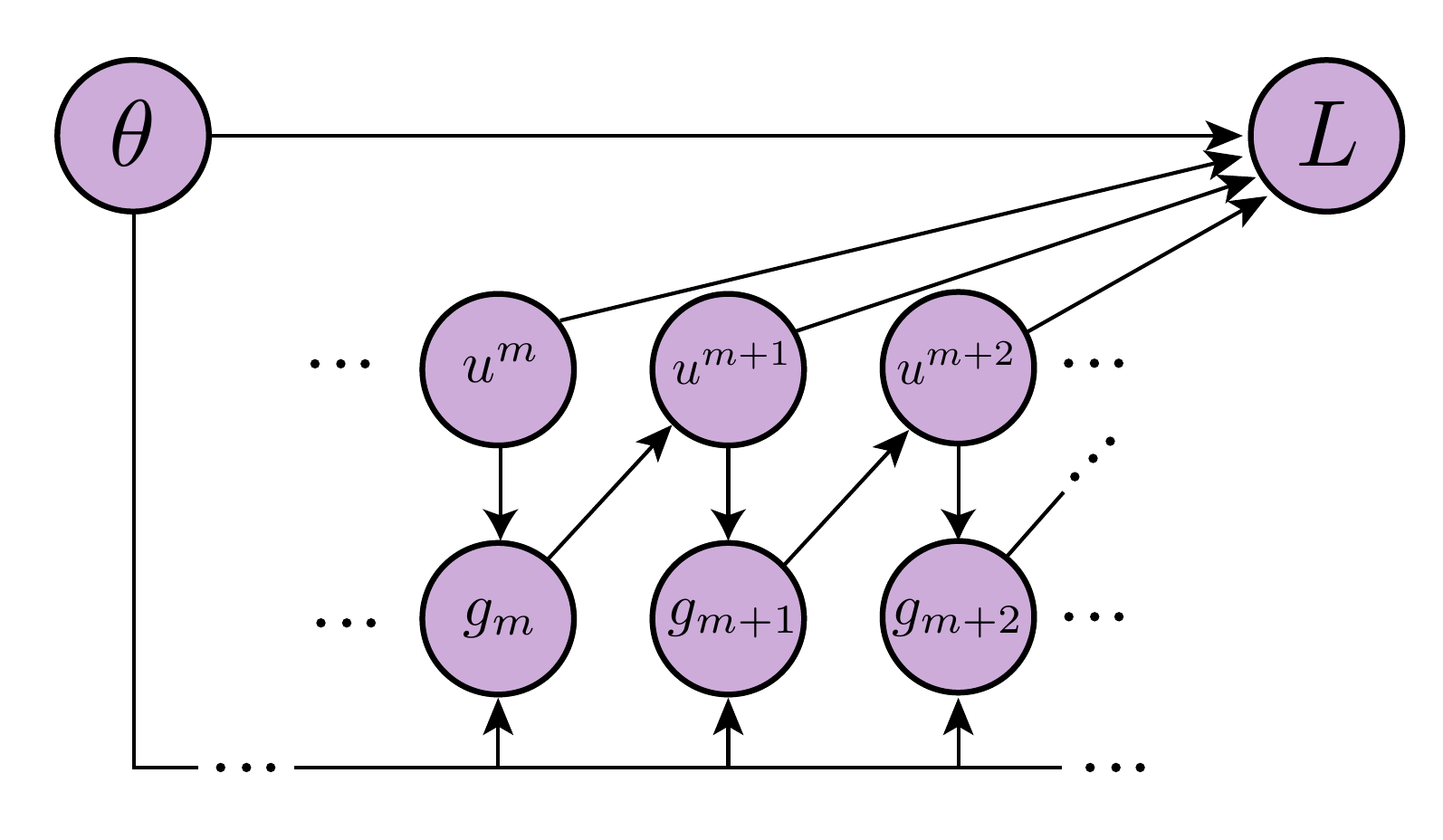}
    \caption{Computational graph associated to the discrete adjoint method. Reverse AD applied on top of the computational graph leads to the update rules for the discrete adjoint. The adjoint variable $\lambda_i$ in the discrete adjoint method coincides with the adjoint variable $\bar g_i$ defined in the backpropagation step.}
    \label{fig:ad-vs-discrete-adjoint}
\end{figure}

Modern numerical solvers use functions $g_m$ that correspond to nested functions, meaning $g_m = g_m^{(k_m)} \circ g_m^{(k_m-1)} \circ \ldots \circ g_m^{(1)}$. 
This is certainly the case for implicit methods when $u^{m}$ is computed as the solution of $g_m(u^m; \theta) = 0$ using an iterative Newton method \cite{SUNDIALS-hindmarsh2005sundials}; or in cases where the numerical solver includes internal iterative sub-routines \cite{Alexe_Sandu_2009}.
If the number of intermediate function is large, reverse AD will result in a large computational graph, potentially leading to excessive memory usage and slow computation \cite{Margossian.2019, Alexe_Sandu_2009}.
A solution to this problem is to introduce a customized \textit{super node} that directly encapsulates the contribution to the full adjoint in $\bar g_m$ without computing the adjoint for each intermediate function $g_m^{(j)}$.
Provided with the value of the Jacobian matrices $\frac{\partial g_m}{\partial u^m}$ and $\frac{\partial g_m}{\partial \theta}$, we can use the implicit function theorem to find $\frac{\partial u^m}{\partial \theta}$ as the solution of the linear system of equations
\begin{equation}
    \frac{\partial g_m}{\partial u^m} \frac{\partial u^m}{\partial \theta} = -\frac{\partial g_m}{\partial \theta}
\end{equation}
and implement AD by directly solving this new system of equations \cite{christianson1994reverse, christianson1998reverse, Bell_Burke_2008, margossian2021efficient}. 
In both cases, the discrete adjoint method can be implemented directly on top of a reverse AD tool that allows customized adjoint calculation \cite{rackauckas2021generalized}. 


\subsubsection{Consistency: forward AD and forward sensitivity equations}
\label{section:forwardAD-sensitivity}

The forward sensitivity equations can also be solved in discrete forward mode by numerically discretizing the original ODE and later deriving the discrete forward sensitivity equations \cite{ma2021comparison}. 
For most cases, this leads to the same result as in the continuous case \cite{FATODE2014}.
We can numerically solve for the sensitivity $s$ by extending the parameter $\theta$ to a multidimensional dual number 
\begin{equation}
    \theta =
    \begin{bmatrix}
    \theta_1 \\
    \theta_2 \\
    \vdots \\
    \theta_p
    \end{bmatrix}
    \longrightarrow
    \begin{bmatrix}
    \theta_1 + \epsilon_1 \\
    \theta_2 + \epsilon_2 \\
    \vdots \\
    \theta_p + \epsilon_p
    \end{bmatrix},
\end{equation}
where $\epsilon_i \epsilon_j = 0$ for all pairs of $i$ and $j$ (see Section \ref{section:dual-numbers}). 
The dependency of the solution $u$ of the ODE on the parameter $\theta$ is now expanded following Equation \eqref{eq:dual-number-function} as 
\begin{equation}
    u =
    \begin{bmatrix}
    u_1 \\
    u_2 \\
    \vdots \\
    u_n
    \end{bmatrix}
    \longrightarrow
    \begin{bmatrix}
    u_1 + \sum_{j=1}^p \frac{\partial u_1}{\partial \theta_j} \epsilon_j \\
    u_2 + \sum_{j=1}^p \frac{\partial u_2}{\partial \theta_j} \epsilon_j \\
    \vdots \\
    u_n + \sum_{j=1}^p \frac{\partial u_n}{\partial \theta_j} \epsilon_j
    \end{bmatrix}
    = 
    u \, + \, s \, 
    \begin{bmatrix}
    \epsilon_1 \\
    \epsilon_2 \\
    \vdots \\
    \epsilon_p
    \end{bmatrix},
\end{equation}
that is, the dual component of the vector $u$ corresponds exactly to the sensitivity matrix $s$. 
This implies forward AD applied to any multistep linear solver will result in the application of the same solver to the forward sensitivity equation (Equation \eqref{eq:sensitivity_equations}).  
For example, for the forward Euler method this gives 
\begin{align}
\begin{split}
    u^{m+1} + s^{m+1} \, \epsilon 
    &= 
    u^m +  s^m \, \epsilon + \Delta t_m \, f (u^m + s^m \, \epsilon, \theta + \epsilon, t_m) \\
    &= 
    u^m + \Delta t_m \, f(u^m, \theta, t_m) 
    + 
    \left [
    s^m
    + 
    \Delta t_m \,
    \left( 
    \frac{\partial f}{\partial u} s^m + 
    \frac{\partial f}{\partial \theta}
    \right) 
    \right ] \epsilon.
\end{split}
\label{eq:sensitivity-equation-AD}
\end{align}
The dual component corresponds to the forward Euler discretization of the forward sensitivity equation, with $s^m$ the temporal discretization of the sensitivity $s(t)$.

The consistency result for discrete and continuous methods also holds for Runge-Kutta methods \cite{Walther_2007}. 
When introducing dual numbers, the Runge-Kutta scheme in Equation \eqref{eq:Runge-Kutta-scheme} gives the following identities
\begin{align}
    u^{m+1} + s^{m+1} \epsilon 
    &= 
    u^m + s^m  \epsilon + \Delta t_m \sum_{i=1}^s b_i (k_i + \dot k_i \epsilon)
    \\
    k_i + \dot k_i \epsilon
    &= 
    f \left(u^m + \sum_{j=1}^s a_{ij} k_j + \left( s^m + \sum_{j=1}^s a_{ij} \dot k_j \right) \epsilon , \theta + \epsilon ,  t_m + c_i \Delta t_m \right), 
    \label{eq:rk-sensitivity-2}
\end{align}
with $\dot k_i$ the dual variable associated to $k_i$.
The partial component in Equation \eqref{eq:rk-sensitivity-2} carrying the coefficient $\epsilon$ gives 
\begin{align}
\begin{split}
    \dot k_i
    &= 
    \frac{\partial f}{\partial u} 
    \left(u^m + \sum_{j=1}^s a_{ij} k_j, \theta,  t_m + c_i \Delta t_m \right)
    \left( s^m + \sum_{j=1}^s a_{ij} \dot k_j \right) \\
    &+ 
    \frac{\partial f}{\partial \theta} 
    \left(u^m + \sum_{j=1}^s a_{ij} k_j, \theta,  t_m + c_i \Delta t_m \right),
\end{split}
\end{align}
which coincides with the Runge-Kutta scheme we would obtain for the original forward sensitivity equation. 
This means that forward AD on Runge-Kutta methods leads to solutions for the sensitivity that have the same convergence properties of the forward solver. 

Note that consistency does not imply that an ODE solver is necessarily correct or stable under such a transformation. 
Consistency of the adjoint may involve other aspects of the solver, such as adaptivity, error control, and the choice of the discretization scheme. 
A common case where continuous methods may fail is when the discretization step is applied without controlling for the join error of the solution of the DE and its sensitivity \cite{Gunzburger_2002}. 
In Section \ref{section:AD-incorrect}, we demonstrate that common implementations of adaptive ODE solvers may not compute the right gradient when forward AD is applied to solver even though the process is mathematically consistent. 
This highlights that additional factors beyond consistency must be considered when investigating whether an implementation is convergent.

\subsubsection{Consistency: discrete and continuous adjoints}

As previously mentioned, the difference between the discrete and continuous adjoint methods is that the former follows the discretize-then-differentiate approach (also known as finite difference of adjoints \cite{Sirkes_Tziperman_1997}).
In contrast, continuous adjoint equations are derived analytically, without a priori consideration of the numerical scheme used to solve it. 
In some sense, we can think of the discrete adjoint $\lambda = (\lambda_1, \lambda_2, \ldots, \lambda_M)$ in Equation \eqref{eq:linea-adjoint-state-equation} as the discretization of the continuous adjoint $\lambda(t)$. 

A natural question then is whether these two methods effectively compute the same gradient, i.e., if the discrete adjoint consistently approximate its continuous counterpart. 
In general, discrete and continuous adjoints will lead to different numerical solutions for the sensitivity, meaning that the discretization and differentiation step do not commute \cite{Jensen_Nakshatrala_Tortorelli_2014, Gunzburger_2002, nadarajah2000comparison}.
However, as the error of the numerical solver decreases, we further expect the discrete and continuous adjoint to lead to the same correct solution. 
Furthermore, since the continuous adjoint method requires to numerically solve the adjoint, we are interested in the relative accuracy of the forward and reverse step. 
It has been shown that for both explicit and implicit Runge-Kutta methods, as long as the coefficients in the numerical scheme given in Equation \eqref{eq:Runge-Kutta-scheme} satisfy the condition $b_i \neq 0$ for all $i=1,2, \ldots, s$, then the discrete adjoint is a consistent estimate of the continuous adjoint with same level of convergence as for the forward numerical solver \cite{Hager_2000,Walther_2007, sandu2006properties, sandu2011solution}.
To guarantee the same order of convergence, it is important that both the forward and backward solver use the same Runge-Kutta coefficients \cite{Alexe_Sandu_2009}.
Importantly, even when consistent, the code generated using the discrete adjoint using AD tools (see Section \ref{section:comparison-discrete-adjoint-AD}) can be sub-optimal and manual modification of the differentiation code is required to guarantee correctness \cite{Eberhard_Bischof_1996, alexe2007denserks}.

\section{Implementation: A computer science perspective}
\label{sec:computational-implementation}
Realizing the full potential of DP requires software that is not only efficient but also sustainable and reproducible. 
Achieving this requires the adoption of established research software engineering (RSE) practices \cite{Combemale_Gray_Rumpe_2023}, including the development of codes that are modular, adaptable, and backend-agnostic.
Such designs enable reuse, composability, and portability across diverse scientific and engineering applications and computational environments \cite{Schäfer_Tarek_White_Rackauckas_2021, naumann2025differentiable, moses_Enzyme}.
In particular, composability enables straightforward generation of differentiable code when combining modular blocks of software.
These requirements reduces the gap between forward and inverse modelling, making new inverse pipelines and the evaluation of sensitivities or gradients easier to prototype, implement, and maintain.   

In this section, we address how the different methods introduced in Section \ref{section:methods} are implemented computationally and how to decide which method to use depending on the scientific task.
In order to address this, it is important to make one further distinction between methods that apply direct differentiation at the algorithmic level and those that are based on numerical solvers.  
The former require a much different implementation since they are agnostic with respect to the mathematical and numerical properties of the ODE.
The latter family of methods that are based on numerical solvers include the forward sensitivity equations and the adjoint methods.
This section is then divided in two parts:
\begin{itemize}
    \item[$ \blacktriangleright$] \textbf{Direct methods.} (Section \ref{section:direct-methods}) Their implementation occurs at a higher hierarchy than the numerical solver software. They include finite differences, AD, complex step differentiation.
    \item[$ \blacktriangleright$] \textbf{Solver-based methods.} Their implementation occurs at the same level of the numerical solver. They include 
    \begin{itemize}
        \item [$\vartriangleright$] Forward sensitivity equations (Section \ref{section:computing-sensitivity-equations})
        \item [$\vartriangleright$] Discrete and continuous adjoint methods (Section \ref{section:computing-adjoints})
    \end{itemize}
\end{itemize}
While these methods can be implemented in different programming languages, we consider examples based on the Julia programming language. 
Julia is a recent but mature programming language that already has a large tradition in implementing packages aiming to advance DP \cite{Bezanson_Karpinski_Shah_Edelman_2012, Julialang_2017}, with a strong emphasis on DE solvers \cite{Rackauckas_Nie_2016, rackauckas2020universal}.
Nevertheless, in reviewing existing work, we also point to applications developed in other programming languages.

The GitHub repository \url{https://github.com/ODINN-SciML/DiffEqSensitivity-Review} contains both text and code used in this manuscript. 
See Appendix \ref{appedix:code} for a complete description of the scripts provided. 
We use the symbol $\clubsuit$ to reference code. 
\subsection{Direct methods}
\label{section:direct-methods}
Direct methods are implemented independent of the structure of the DE and the numerical solver used. 
These include finite differences, complex step differentiation, and both forward and reverse mode AD. 

\subsubsection{Finite differences}
\label{section:software-finite-differences}

\begin{figure}[tb]
    \centering
    \includegraphics[width=1.0\textwidth]{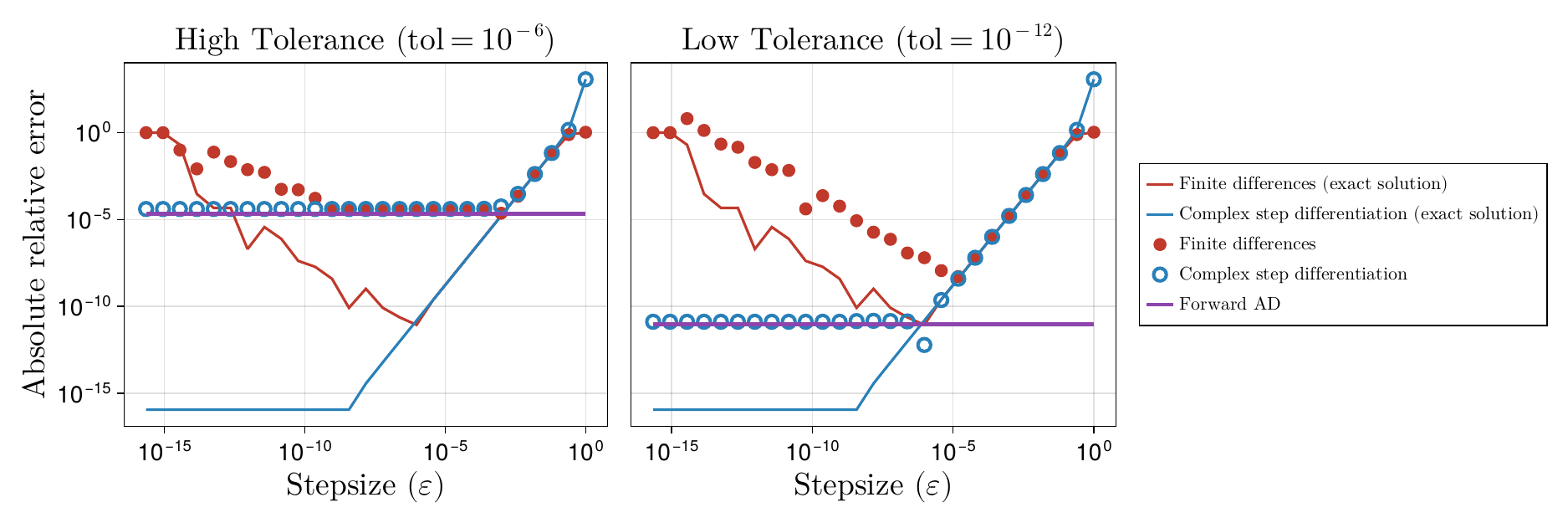}
    \caption{Absolute relative error when computing the gradient of the function $u(t) = \sin (\theta t)/\theta$ with respect to $\theta$ at $t=10.0$ as a function of the stepsize $\varepsilon$ for different direct methods. Here, $u(t)$ corresponds to the solution of the differential equation $u'' + \theta^2 u = 0$ with initial condition $u(0)=0$ and $u'(0)=1$. The red and blue lines correspond to the case where differentiation is applied on the analytical solution, then providing the baseline independent of the numerical solver. Dots correspond to the differentiation of the numerically computed solution using the default Tsitouras solver \cite{Tsitouras_2011} from \texttt{OrdinaryDiffEq.jl} using different solver tolerances, in this case $10^{-6}$ in the left panel, and $10^{-12}$ for the right panel (see Equation \eqref{eq:internal-norm-wrong}). The violet horizontal lines correspond to AD (either forward or reverse). The error when using a numerical solver is larger and it is dependent on the numerical tolerance of the numerical solver. $\clubsuit_\text{\code{code:figure-comparison}}$}
    \label{fig:direct-methods}
\end{figure}

Finite differences are easy to implement manually, do not require much software support, and provide a direct way of approximating a directional derivative. 
In Julia, these methods are implemented in \texttt{FiniteDiff.jl} and \texttt{FiniteDifferences.jl}, which already include subroutines to determine optimal step-sizes.
However, finite differences are less accurate and as costly as forward AD \cite{Griewank_1989} and complex-step differentiation. 
Figure \ref{fig:direct-methods} illustrates the error in computing the derivative of a simple loss function for both true analytical solution and numerical solution of a system of ODEs as a function of the stepsize $\varepsilon$ using finite differences.
Here we consider the solution of the simple harmonic oscillator $u'' + \theta^2 u = 0$ with initial condition $u(0)=0$ and $u'(0)=1$, which has analytical solution $u^\text{true}_\theta(t) = \sin(\theta t) / \theta$.
The numerical solution $u_\theta^\text{num}(t)$ can be obtained by solving the following ODE:
\begin{equation}
\begin{cases}
    \frac{du_1}{dt} = u_2 \,   & \qquad u_1(0) = 0 \\
    \frac{du_2}{dt} = - \theta^2 u_1 \,   & \qquad u_2(0) = 1.
    \label{eq:example-ode-direct-methods}
\end{cases}
\end{equation}
We use $L(\theta) = u_\theta(t_1)$ as our loss function, so that $\frac{dL}{d\theta} = (t_1 / \theta) \cos(\theta t_1)  - \sin (\theta t_1) / \theta^2$ for $t_1=10$.
Finite differences are inaccurate for computing the derivative of $u_\theta^\text{true}$ with respect to $\theta$ when the stepsize $\varepsilon$ is both too small and too large (red line), with a minimum error for $\varepsilon \approx 10^{-6}$.
This case is idealistic as $u_\theta^\text{true}$ cannot generally be obtained analytically, so its derivative obtained using finite differences just serves as a lower bound of the error we expect to see when performing sensitivity analysis on top of the numerical solver. 
When the derivative is instead computed using the numerical solution $u_\theta^\text{num}(t)$ (red circles), the accuracy of the derivative further deteriorates due to approximation errors in the solver. 
This effect is dependent on the numerical solver tolerance. 
For this experiment, both relative and absolute tolerances of the numerical solver had been set to $10^{-6}$ (high tolerance) and $10^{-12}$ (low tolerance) (see Section \ref{section:intro-numerical-solvers}).

\subsubsection{Automatic differentiation}

The AD algorithms described in Section \ref{section:automatic-differentiation} can be implemented using different strategies, namely operator overloading for AD based on dual numbers, and source code transformation for both forward and reverse AD based on the computational graph \cite{martins2001connection}.
In Section \ref{section:software-Forward-AD} we first cover how forward AD is implemented using dual numbers, postponing the discussion about the implementation using computational graphs for reverse AD in Section \ref{sec:software-reverse-AD}. 

\paragraph{Forward AD based on dual numbers}
\label{section:software-Forward-AD}

Implementing forward AD using dual numbers is usually carried out using operator overloading \cite{Neuenhofen_2018}. 
This means expanding the object associated with a numerical value to include the tangent and extending the definition of atomic algebraic functions. 
In Julia, this can be done by relying on multiple dispatch \cite{Julialang_2017}. 
The following example illustrates how to define a dual number and its associated binary addition and multiplication extensions $\clubsuit_\text{\code{code:dual-number}}$. 
\begin{jllisting}
@kwdef struct DualNumber{F <: AbstractFloat}
    value::F
    derivative::F
end

# Binary sum
Base.:(+)(a::DualNumber, b::DualNumber) = DualNumber(value = a.value + b.value, derivative = a.derivative + b.derivative)

# Binary product 
Base.:(*)(a::DualNumber, b::DualNumber) = DualNumber(value = a.value * b.value, derivative = a.value*b.derivative + a.derivative*b.value)
\end{jllisting}
We further overload base operations for this new type to extend the definition of standard functions by simply applying the chain rule and storing the derivative in the dual variable following Equation \eqref{eq:dual-number-function}:
\begin{jllisting}
function Base.:(sin)(a::DualNumber)
    value = sin(a.value)
    derivative = a.derivative * cos(a.value)
    return DualNumber(value=value, derivative=derivative)
end
\end{jllisting}

In the Julia ecosystem, \texttt{ForwardDiff.jl} implements forward mode AD with multidimensional dual numbers \cite{RevelsLubinPapamarkou2016}. 
While \texttt{ForwardDiff.jl} defines the interface for defining the tangent of primitive operations, the tangent of different operations are implemented in \texttt{DiffRules.jl}. 
Figure \ref{fig:direct-methods} shows the result of performing forward AD inside the numerical solver.
We can see that for this simple example forward AD performs as good as the best output of finite differences and complex step differentiation (see Section \ref{section:software-complex-step}) when optimizing by the stepsize $\varepsilon$. 
Implementations of forward AD using dual numbers and computational graphs require a number of operations that increases with the number of variables to differentiate, since each computed quantity is accompanied by the corresponding derivative calculations \cite{Griewank_1989}. 
This consideration also applies to the other forward methods, including finite differences and complex-step differentiation.

\paragraph{Reverse AD based on computational graph}
\label{sec:software-reverse-AD}

In contrast to finite differences, forward AD, and complex-step differentiation, reverse AD is the only of this family of methods that propagates the gradient in reverse mode by relying on analytical derivatives of primitive functions.
The interface for defining primitives in implemented in \texttt{ChainRulesCore.jl}, while the primitives themselves are defined in different libraries (eg, \texttt{ChainRules.jl}, \texttt{SciMLSenstivity.jl}, \texttt{NNlib.jl}). 
Reverse AD can be implemented via pullback functions \cite{Innes_2018}, a method also known as continuation-passing style \cite{Wang_Zheng_Decker_Wu_Essertel_Rompf_2019}.
In the backward step, it executes a series of function calls, one for each elementary operation.
If one of the nodes in the graph $w$ is the output of an operation involving the nodes $v_1, \ldots, v_m$, where $v_i \rightarrow w$ are all edges in the graph, then the pullback $\bar v_1, \ldots, \bar v_m = \mathcal B_w(\bar w)$ is a function that accepts gradients with respect to $w$ (defined as $\bar w$) and returns gradients with respect to each $v_i$ (defined as $\bar v_i$) by applying the chain rule. 
Consider the example of the multiplication $w = v_1 \times v_2$. 
Then
\begin{equation}
 \bar v_1, \, \bar v_2 
 \,=\,
 v_2 \times \bar w , \,
 v_1 \times \bar w 
 \,=\,
 \mathcal{B}_w (\bar w),
\end{equation}
which is equivalent to using the chain rule as
\begin{equation}
 \frac{\partial \ell}{\partial v_1} 
 = \frac{\partial}{\partial v_1}(v_1 \times v_2) \frac{\partial \ell}{\partial w}
 =
 v_2 \times \bar w \, , \qquad 
 \frac{\partial \ell}{\partial v_2} = v_1 \times \bar \omega \, .
\end{equation}

A crucial distinction between AD implementations based on computational graphs is between \textit{static} and \textit{dynamic} methods \cite{Baydin_Pearlmutter_Radul_Siskind_2015}. 
In the case of a static implementation, the computational graph is constructed before any code is executed, which is encoded and optimized for performance within the graph language. 
For static structures such as neural networks, this is ideal, as it simplifies performance optimizations to be applied \cite{abadi-tensorflow}. 
However, two major drawbacks of static methods are composability with existing code, including support of custom types, and adaptive control flow, which is a common feature of numerical solvers. 
In the case of dynamic methods, these issues are addressed using \textit{tracing} or tape-based implementations, where the program structure is transformed into a list of pullback functions that build the graph dynamically at runtime. 
Popular Julia libraries falling in this category are \texttt{Tracker.jl} and \texttt{ReverseDiff.jl}. 
A major drawback of tracing systems is that the pullbacks are constructed with respect to the control flow of the input value and thus do not necessarily generalize to other inputs. 
This means that the pullback must be reconstructed for each forward pass, limiting the reuse of computational optimizations and inducing higher overhead. 
Source-to-source AD systems can achieve higher performance by giving a static derivative representation to arbitrary control flow structure, thus allowing for the construction and optimization of pullbacks independent of the input value. 
These include \texttt{Zygote.jl} \cite{Innes_Zygote}, \texttt{Enzyme.jl} \cite{moses_Enzyme, Moses.2021}, and \texttt{Diffractor.jl}.
The existence of these multiple AD packages lead to the development of \texttt{AbstractDifferentiation.jl} \cite{Schäfer_Tarek_White_Rackauckas_2021} and \texttt{DifferentiationInterface.jl} \cite{dalle_2024_11573435}, which allows one to combine different methods under the same framework. 


\paragraph{Discrete checkpointing}
\label{section:checkpointing}

In contrast to forward methods, all reverse methods, including backpropagation and adjoint methods, require accessing the value of intermediate variables during the propagation of the gradient. 
For a numerical solver or for time-stepping codes, the amount of memory required to accomplish this can be very large, involving a total of at least $\mathcal O(nk)$ terms, with $k$ the number of steps of the numerical solver (or the number of time steps). 
Checkpointing is a technique that can be used for all reverse methods. 
It avoids storing all the intermediate states by balancing  storing and recomputation to recover the required state exactly \cite{Griewank:2008kh}.
This is achieved by saving intermediate states of the solution in the forward pass and recalculating the solution between intermediate states in the reverse mode. 
Different checkpointing algorithms have been proposed, ranging from static or uniform, multi-level \cite{Giering_Kaminski_1998,Heimbach.2005} to optimized, binomial checkpointing algorithms  \cite{Griewank.2000,Walther.2004,Bockhorn.2020,Checkpoiting_2023}.


\paragraph{When AD is algorithmically correct but numerically wrong}
\label{section:AD-incorrect}

Although AD is always algorithmically correct, when combined with a numerical solver \textit{AD can be numerically incorrect} and result in wrong gradient calculations \cite{Eberhard_Bischof_1996}. 
In this section we are going to show an example where AD fails when directly applied to an unmodified solution computed with an adaptive stepsize numerical solver (see Section \ref{section:intro-numerical-solvers}).
When performing forward AD though numerical solver, the error used in the stepsize controller needs to naturally account for both the errors induced in the numerical solution of the original ODE and the errors in the dual component carrying the value of the sensitivity. 
This relation between the numerical solver and AD has been made explicit when we presented the relationship between forward AD and the forward sensitivity equations (Section \ref{section:forwardAD-sensitivity}). 

To illustrate this, let us consider the following first-order ODE:
\begin{equation}
\begin{cases}
 \frac{du_1}{dt} = a u_1 - u_1 u_2 & \quad u_1(0) = 1  \\ 
 \frac{du_2}{dt} = - a u_2 + u_1 u_2 & \quad u_2(0) = 1.
\end{cases}
\end{equation}
Notice that for the value of the parameter $a = 1$, this ODE admits an analytical solution $u_1(t) \equiv u_2(t) \equiv 1$, making this problem very simple to solve numerically.
The following code solves for the derivative with respect to the parameter $a$ using two different methods. 
The second method using forward AD with dual numbers declares the \texttt{internalnorm} argument for the stepsize controller according to Equation \eqref{eq:internal-norm-wrong} $\clubsuit_\text{\code{code:AD-wrong}}$.
\begin{jllisting}
using SciMLSensitivity, OrdinaryDiffEq, Zygote, ForwardDiff

function fiip(du, u, p, t)
    du[1] =  p[1] * u[1] - u[1] * u[2]
    du[2] = -p[1] * u[2] + u[1] * u[2]
end

p = [1.]
u0 = [1.0;1.0]
prob = ODEProblem(fiip, u0, (0.0, 10.0), p);

# Correct gradient computed using 
grad0 = Zygote.gradient(p->sum(solve(prob, Tsit5(), u0=u0, p=p, sensealg = ForwardSensitivity(), saveat = 0.1, abstol=1e-12, reltol=1e-12)), p)
# grad0 = ([212.71042521681443],)

# Original AD with wrong norm 
grad1 = Zygote.gradient(p->sum(solve(prob, Tsit5(), u0=u0, p=p, sensealg = ForwardDiffSensitivity(), saveat = 0.1, internalnorm = (u,t) -> sum(abs2,u/length(u)), abstol=1e-12, reltol=1e-12)), p)
# grad1 = ([6278.15677493293],)
\end{jllisting}
The reason why the two methods give different answers is because the error estimation by the stepsize controller is ignoring numerical errors in the dual component. 
In the later case, since the numerical solution of the original ODE is constant, the local estimated error is drastically underestimated to $\text{err}_i^{m} = 0$, which makes the stepsize $\Delta t_{m}$ to increase by a multiplicative factor at every step (see Equations \eqref{eq:internal-norm-wrong} and \eqref{eq:PIC}). 
This can be fixed by instead considering a norm that accounts for both the primal and dual components in the forward pass, 
\begin{align}
    \text{Err}_\text{scaled}^{m}
    =
    \Bigg [ \frac{1}{n(p+1)} \Bigg( 
    &\sum_{i=1}^n \left( \frac{u_i^{m} - \hat u_i^{m}}{\mathfrak{abstol} + \mathfrak{reltol} \times \max \{ u_i^{m}, \hat u_i^{m}\}} \right)^2 \nonumber \\
    + 
    &\sum_{i=1}^n \sum_{j=1}^p  
    \left( \frac{s_{ij}^{m} - \hat s_{ij}^{m}}{\mathfrak{abstol} + \mathfrak{reltol} \times \max \{ s_{ij}^{m}, \hat s_{ij}^{m}\}} \Bigg)^2 \right)
    \Bigg ]^{\frac{1}{2}},
    \label{eq:internal-norm-correct} 
\end{align}
with $s^m$ and $\hat s^m$ two different numerical approximations of the sensitivity matrix. 
This correction now gives the right answer:
\begin{jllisting}
sse(x::Number) = x^2
sse(x::ForwardDiff.Dual) = sse(ForwardDiff.value(x)) + sum(sse, ForwardDiff.partials(x))

totallength(x::Number) = 1
totallength(x::ForwardDiff.Dual) = totallength(ForwardDiff.value(x)) + sum(totallength, ForwardDiff.partials(x))
totallength(x::AbstractArray) = sum(totallength,x)

grad3 = Zygote.gradient(p->sum(solve(prob, Tsit5(), u0=u0, p=p, sensealg = ForwardDiffSensitivity(), saveat = 0.1, internalnorm = (u,t) -> sqrt(sum(x->sse(x),u) / totallength(u)), abstol=abstol, reltol=reltol)), p)
# grad3 = ([212.71042521681392],)
\end{jllisting}
This is an example where the form of the numerical solver for the original ODE is affected by the fact we are simultaneously solving for the sensitivity. 
Notice that current implementations of forward AD inside \texttt{SciMLSensitivity.jl} already account for this and there is no need to specify the internal norm $\clubsuit_\text{\ref{code:AD-wrong}}$. 
To highlight the pervasiveness of this issue with respect to AD, we further provide a script with an example in \texttt{Diffrax} where the derivative that does not converge to the correct answer as tolerance is decreased to zero $\clubsuit_\text{\code{code:AD-wrong-JAX}}$.



\subsubsection{Complex step differentiation}
\label{section:software-complex-step}

Modern software already have support for complex number arithmetic, making complex step differentiation very easy to implement.
In Julia, complex analysis arithmetic can be easily carried inside the numerical solver.
The following example shows how to extend the numerical solver used to solve the ODE in Equation \eqref{eq:example-ode-direct-methods} to support complex numbers $\clubsuit_\text{\code{code:complex-step}}$.
\begin{jllisting}
function dyn!(du::Array{Complex{Float64}}, u::Array{Complex{Float64}}, p, t)
    ω = p[1]
    du[1] = u[2]
    du[2] = - ω^2 * u[1]
end

tspan = [0.0, 10.0]
du = Array{Complex{Float64}}([0.0])
u0 = Array{Complex{Float64}}([0.0, 1.0])

function complexstep_differentiation(f::Function, p::Float64, ε::Float64)
    p_complex = p + ε * im
    return imag(f(p_complex)) / ε
end

complexstep_differentiation(x -> solve(ODEProblem(dyn!, u0, tspan, [x]), Tsit5()).u[end][1], 20., 1e-3)
\end{jllisting}

Figure \ref{fig:direct-methods} further shows the result of performing complex step differentiation using the same example as in Section \ref{section:software-finite-differences}.
We can see from both exact and numerical solutions that complex-step differentiation does not suffer from small values of $\varepsilon$, meaning that $\varepsilon$ can be chosen arbitrarily small \cite{martins2001connection} as long as it does not reach the underflow threshold \cite{Goldberg_1991_floatingpoint}. 
Notice that for large values of the stepsize $\varepsilon$ complex step differentiation gives similar results to finite differences, while for small values of $\varepsilon$ the performance of complex step differentiation is slightly worse than AD. 
This result emphasizes the observation made in Section \ref{section:comparison-discrete-adjoint-AD}, namely that complex step differentiation has many aspects in common with finite differences and AD based on dual numbers. 

However, the difference between the methods also makes the complex step differentiation sometimes more efficient than both finite differences and AD \cite{Lantoine_Russell_Dargent_2012}, an effect that can be counterbalanced by the number of extra unnecessary operations that complex arithmetic requires (see last column in Figure \ref{fig:complex-step-AD}) \cite{Martins_Sturdza_Alonso_2003_complex_differentiation}.
Further notice that complex-step differentiation will work as long as every function involved in the computation is locally analytical. 
This is a problem with implementations of functions that rely on the absolute function, for example using complex step differentiation on the square function implemented as $f(z) = \text{abs}(z)^2$ will always return zero ($\text{Im}(f(x + i \varepsilon))= \text{Im}((x + i \epsilon)(x - i \epsilon)) = \mathcal O (\varepsilon^2)$).
\subsection{Solver-based methods}
\label{section:solver-methods}
We now move our discussion to DP methods based on numerical solvers.
These need to deal with some numerical and computational considerations, including:
\begin{enumerate}
    \item[$ \blacktriangleright$] How to handle JVPs and/or VJPs 
    \item[$ \blacktriangleright$] Stability of the numerical solver, including the original ODE but also the sensitivity/adjoint equations
    \item[$ \blacktriangleright$] Memory-time tradeoff
\end{enumerate}
These factors are further exacerbated by the size $n$ of the ODE and the number $p$ of parameters in the model. 
Just a few modern scientific software implementations have the capabilities of solving ODE and computing their sensitivities at the same time. 
These include 
\texttt{CVODES} within \texttt{SUNDIALS} in C \cite{serban2005cvodes, SUNDIALS-hindmarsh2005sundials}; 
\texttt{ODESSA} \cite{ODESSA} and \texttt{FATODE} (discrete adjoints) \cite{FATODE2014} both in Fortram; 
\texttt{SciMLSensitivity.jl} in Julia \cite{rackauckas2020universal}; 
\texttt{Dolfin-adjoint} based on the \texttt{FEniCS} Project \cite{dolfin2013, dolfin2018};
\texttt{DENSERKS} in Fortran \cite{alexe2007denserks}; 
\texttt{DASPKADJOINT} \cite{Cao_Li_Petzold_2002};
and \texttt{Diffrax} \cite{kidger2021on} and \texttt{torchdiffeq} \cite{torchdiffeq} in Python. 


\subsubsection{Forward sensitivity equation}
\label{section:computing-sensitivity-equations}

For systems of equations with few number of parameters, the forward sensitivity equation is useful since the system of $n(p+1)$ equations composed by Equations \eqref{eq:original_ODE} and \eqref{eq:sensitivity_equations} can be solved using the same precision for both solution and sensitivity numerical evaluation. 
Furthermore, it does not required saving the solution in memory. 
The following example illustrates how Equation \eqref{eq:example-ode-direct-methods} and the forward sensitivity equation can be solved simultaneously using the simple explicit Euler method $\clubsuit_\text{\code{code:sensitivity-equation}}$: 
\begin{jllisting}
p = [0.2]
u0 = [0.0, 1.0]
tspan = [0.0, 10.0]

# Dynamics
function f(u, p, t)
    du₁ = u[2]
    du₂ = - p[1]^2 * u[1]
    return [du₁, du₂]
end

# Jacobian ∂f/∂p
function ∂f∂p(u, p, t)
    Jac = zeros(length(u), length(p))
    Jac[2,1] = -2 * p[1] * u[1]
    return Jac
end

# Jacobian ∂f/∂u
function ∂f∂u(u, p, t)
    Jac = zeros(length(u), length(u))
    Jac[1,2] = 1
    Jac[2,1] = -p[1]^2
    return Jac
end

# Explicit Euler method
function sensitivityequation(u0, tspan, p, dt)
    u = u0
    sensitivity = zeros(length(u), length(p))
    for ti in tspan[1]:dt:tspan[2]
        sensitivity += dt * (∂f∂u(u, p, ti) * sensitivity + ∂f∂p(u, p, ti))
        u += dt * f(u, p, ti) 
    end
    return u, sensitivity   
end

u, s = sensitivityequation(u0, tspan , p, 0.001)
\end{jllisting}
The simplicity of the sensitivity method makes it available in most software for sensitivity analysis. 
In \texttt{SciMLSensitivity} in the Julia SciML ecosystem, the \texttt{ODEForwardSensitivityProblem} method implements continuous sensitivity analysis, which generates the JVPs required as part of the forward sensitivity equations via \texttt{ForwardDiff.jl} (see Section \ref{section:forwardAD-sensitivity}) or finite differences.
Using \texttt{SciMLSensitivity} reduces the code above to
\begin{jllisting}
using SciMLSensitivity

function f!(du, u, p, t)
    du[1] = u[2]
    du[2] = - p[1]^2 * u[1]
end

prob = ODEForwardSensitivityProblem(f!, u0, tspan, p)
sol = solve(prob, Tsit5())
\end{jllisting}

For stiff systems of ODEs the use of the forward sensitivity equations can be computationally unfeasible \cite{kim_stiff_2021}.
This is because stiff ODEs require the use of stable solvers with cubic cost with respect to the number of ODEs \cite{hairer-solving-2}, making the total complexity of the sensitivity method $\mathcal{O}(n^3p^3)$. 
This complexity makes this method expensive for models with large $n$ and/or $p$ unless the solver is able to further specialize on sparsity or properties of the linear solver (i.e. through Newton-Krylov methods). 

\paragraph{Computing JVPs and VJPs inside the solver}
\label{section:computing-vjp-inside-solver}

An important consideration is that all solver-based methods have subroutines to compute the JVPs and VJPs involved in the sensitivity and adjoint equations, respectively. 
This calculation is carried out by another sensitivity method, usually finite differences or AD, which plays a central role when analyzing the accuracy and stability of the adjoint method. 
In the case of the forward sensitivity equation, this correspond to the JVPs resulting form the product $\frac{\partial f}{\partial u} s $ in Equation \eqref{eq:sensitivity_equations}.
For the adjoint equations, we need to evaluate the term $\lambda^T \frac{\partial f}{\partial \theta}$ for the continuous adjoint method in Equation \eqref{eq:casa-final-loss-gradient}, while for the discrete adjoint method we need to compute the term $\lambda^T \frac{\partial G}{\partial \theta}$ in Equation \eqref{eq:def_adjoint} (which may further coincide with $\lambda^T \frac{\partial f}{\partial \theta}$ for some numerical solvers, but not in the general case). 
Therefore, the choice of the algorithm to compute JVPs/VJPs can have a significant impact in the overall performance \cite{Schäfer_Tarek_White_Rackauckas_2021}. 

In \texttt{SUNDIALS}, the JVPs/VJPs involved in the sensitivity and adjoint method are handled using finite differences unless specified by the user \cite{SUNDIALS-hindmarsh2005sundials}.
In \texttt{FATODE}, they can be computed with finite differences, AD, or it can be provided by the user.
In the Julia SciML ecosystem, the options \texttt{autodiff} and \texttt{autojacvec} allow one to customize if JVPs/VJPs are computed using AD, finite differences, or alternatively these are provided by the user. 
Different AD packages with different performance trade-offs are available for this task (see Section \ref{sec:software-reverse-AD}), including \texttt{ForwardDiff.jl} \cite{RevelsLubinPapamarkou2016}, \texttt{ReverseDiff.jl}, \texttt{Zygote.jl} \cite{Innes_Zygote}, \texttt{Enzyme.jl} \cite{moses_Enzyme}, \texttt{Tracker.jl}.

\subsubsection{Adjoint methods}
\label{section:computing-adjoints}

For complex and large systems (e.g. for $n + p > 50$, as we will discuss in Section \ref{section:recomendations}), direct methods for computing the gradient on top of the numerical solver can be memory expensive due to the large number of function evaluations required by the solver and the later store of the intermediate states. 
For these cases, adjoint-based methods allow us to compute the gradients of a loss function by instead computing the adjoint that serves as a bridge between the solution of the ODE and the final gradient. 
Since adjoint methods rely on an additional set of ODEs which are solved, numerical efficiency and stability must further be taken into account at the moment of implementing adjoint methods.

\paragraph{Discrete adjoint method}

In order to illustrate how the discrete adjoint method can be implemented, the following example shows how to manually solve for the gradient of the solution of \eqref{eq:example-ode-direct-methods} using an explicit Euler method $\clubsuit_\text{\code{code:discrete-adjoint}}$. 
\begin{jllisting}
function discrete_adjoint_method(u0, tspan, p, dt)
    u = u0
    times = tspan[1]:dt:tspan[2]

    λ = [1.0, 0.0]
    ∂L∂θ = zeros(length(p))
    u_store = [u]

    # Forward pass to compute solution
    for t in times[1:end-1]
        u += dt * f(u, p, t)
        push!(u_store, u)
    end

    # Reverse pass to compute adjoint
    for (i, t) in enumerate(reverse(times)[2:end])
        u_memory = u_store[end-i+1]
        λ += dt * ∂f∂u(u_memory, p, t)' * λ
        ∂L∂θ += dt * λ' * ∂f∂p(u_memory, p, t)
    end

    return ∂L∂θ
end

∂L∂θ = discrete_adjoint_method(u0, tspan, p, 0.001) 
\end{jllisting}
In this case, the full solution in the forward pass is stored in memory and used to compute the adjoint and integrate the loss function during the reverse pass. 
While the previous example illustrates a manual implementation of the adjoint method, the discrete adjoint method can be directly implemented using reverse AD (Section \ref{section:comparison-discrete-adjoint-AD}).
In the Julia SciML ecosystem, \texttt{ReverseDiffAdjoint} performs reverse AD on the numerical solver via \texttt{ReverseDiff.jl},  and \texttt{TrackerAdjoint} via \texttt{Tracker.jl}. 
As in the case of reverse AD, checkpointing can be used here.

\paragraph{Continuous adjoint method}

\begin{table}[tb]
    \centering
    \includegraphics[width=1.0\textwidth]{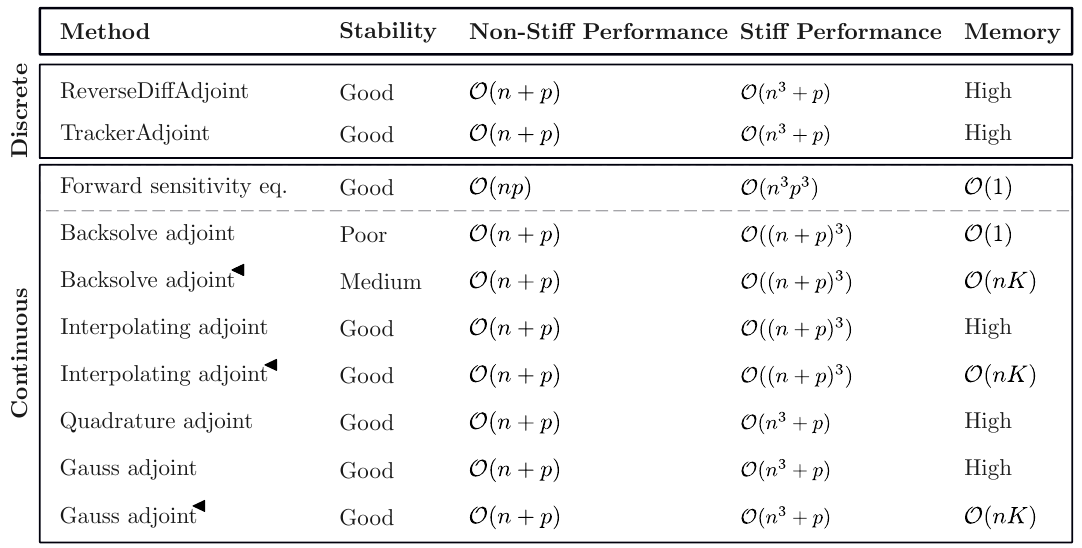}
    \caption{Comparison in performance and cost of solver-based methods. Methods that can be checkpointed are indicated with the symbol $\blacktriangleleft$, with $K$ the total number of checkpoints. The nomenclature of the different adjoint methods here follows the naming in the documentation of \texttt{SciMLSensitivity.jl} \cite{rackauckas2020universal}.}
    \label{table:adjoint}
\end{table}

The continuous adjoint method offers a series of advantages over the discrete method and the rest of the forward methods previously discussed. 
Just as with the discrete adjoint methods and reverse AD, the bottleneck is how to solve for the adjoint $\lambda(t)$ due to its dependency with VJPs involving the state $u(t)$.
Effectively, notice that Equation \eqref{eq:casa-adjoint-equation} involves the terms $f(u, \theta, t)$ and $\frac{\partial h}{\partial u}$, which are both functions of $u(t)$. 
However, in contrast to the discrete adjoint methods, here the full continuous trajectory $u(t)$ is needed, instead of its discrete pointwise evaluation. 
There are two solutions for addressing the evaluation of $u(t)$ during the computation of $\lambda (t)$:
\begin{enumerate}
    \item[$ \blacktriangleright$] \textbf{Interpolation.} During the forward model, we can store in memory intermediate states of the numerical solution allowing the dense evaluation of the numerical solution at any given time. 
    This can be done using dense output formulas, for example by adding extra stages to the Runge-Kutta scheme that allows to define a continuous interpolation, a method known as continuous Runge-Kutta \cite{hairer-solving-2, Alexe_Sandu_2009}. 
    \item[$ \blacktriangleright$] \textbf{Backsolve.} Solve again the original ODE together with the adjoint as the solution of the following reversed augmented system \cite{chen_neural_2019}:
    \begin{equation}
    \frac{d}{dt}
    \begin{bmatrix}
       u \\
       \lambda \\
       \frac{dL}{d\theta}
    \end{bmatrix}
    = 
    \begin{bmatrix}
       -f \\
       - \frac{\partial f}{\partial u}^T \lambda - \frac{\partial h}{\partial u}^T \\
       - \lambda^T \frac{\partial f}{\partial \theta} - \frac{\partial h}{\partial \theta}
    \end{bmatrix}
    \qquad 
    \begin{bmatrix}
       u \\
       \lambda \\
       \frac{dL}{d\theta}
    \end{bmatrix}(t_1)
    = 
    \begin{bmatrix}
       u(t_1) \\
       \frac{\partial L}{\partial u(t_1)} \\
       \lambda(t_0)^T s(t_0)
    \end{bmatrix}.
    \end{equation}
    An important problem with this approach is that computing the ODE backwards $\frac{du}{dt} = - f(u,\theta, t)$ can be unstable and lead to large numerical errors \cite{kim_stiff_2021, Zhuang_2020}. 
    Implicit methods may be used to ensure stability when solving this system of equations. 
    However, this requires cubic time in the total number of ordinary differential equations, leading to a total complexity of $\mathcal O((n+p)^3)$ for the adjoint method. In practice, this method is hardly stable for most complex (even non-stiff) differential equations \cite{kim_stiff_2021}. 
\end{enumerate} 

The following example shows how to implement the continuous adjoint method of the solution of Equation \eqref{eq:example-ode-direct-methods} using the backsolve strategy $\clubsuit_\text{\code{code:continuous-adjoint}}$. 
\begin{jllisting}
using RecursiveArrayTools

# Augmented dynamics
function f_aug(z, p, t)
    u, λ, L = z
    du = f(u, p, t)
    dλ = ∂f∂u(u, p, t)' * λ
    dL = λ' * ∂f∂p(u, p, t)
    VectorOfArray([du, vec(dλ), vec(dL)])
end

# Solution of original ODE
prob = ODEProblem(f, u0, tspan, p)
sol = solve(prob, Euler(), dt=0.001)

# Final state 
u1 = sol.u[end]
z1 = VectorOfArray([u1, [1.0, 0.0], zeros(length(p))])

aug_prob = ODEProblem(f_aug, z1, reverse(tspan), p)
u0_, λ0, dLdp_cont = solve(aug_prob, Euler(), dt=-0.001).u[end]
\end{jllisting}
Notice that here we used the final state of the solution $u$ of the ODE as starting point, which is then recalculated in backwards direction (implemented via a negative stepsize \texttt{dt=-0.001}). 

When dealing with stiff DEs, special considerations need to be taken into account.
Two alternatives are proposed in \textcite{kim_stiff_2021}, the first referred to as quadrature adjoint produces a high order interpolation of the solution $u(t)$, then solve for $\lambda$ in reverse using an implicit solver and finally integrating $\frac{dL}{d\theta}$ in a forward step.
This reduces the complexity to $\mathcal O (n^3 + p)$, where the cubic cost in the size $n$ of the ODE comes from the fact that we still need to solve the original stiff ODE in the forward step. 
A second similar approach is to use an implicit-explicit (IMEX) solver, where we use the implicit part for the original equation and the explicit for the adjoint. 
This method also has a complexity of $\mathcal O (n^3 + p)$. 

\paragraph{Continuous checkpointing}
\label{section:checkpointint-cont}

Both interpolating and backsolve adjoint methods can be implemented along with a checkpointing scheme. 
This can be done by choosing saved points in the forward pass. 
For the interpolating methods, the interpolation is reconstructed in the backwards pass between two save points. 
This reduces the total memory requirement of the interpolating method to simply the maximum cost of holding an interpolation between two save points, but requires a total additional computational effort equal to one additional forward pass. 
In the backsolve variation, the value $u$ in the reverse pass can be corrected to be the saved point, thus resetting the numerical error introduced during the backwards evaluation and thus improving the accuracy.

\paragraph{Solving the quadrature}

Another computational consideration is how the integral in Equation \eqref{eq:casa-final-loss-gradient} is numerically evaluated. 
While one can solve the integral simultaneously with the other equations using an ODE solver, this is only recommended with explicit methods as with implicit methods these additional ODEa is of size $p$ and thus can increase the complexity of an implicit solve by $O(p^3)$. 
The interpolating adjoint and backsolve adjoint methods use this ODE solver approach for computing the integrand. 
On the other hand, the quadrature adjoint approach avoids this computational cost by computing the dense solution $\lambda(t)$ and then computing the quadrature 
\begin{equation}
    \inttime
    \left( \frac{\partial h}{\partial \theta} + \lambda^T \frac{\partial f}{\partial \theta} \right) dt 
    \approx
    \sum_{j=1}^N \omega_j \left( \frac{\partial h}{\partial \theta} + \lambda^T \frac{\partial f}{\partial \theta} \right) (\tau_i),
\end{equation}
where $\omega_i$, $\tau_i$ are the weights and knots of a Gauss-Kronrod quadrature method for numerical integration from \texttt{QuadGK.jl} \cite{laurie1997calculation, gonnet2012review}. 
This method results in global error control on the integration and removes the cubic scaling within implicit solvers. Nonetheless, it requires a larger memory cost by storing the adjoint pass continuous solution.

Solvers designed for large implicit systems allow for solving explicit integrals based on the ODE solution simultaneously without including the equations in the ODE evaluation in order to avoid this expense. 
The Sundials CVODE solver introduced this technique specifically for BDF methods \cite{SUNDIALS-hindmarsh2005sundials}. 
In the Julia \texttt{DifferentialEquations.jl} solvers, this can be done using a callback (specifically the numerical integration callbacks form the \texttt{DiffEqCallbacks.jl} library). 
The Gauss adjoint method uses the callback approach to allow for a simultaneous explicit evaluation of the integral using Gaussian quadrature, similar to \textcite{Norcliffe_gaussquadrature_2023} but using a different approximation to improve convergence.
\\ \\ 
These differences in the strategies for computing $u(t)$ and the final quadrature give rise to the set of methods in Table \ref{table:adjoint}. 
Excluding the forward sensitivity equation which was added for reference, all of these are \textit{the adjoint method} with differences being in the way steps of the adjoint method are approximated, and notably Table \ref{table:adjoint} shows a general trade-off in stability, performance, and memory across the methods. 
While the Gauss adjoint achieves good properties according to this chart, the quadrature adjoint notably uses a global error control of the quadrature as opposed to the local error control of the Gauss adjoint, and thus can achieve more robust bounding of the error with respect to user chosen tolerances.

\section{Generalizations}
\label{section:generalization}
In this section, we briefly discuss how the ideas covered in Sections \ref{section:methods} and \ref{sec:computational-implementation} for first-order ODEs generalize to more complicated systems of DEs. 

Notice that the application of all the direct methods (finite differences, AD, complex step differentiation, symbolic differentiation) applies to more general systems of DEs.
The fundamental behaviour and implementation of these DP methods does not change, although new considerations about numerical accuracy may need to be taken into account, especially for discrete methods based on unmodified solution processes.  
The mathematical derivation of continuous methods (forward sensitivity equations and continuous adjoint) covered in Sections \ref{section:sensitivity-equation} and \ref{section:continuous-adjoint} still applies, although more specific details of the DEs may apply (e.g., inclusion for boundary conditions or other constraints). 
Regarding discrete adjoint methods, the mathematical formulation covered in Section \ref{section:discrete-adjoint} and its connection with reverse AD (Section \ref{section:comparison-discrete-adjoint-AD}) applies to more general solvers for DEs. 

In the next section we are going to consider the cases of higher-order ODEs, PDEs, chaotic systems of ODEs, and stochastic differential equations (SDEs). 
Further generalizations of sensitivity methods to other families of DEs include differential algebraic equations (DAEs) (see \textcite{Cao_Li_Petzold_2002, margossian2021efficient}), delay differential equations (DDEs) (see \textcite{calver2017numerical}, \textcite{rihan2018inverse}), among others.

\subsection{Higher-order ODEs}

Higher-order ODEs are characterized by the presence of second and higher-order time derivatives in the differential equation. 
A simple example popular in structural design consists of the linear dynamic equations used to model elastic structures given by 
\begin{equation}
    M \frac{d^2 u}{dt^2}
    +
    C \frac{du}{dt}
    + 
    K u
    = 
    F(t, \theta),
\end{equation}
subject to the initial condition $u(t_0) = u_0 \in \R^n$, $\frac{du}{dt}(t_0) = v_0 \in \R^n$, where $M, C, K \in \R^{n \times n}$ are the mass, damping, and stiffness matrices function of some design parameter $\theta$, respectively, and $F(t, \theta)$ is an external forcing \cite{min1999optimal, Jensen_Nakshatrala_Tortorelli_2014}.  

Just as we did in Section \ref{section:software-finite-differences}, higher-order ODEs can be transformed to first-order ODEs, after which the same sensitivity methods we discussed in this review can be used. 
However, there may be reasons why we would prefer to avoid this, such as the existence of more efficient higher-order ODEs solvers, including Nystr\"{o}m methods for the case when $C = 0$ \cite{Butcher_Wanner_1996, hairer-solving-1}. 
In this case, the forward sensitivity equations can be derived using the same strategy explained in Section \ref{section:sensitivity-equation}: 
\begin{equation}
    \frac{d}{d\theta} 
    \left(
    M \frac{d^2 u}{dt^2}
    + 
    K u
    -
    F(t, \theta)
    \right) = 0,
\end{equation}
which results in the forward sensitivity equation for the sensitivity $s(t)$:
\begin{equation}
    M \frac{d^2 s}{dt^2}
    + 
    K s
    = 
    \frac{\partial F}{\partial \theta} 
    - 
    \frac{dM}{d\theta} \frac{d^2 u}{dt^2}
    - 
    \frac{dK}{d\theta} u. 
\end{equation}
Similarly, the same strategy introduced Section \ref{section:continuous-adjoint} can be followed to derive the continuous adjoint equation \cite{kang2006review}. 

\subsection{Partial differential equations}




Systems of partial differential equations (PDEs) include derivatives with respect to more than one independent variable. 
As we discussed in Section \ref{sectopn:motivation}, PDEs play a central role in mathematics, physics, and engineering, where these variables are usually associated to time and space. 
Due to the spatial characteristics of such systems, generally boundary conditions need to be provided in addition to initial conditions. 
While in Section \ref{section:intro-numerical-solvers} we briefly introduced the fundamentals for numerical solvers of ODEs, there is a broader family of numerical methods to solve PDEs. 
These include the finite element method and the finite volume method, among others \cite{tadmor2012review}. 
For these methods, a required ingredient is the spatial mesh used to discretize the spatial dimension \cite{thompson1998handbook}. 
In the case of the discrete adjoint method, all these methods will result in a series of discrete equations where the adjoint method introduced in Section \ref{section:discrete-adjoint} will still apply. 
Continuous methods require a more careful manipulation of the PDE in order to derive correct sensitivity and adjoint equations. 

The method of lines can be used to solve PDEs by applying a semi-discretization in the spacial coordinate and then numerically solve a new system of ODEs \cite{ascher2008numerical}. 
This implies that all sensitivity methods for ODEs also apply to PDEs. 
Let us consider the case of the one-dimensional heat equation
\begin{align}
\begin{split}
 \frac{\partial u}{\partial t}
 &= 
 D(x,t) \, 
 \frac{\partial^2 u}{\partial x^2} \qquad x \in [0,1], \, t \in [t_0, t_1]\\
 u(x, t_0) &= v(x) \\
 u(0, t) &= \alpha(t) \\
 u(1, t) &= \beta(t),
 \label{eq:heat-equation}
\end{split}
\end{align}
with $D(x, t) > 0$ a global diffusivity coefficient.
In order to numerically solve this equation, we can define a uniform spatial mesh with coordinates $m \Delta x$, $m=0, 1, 2, \ldots, N$ and $\Delta x = 1 / N$.
If we call $u_m(t) = u(m \Delta x, t)$ and $D_m(t) = D(m\Delta x, t)$ the values of the solution and the diffusivity evaluated in the fixed points in the mesh, respectively, then we can replace the second order partial derivative in Equation \eqref{eq:heat-equation} by the corresponding second order finite difference
\begin{equation}
 \frac{d u_m}{dt} 
 = 
 D_m(t) \, 
 \frac{u_{m-1} - 2u_m + u_{m+1}}{\Delta x^2}
 \label{eq:heat-equation-discrete}
\end{equation}
for $m = 1, 2, \ldots, N-1$ (in the boundary we simply have $u_0(t) = \alpha(t)$ and $u_N(t)=\beta(t)$).
Now, following this semi-discretization, equation \eqref{eq:heat-equation-discrete} is a system of first-order ODEs of size $N-1$.
Semi-discretized PDEs typically involve large systems of coupled and possibly stiff ODEs subject to some suitable boundary conditions. 
Explicit calculation of the Jacobian quickly becomes cumbersome and eventually intractable as the spatial dimension and the complexity of the PDE increase.
Further improvements can be made by exploiting the fact that the coupling in the ODE is sparse, that is, the temporal derivative depends on the state value of the solution in the neighbouring points in the mesh.
PDEs are often also subject to additional time stepping constraints, such as the Courant-Fredrichs-Lewy (CFL) condition, which may limit the maximum time step size and thus increase the number of time steps required to obtain a valid solution \cite{courantPartialDifferenceEquations1967}. 

Besides the methods of lines which already involves a first discretization of the original PDE, the same recipe introduced in this review to derive the forward sensitivity equations and the continuous adjoint method for ODEs can be employed for PDEs \cite{Giles_Pierce_2000}. 
Assuming that the diffusivity $D = D(x, t; \theta)$ depends on some design parameter $\theta$, the forward sensitivity equation is obtained by differentiating Equation \eqref{eq:heat-equation} with respect to $\theta$. 
This defines a new PDE
\begin{equation}
 \frac{\partial s}{\partial t}
 = 
 D \, 
 \frac{\partial^2 s}{\partial x^2}
 + 
 \frac{\partial D}{\partial \theta} \frac{\partial^2 u}{\partial x^2}
\end{equation}
for the sensitivity $s(x,t) = \frac{d}{d\theta}u(x,t;\theta)$. 
The continuous adjoint equation can be derived using the strategy followed in Section \ref{section:continuous-adjoint} by multiplying the forward sensitivity equation by the transpose of the adjoint $\lambda (x,t)$ so that we can efficiently compute gradients of the objective function
\begin{equation}
    L(\theta) = \int_{t_0}^{t_1} \int_0^1 h(u(x,t;\theta); \theta) dx \, dt. 
\end{equation}
For the one dimensional heat equation in Equation \eqref{eq:heat-equation}, it is easy to derive using integration by parts the adjoint PDE given by 
\begin{equation}
    \frac{\partial \lambda}{\partial t}
    = 
    - 
    D \, \frac{\partial^2 \lambda}{\partial x^2}
    - 
    \frac{\partial h^T}{\partial u}
\end{equation}
with zero final condition $\lambda(x, t_1) \equiv 0$ and boundary conditions $\lambda(0, t) \equiv \lambda(1, t) \equiv 0$ \cite{duchateau1996introduction}.

An important consideration when working with PDEs is that meshing may be sensitive to model parameters which can lead to errors in the calculation of derivatives \cite{nadarajah2000comparison}. 
For example, this is a problem in finite differences since differences in the values of the objective function evaluated at $\theta$ and $\theta + \delta \theta$ can be affected by the choice of different meshes. 
The same errors induced by the adaptive stepsize controller in the case of AD (Section \ref{section:AD-incorrect}) can appear in cases of meshes that do not account for the joint error of the original PDE and its sensitivity. 
This can produce inaccurate gradients in the case of coarser meshes where the mesh or the numerical solver have an impact on the accuracy of the solution of the PDE \cite{economon2017adjoint, KENWAY2019100542}

Independently of how DP is implemented, PDEs remain some of the most challenging problems for computing sensitivities due to the frequent combination of a large number of discretized possibly stiff ODEs, with a large memory footprint. 
This makes it difficult to strike a balance between memory usage and computational performance. 
There are, however, numerous recent developments that have made solutions to these challenges more accessible. 
As will also be discussed in Section \ref{section:recomendations}, sensitivity methods that require the storage of a dense forward solution need special treatment, such as reverse AD and adjoint methods. 
Their huge memory footprint can be mitigated by using checkpointing (see Sections \ref{section:checkpointing} and \ref{section:checkpointint-cont}).
However, the memory requirements for even moderate size PDEs (e.g. $10^2$ to $10^3$ equations) over long time spans can still incur a large memory cost in cases where many checkpoints are required for stability in the reverse pass. 
This again can be mitigated by a multi-level checkpointing approach that enables checkpointing to either memory or to disk.
Another practical consideration when differentiating numerical PDE solvers arises from the way they are typically implemented. 
Due to the large size of the system, numerical calculations for PDEs are typically performed in-place, i.e. large memory buffers are often used to store intermediate calculations and system state thereby avoiding the need to repeatedly allocate large amounts of memory for each array operation. 
This can preclude the use of reverse AD implementations that do not support in-place mutation of arrays.
Automated sparsity detection \cite{gowdaSparsityProgrammingAutomated2019} and Newton-Krylov methods \cite{knollJacobianfreeNewtonKrylov2004,montoisonKrylovJlJulia2023} can drastically decrease both the time and space complexity of calculating JVPs or VJPs for large systems. 
Recent advances in applying AD to implicit functions, i.e. functions which require the solution of a nonlinear system, also provide a promising path forward for many complex PDE problems that often involve multiple nested numerical solvers \cite{blondelEfficientModularImplicit2022a}. 
Finally, some state-of-the-art AD tools such as Enzyme \cite{moses_Enzyme} are able to support both in-place modification of arrays as well as complex control flow, making them directly applicable to many high efficiency numerical codes for solving PDEs.

\subsection{Chaotic systems}

Continuous (nonlinear or infinite-dimensional) dynamical systems described by ODEs or PDEs can exhibit chaotic behavior~\cite{strogatz2018nonlinear}. 
In contrast to other systems we discussed previously, chaotic systems appear to become random after a certain, system-specific time scale, called the \textit{Lyapunov time}, making precise future predictions infeasible even though the underlying dynamical description might be completely deterministic.
In particular, such systems are characterized by their strong sensitivity to small perturbations of the parameters or initial conditions, i.e. small changes in the initial state or parameter can result in large differences in a later state, which is popularly known under the term \textit{butterfly effect}~\cite{PalmerButterfly2024}. 
As a consequence, all the sensitivity methods discussed in the previous sections become less useful when applied to chaotic systems and special considerations need to be taken under account \cite{Wang2012-chaos-adjoint}.

The butterfly effect makes inverse modelling based on point evaluations of the trajectory impractical. 
Therefore, we here resort to the loss function consisting in the long-time-averaged quantity 
\begin{equation}\label{eq:long_time_averaged_quantities}
    \langle L(\theta) \rangle_T = \frac{1}{T} \int_0^T h(u(t; \theta), \theta) \, dt, 
\end{equation}
where $h(u(t; \theta), \theta)$ is the instantaneous loss and $u(t; \theta)$ denotes the state of the dynamical system at time $t$.
In the presence of positive Lyapunov exponents, errors in solutions of the forward sensitivity equations and adjoint method to compute the gradient of $\langle L(\theta) \rangle_T$ with respect to $\theta$ blow up (exponentially fast) instead of converging to the actual gradient.
To address these issues, various modifications and methods have been proposed, including approaches based on ensemble averages~\cite{lea2000sensitivity, eyink2004ruelle}, the Fokker-Planck equation~\cite{thuburn2005climate, blonigan2014probability}, the fluctuation-dissipation theorem~\cite{leith1975climate, abramov2007blended, abramov2008new}, shadowing lemma~\cite{wang2013forward, wang2014least, wang2014convergence, ni2017sensitivity, blonigan2017adjoint, blonigan2018multiple, ni2019adjoint, ni2019sensitivity}, and modifications of Ruelle's formula~\cite{chandramoorthy2022efficient, ni2020fast}, which provides closed-form expressions and differentiability conditions for $\langle L(\theta) \rangle_T$ under the assumption of uniform hyperbolic systems \cite{ruelle1997differentiation,ruelle2009review}.


In Julia, the following methods based on the shadowing lemma are currently supported in the packages \texttt{AdjointLSS}, \texttt{ForwardLSS}, \texttt{NILSAS}, and \texttt{NILSS}. 
Standard derivative approximations are inappropriate for chaotic systems and will not give convergent estimates when the simulation time is a multiple of the Lyapunov time.

\subsection{Stochastic differential equations}
\label{section:sde}


DP plays an important role in stochastic differential equations (SDEs), which are used to describe dynamical systems with intrinsic source of randomness, including applications to particle systems (e.g., \textcite{pavliotis2014stochastic}), computational finance (see Section \ref{section:comp_finance}), biology and ecology (see Section \ref{section:biology}), and machine learning as generative models (e.g., \textcite{kidger2021neural}, \textcite{li2020scalable}, \textcite{wu2020stochastic}).
Given the randomness of the solution itself, SDEs control theory focuses on model parameter calibration and state estimation based on expected process quantities \cite{bertsekas2012dynamic}. 
For example, here we consider objective functions of the form $J(\theta) = \mathbb{E}[L(X_T)]$ where $X_t \in \mathbb{R}^d$ is a stochastic process solution of the SDE given by 
\begin{equation}
    dX_t = b(X_t, t, \theta) dt + \sigma(X_t, t, \theta) dW_t, \quad X_0 = x_0,
    \label{eq:ito}
\end{equation}
where $b: \mathbb{R}^d \times [0,T] \times \mathbb{R}^p \to \mathbb{R}^d$, $\sigma: \mathbb{R}^d \times [0,T] \times \mathbb{R}^p \to \mathbb{R}^{d \times m}$, $W_t \in \mathbb{R}^m$ is a Wiener process, and $\theta \in \mathbb{R}^p$ are parameters.

As in the case of ODEs, a number of different approaches enable gradient computation.
The \textit{pathwise continuous forward differentiation} method is the analogue of the sensitivity equations for SDEs. 
Assuming interchangeability of the differentiation and expectation~\cite{mohamed2020monte},
the pathwise sensitivity defined as $S_t = \partial X_t/\partial \theta \in \mathbb{R}^{d \times p}$ also follows an SDE given by
\begin{equation}
    dS_t = \left[\frac{\partial b}{\partial x}(X_t, t, \theta) S_t + \frac{\partial b}{\partial \theta}(X_t, t, \theta)\right] dt + \sum_{j=1}^m \left[\frac{\partial \sigma_{\cdot j}}{\partial x}(X_t, t, \theta) S_t + \frac{\partial \sigma_{\cdot j}}{\partial \theta}(X_t, t, \theta)\right] dW_{t,j},
\end{equation}
where all individual Jacobians and JVPs can be computed based on AD~\cite{yang1991monte, glasserman2004monte, tzen2019neural}.

On the other hand, the \textit{pathwise continuous adjoint differentiation} method defines an SDE for the adjoint variable \cite{li2020scalable}.
For Stratonovich SDEs given by
\begin{equation}
    dX_t = b(X_t, t, \theta) dt + \sigma(X_t, t, \theta) \circ dW_t,
\end{equation}
where the $\circ$ indicates the Stratonovich integral, the adjoint $a_t \in \mathbb{R}^d$ satisfies the Stratonovich SDEs given by
\begin{equation}
    da_t = -\frac{\partial b}{\partial x}(X_t, t, \theta)^T a_t dt - \sum_{j=1}^m \frac{\partial \sigma_{\cdot j}}{\partial x}(X_t, t, \theta)^T a_t \circ dW_{t,j},
\end{equation}
with $a_T = \nabla_x L(X_T)$, yielding to the final computation of the gradient
\begin{equation}
    \nabla_\theta J = \mathbb{E}\left[\int_0^T a_t^T \frac{\partial b}{\partial \theta} dt + \int_0^T \sum_{j=1}^m a_t^T \frac{\partial \sigma_{\cdot j}}{\partial \theta} \circ dW_{t,j}\right].
\end{equation}
We can convert between Stratonovich and It\^o SDEs using a standard conversion rule~\cite{Kloeden_Platen_1992}. 
Thus, the continuous adjoint method can also be applied to It\^o SDEs (Equation \eqref{eq:ito}) using the standard conversion rule twice, resulting in an additional term in the drift of the adjoint for It\^o SDEs.
This additional complexity in computing the adjoint in the It\^o sense makes the Stratonovich formulations usually preferable for adjoint-based optimization. 
In contrast to ODEs, we have to recompute or store the noise values used during the forward pass to compute the adjoint.

As in the case of ODEs, we can use AD to differentiate the SDE solver operations and obtain a sample-based estimate of the gradient. 
For instance, \textcite{Giles:2006vr} and \textcite{Innes_Zygote} applied reverse-mode AD through all intermediate steps of the SDE solver.
The time complexity of forward and adjoint differentiation scales roughly as $O(pd)$ and $O(p+d)$, respectively, per Monte Carlo path for non-stiff problems. 
Consequently, the adjoint method exhibits superior efficiency when $p$ is large.
A crucial advantage for DP is the independence of Monte Carlo paths, 
$\nabla_\theta J \approx \frac{1}{M} \sum_{i=1}^M g^{(i)}(\theta)$, 
which enables $O(d + p)$ memory per path rather than $O(M(d + p))$ total, parallel computation without inter-path communication, and streaming gradient accumulation.

For non-smooth objectives or when pathwise methods fail~\cite{mohamed2020monte}, the likelihood ratio method can be used ~\cite{glynn1990likelihood}.
This method is based on the decomposition of the gradient of the loss function as $\nabla_\theta J = \mathbb{E}[L(X_{T}) \nabla_\theta \log p_\theta(X_{T})]$.
A common strategy for evaluating $\log p_\theta(X_{T})$ is to use the Euler approximation, for which the (Gaussian) transition laws are explicitly available.
However, compared with pathwise estimates, the likelihood ratio method usually has a larger variance~\cite{glasserman2004monte}.

Common strategies to reduce variance for gradient estimators include control variates, multilevel Monte Carlo methods, and changes of measure.
Control variates exploit path independence:
\begin{equation}
    \nabla_\theta J = \mathbb{E}[(L(X_T) - c^T Y) \nabla_\theta \log p_\theta] + c^T \mathbb{E}[Y \nabla_\theta \log p_\theta],
\end{equation}
where $Y$ is an auxiliary random variable and $c$ minimizes variance, computable via least-squares regression across paths~\cite{glasserman2004monte}.
Multilevel Monte Carlo methods~\cite{giles2015multilevel, giles2018multilevel} are instead based on a telescoping sum representation:
\begin{equation}
    \mathbb{E}[L(X_T)] = \mathbb{E}[L(X_T^{(0)})] + \sum_{\ell=1}^L \mathbb{E}[L(X_T^{(\ell)}) - L(X_T^{(\ell-1)})],
\end{equation}
where $X_T^{(\ell)}$ uses time step $\Delta t_\ell = 2^{-\ell}T$, i.e., a sequence of approximations with increasing accuracy and cost, but fewer and fewer samples are required to accurately approximate the expectation of the difference as $\ell \rightarrow \infty$. Gradients inherit the same telescoping structure.
Finally, change of measure techniques via Girsanov's theorem provide a way for gradient estimation in stochastic systems by decoupling the randomness from parameter dependence. 
The fundamental advantage lies in transforming expectations $\mathbb{E}_{P_\theta}[L(X)]$ over parameter-dependent probability measures $P_\theta$ into expectations $\mathbb{E}_{P_0}[L(X) \frac{dP_\theta}{dP_0}]$ over a fixed reference measure $P_0$, enabling gradient computation without differentiating through stochastic dynamics. 
This approach proves particularly valuable when combined with variance reduction by choosing reference measures close to the target and then maintaining tractable variance in the likelihood ratio $\frac{dP_\theta}{dP_0}$. 
In stochastic optimal control, \textcite{hua2024simulation} exploit a reformulation using Girsanov's theorem~\cite{yang1991monte} to derive a computationally efficient gradient estimator for specific control problems that avoid both the high variance of naive likelihood ratios and the computational cost of differentiating through SDE solutions.

The DP for SDEs framework can be easily extended to more generic objectives, such as
\begin{equation}
    J(\theta) = \mathbb{E}\left[\int_0^T h(X_t) dt + L(X_T)\right].
\end{equation}
For discrete objectives, a common strategy is to combine the likelihood ratio method to handle discontinuities (since it avoids direct differentiation of the discontinuity) and pathwise methods to differentiate smooth parts, achieving a smaller variance than the pure likelihood ratio while maintaining unbiasedness~\cite{giles2009vibrato, mohamed2020monte}. 



\section{Recommendations}
\label{section:recomendations}

\begin{figure}[tb]
    \centering
    \includegraphics[width=1\textwidth]{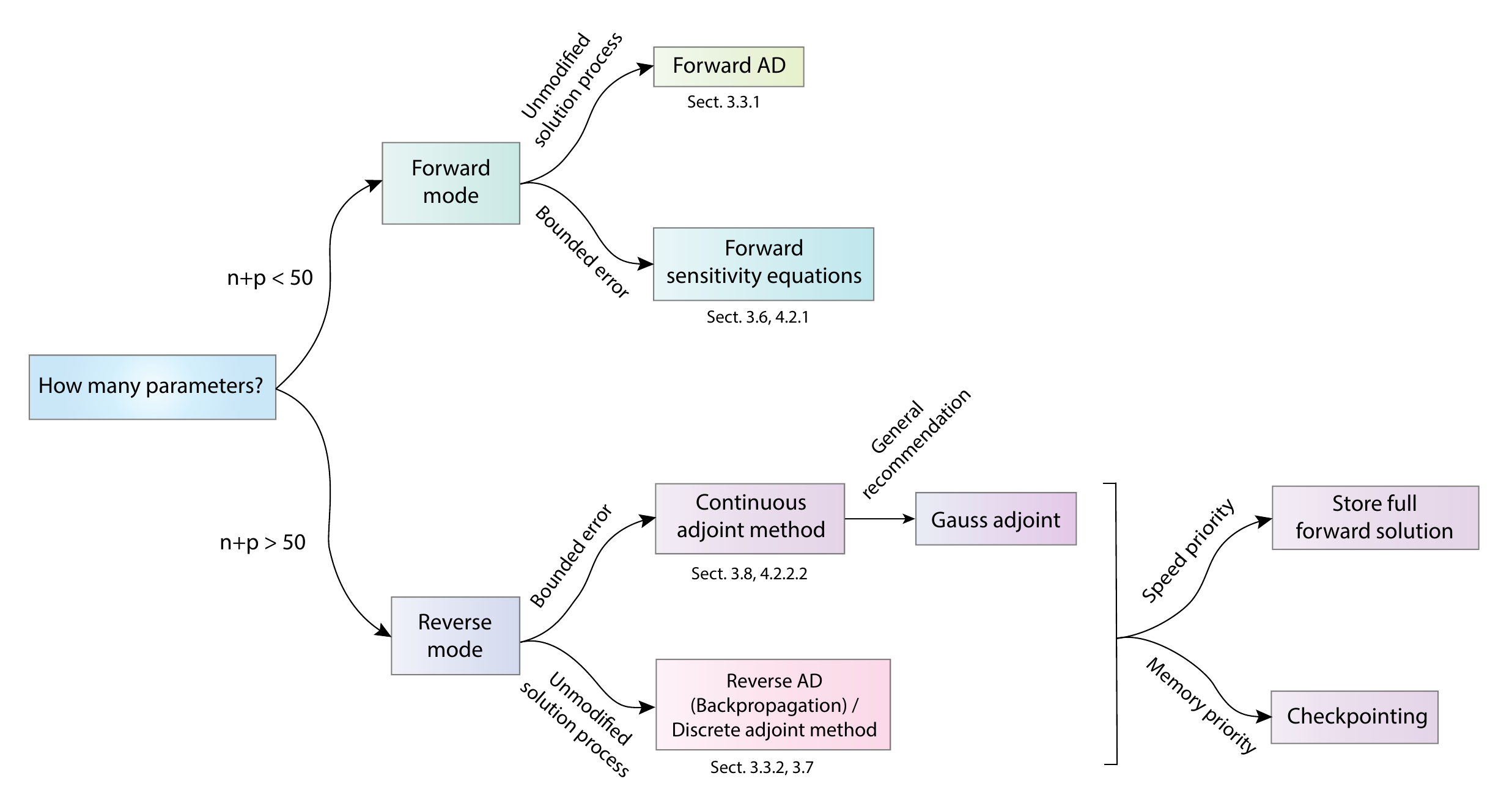}
    \caption{Decision-making tree summarizing the choice of sensitivity methods for different problems depending on: the number of parameters $p$, the number of ODEs $n$, the need for an unmodified solution during differentiation vs a bounded error (e.g. in the presence of a numerical solver to ensure correct gradients) and memory-speed trade-off.}
    \label{fig:roadmap}
\end{figure}

There is no sensitivity method that is universally suitable for all types of DE problems and that performs better under all conditions. 
However, in light of the methods we explore in this work, we can give general guidelines on which methods to use in specific circumstances.
In this section we provide practical guidance of which methods are the most suitable for different situations. 
A simplified overview of this decision-making process is depicted in Figure \ref{fig:roadmap}. 

\subsection*{Size of the system}

\subsubsection*{\textit{Working with small systems}}

For sufficiently small systems of less than $50$ parameters and ODEs, that is $n + p < 50$, it has been shown that forward AD and forward sensitivity equations are the most efficient methods, outperforming adjoint methods. 
The original benchmark of these methods is included in \textcite{ma2021comparison}, though the \texttt{SciMLBenchmarks} system continually updates the benchmarks and has revised the cutoff point as reverse AD engines improved. 
See \url{https://docs.sciml.ai/SciMLBenchmarksOutput/stable/} for continued updates.
Furthermore, as we have shown in Section \ref{section:direct-methods}, AD outperforms other forms of direct differentiation (finite differences, complex-step differentiation). 
Modern scientific software commonly supports AD, making forward AD the best choice for small problems. 

\subsubsection*{\textit{Working with large systems}}

For larger systems with more than 50 parameters plus the size of the ODE, reverse techniques are required. 
As explained in section \ref{section:solver-methods}, the continuous adjoint method, particularly the Gauss adjoint, and for very specific cases, the interpolating adjoint and quadrature adjoint, are the most suitable methods to tackle large stiff systems. 
The choice between these three types of adjoints will be problem-specific and will depend on the trade-off between numerical stability, performance and memory usage. 
Adjoint methods supporting checkpointing present more flexibility in this sense, and can allow modulating the method depending on the performance vs memory or input/output constraints of each problem. 

Unlike for small systems of ODEs and a reduced number of parameters, differentiating large ODEs (e.g. stiff discretized PDEs) with respect to a large number of parameters (e.g., in a neural network or in large-scale inversion), is a much more complex problem. 
Current state-of-the-art tools can easily work for a wide array of small systems, whereas methods for large systems are still under heavy development or require tailored approaches and are likely to see many changes and improvements in the future.
We refer to Table \ref{table:adjoint} for considerations of which adjoint method to use depending on the stability, performance, and memory trade-offs. 

~\\
\noindent
\textit{Special considerations for neural networks in DEs.} While the general guidelines for large systems hold for computing sensitivities for problems involving neural networks, it might be beneficial to use discrete methods when the training cost for the former is prohibitively high. 
For example, \textcite{Onken_Ruthotto_2020} demonstrates that discrete methods speed up neural ODE training by 6x. 
Additionally, \textcite{pal2021opening} show that discrete methods enable using information from the DE solver to speed up training and inference of neural ODEs and SDEs by 1.4x and 1.8x, respectively. 
These methods are available in \texttt{SciMLSensitivity.jl} via \texttt{TrackerAdjoint} and \texttt{DiffEqCallbacks.jl}.

\subsection*{Efficiency vs stability}

When using discrete methods, it is important to be aware that the differentiation machinery is applied after the numerical solver for the differential equation has been specified, meaning that the derivatives are computed with respect to the time discretization instead of the solution \cite{Eberhard_Bischof_1996}.
As discussed in Section \ref{section:software-Forward-AD}, this can mean the method is non-convergent in the case where the iterative solver has adaptive stepsize controllers that depend on the parameter to differentiate.
Although some solutions have been proposed and implemented in Julia to solve this in the case of discrete methods \cite{Eberhard_Bischof_1996}, this is a problem that continuous methods do not have since they apply the differentiation step before the numerical algorithm has been specified. 
Using many of the aforementioned tricks, such as continuous checkpointing and Gaussian quadrature approximations, continuous sensitivity analysis tends to be more memory and computationally efficient. 
However, the discrete adjoint method's derivative error may better represent the actual code being evaluated. 
For this reason, discrete adjoint method has been found in some instances to lead to more stable optimizations. 

In a nutshell, continuous adjoint methods tend to be more efficient while discrete adjoint methods tend to be more stable, though the opposite can apply and as such the choice ultimately depends on the nuances of each problem.
This is reflected in the fact that discrete methods usually differentiate the unmodified solution of the original ODE, while continuous methods adapt the solution of the original ODE and the sensitivity/adjoint to control for their joint numerical error.

\subsection*{Choosing a direct method}

When computing the gradient of a generic function other than a numerical solver, we further recommend the use of AD (reverse or forward depending the number of parameters) as the direct method of choice, outperforming finite differences, complex step differentiation, and symbolic differentiation. 
This recommendation also applies for the inner JVPs and VJPs calculations performed inside the numerical solver (Section \ref{section:computing-vjp-inside-solver}).
As discussed in Section \ref{section:direct-methods}, finite differences and complex step differentiation do not really provide an advantage over AD in terms of precision and require the tuning of the stepsize $\varepsilon$. 
On the other hand, if symbolic differentiation can be more efficient in nested cases or when the sparsity pattern of the Jacobian is known, in general this advantage is not drastic in most real cases and can generate difficulties when used inside the numerical solver. 

However, this recommendation is constrained by the availability and interoperability of different AD and sensitivity software. 
For example, when computing higher-order derivatives multiple layers of direct methods became more difficult to implement and may result in complicated computer programs.
In this case, complex step differentiation may offer an interesting alternative with similar performance than AD for small stepsizes. 
It is important to mention that incorrect implementations of both forward and reverse AD can lead to \textit{perturbation confusion}, an existing problem in some AD software where either repeated applications of AD or differentiation with respect to different dual variables result indistinguishable \cite{siskind2005perturbation, manzyuk2019perturbation}. 

\subsubsection*{Taking into account model architecture}

Code structure and characteristics have a very strong impact on the choice of which packages to use to compute the sensitivities. 
Within the Julia and Python ecosystems, each available AD package implements a specific AD technique that will face certain limitations.
Current limitations include:
\begin{enumerate}
    \item[$ \blacktriangleright$] The use of control flow (i.e. \texttt{if}/\texttt{else} statements; \texttt{for} and \texttt{while} loops) presents issues for dynamic (tape-based) AD methods (see Section \ref{sec:software-reverse-AD}). 
    This is currently not supported by \texttt{ReverseDiff.jl} (with tape compilation) and partially supported by \texttt{JAX} in Python. 
    Non-tape-based AD methods tend to support this, like \texttt{Enzyme.jl} and \texttt{Zygote.jl}.
    \item[$ \blacktriangleright$] Mutation of arrays (i.e. in-place operations) is sometimes problematic, since it does not allow the preservation of the chain rule during reverse differentiation. As such, mutations are not possible for packages like \texttt{Zygote.jl} or \texttt{JAX}. It is however currently supported by \texttt{ReverseDiff.jl} and \texttt{Enzyme.jl}.
    \item[$ \blacktriangleright$] Compatibility with GPUs (Graphical Processing Units) is still greatly under development for sensitivity methods. 
    Certain AD packages like \texttt{ReverseDiff.jl} do not support GPU operations, while others like \texttt{JAX}, \texttt{Enzyme.jl} and \texttt{Zygote.jl} support it. 
    This makes the former unsuitable for problems involving large neural networks (e.g, neural ODEs \cite{chen_neural_2019}) that rely on GPUs for scalability.
\end{enumerate}

It is important to bear in mind that direct methods are easier to implement in programming languages where AD already exists and sometimes does not require any special package, like for the Julia programming language.
Nonetheless, users must be aware of the aforementioned convergence issues of AD naively applied to solvers. 
Thus, we recommend the use of robust and tested software when available (e.g., the Julia SciML ecosystem or Diffrax in Python), as the solvers must apply corrections to AD implementations in order to guarantee numerically correct derivatives.

\section{Conclusions}
We have presented a comprehensive overview of the different existing methods for calculating the sensitivity or gradients of functions, including loss functions, involving numerical solutions to differential equations.
This task has been approached from three different angles.
First, we surveyed the existing literature in different scientific communities where differential programming tools have been used before and play a central modelling role, especially for inverse modeling.
Next, we reviewed the mathematical foundations of these methods and their classification as forward vs reverse, discrete vs continuous, and mixed approaches.
We further compared the mathematical and computational foundations of these methods with an aim to enlighten the discussion on sensitivity methods and to demystify  misconceptions around the sometimes apparent differences between methods.  
We showed how these methods can be translated to software implementations, evaluating considerations that we must take into account when implementing or using a sensitivity algorithm. 
We further exemplified how these methods are implemented in the easy-to-read Julia programming language. 

There are challenges that the next generation of differentiable programming methods will have to address, in particular in the context of large-scale DE-based modelling. 
Among them, we highlight the development of general-purpose source-to-source AD,
the use of control flows and in-place operations, GPU support, parallelization of algorithms, dealing with strong linearities, and hybrid machine learning-PDE based approaches. 
These features are of particular importance for PDE-based inverse modelling due to the combination of complex systems of equations with a large memory footprint.  
These improvements also have an impact on the performance of reverse methods, which is why benchmarks for these continue to improve over time. 
Furthermore, we expect the relative performance between  methods to change over time due to the development of new reverse-mode (adjoint or backpropagation) methods that trade off accuracy, time, and memory usage. 

There exist a myriad of options and combinations to compute sensitivities of functions involving differential equations, further complicated by the jargon and scientific culture in different communities. 
We hope this review provides a clearer overview of the subject, provides a bridge across different communities, and can serve as an entry point to navigate this field and guide researchers in choosing the most appropriate method for their scientific application.

Differentiable programming is opening new ways of doing research across different domains of science and engineering. 
Arguably, its potential has so far been somewhat under-explored but is resurging in the age of data-driven science. 
Realizing its full potential requires collaboration between domain scientists, computational scientists, computer scientists, and applied mathematicians in order to develop successful, scalable, practical, and efficient frameworks for real-world applications. 
As we make progress in the development and use of these tools, new methodological challenges and opportunities will emerge. 


\vspace{30px}
\section*{Software availability}
All the scripts and code shown in this paper can be found in the GitHub repository \url{https://github.com/ODINN-SciML/DiffEqSensitivity-Review}. 
Examples of available code are indicated in the manuscript with the symbol $\clubsuit$.
See Appendix \ref{appedix:code} for a complete description of the scripts provided. 
\section*{Contribution statement}
\noindent
The following categories are based on the \href{https://credit.niso.org/}{Contributor Roles Taxonomy (CRediT)}. 
FSap: conceptualization, investigation, project administration, software, visualization, writing - original draft.
JB: conceptualization, visualization, writing - review and editing.
FSch: conceptualization, investigation, software, writing - review and editing.
BG: conceptualization, software, writing - review and editing.
AP: software, writing - review and editing.
VB: conceptualization, writing - review and editing.
PH: conceptualization, investigation, supervision, writing - review and editing.
GH: conceptualization, supervision, writing - review and editing.
FP: conceptualization, funding acquisition, supervision, writing - review and editing.
PP: conceptualization, supervision, writing - review and editing.
CR: conceptualization, funding acquisition, investigation, software, supervision, writing - review and editing.
\section*{Acknowledgments}

FSap would like to thank Jonathan Taylor, Ryan Giordano, Alexander Strang, and Olivier Bonte for useful comments and feedback. 
All authors acknowledge the constructive feedback of anonymous reviewers. 
FSap and FP benefits from the Jupyter meets the Earth project supported by the NSF Earth Cube Program under awards 1928406, 1928374. 
JB has been supported by the Nederlandse Organisatie voor Wetenschappelijk Onderzoek, Stichting voor de Technische Wetenschappen (Vidi grant 016.Vidi.171.063).
BG acknowledges the support of the Helmholtz Einstein International Berlin Research School in Data Science as well as the Center for Scalable Data Analytics and Artificial Intelligence Dresden/Leipzig.
PH was supported in part by NSF CSSI \#OAC-2103942, ONR \#N00014-20-1-2772, DOE \#DE-SC002317 \& \#DE-AC02-06CH11357, and a JPL/Caltech subcontract of the NASA Estimating the Circulation and Climate of the Ocean (ECCO) project.
GH was supported by the NSF grant DEB-1933497. 
PP was supported in part by the U.S. Department of Energy, Office of Science, Office of Advanced Scientific Computing Research's Applied Mathematics Competitive Portfolios program under Contract No. AC02-05CH11231, and in part by the National Science Foundation under Grant DMS-2309596.

Fsch, AP, and CR works was supported in part by the Director, Office of Science, Office of Advanced Scientific Computing Research, U.S. Department of Energy under Contract No. DE-AC02-05CH11231, and in part by the National Science Foundation under Grant DMS-2309596.
This material is based upon work supported by the Department of Energy, National Nuclear Security Administration under Award Number DE-NA0003965.
This report was prepared as an account of work sponsored by an agency of the United States Government. 
Neither the United States Government nor any agency thereof, nor any of their employees, makes any warranty, express or implied, or assumes any legal liability or responsibility for the accuracy, completeness, or usefulness of any information, apparatus, product, or process disclosed, or represents that its use would not infringe privately owned rights. 
Reference herein to any specific commercial product, process, or service by trade name, trademark, manufacturer, or otherwise does not necessarily constitute or imply its endorsement, recommendation, or favoring by the United States Government or any agency thereof. 
The views and opinions of authors expressed herein do not necessarily state or reflect those of the United States Government or any agency thereof. 
This material is based upon work supported by the National Science Foundation under grant no. OAC-2103804 , no. OSI-2029670, no. DMS-2325184, no. PHY-2028125. 
This material was supported by The Research Council of Norway and Equinor ASA through the Research Council project ”308817 - Digital wells for optimal production and drainage”. 
Fsch, AP, and CR research was sponsored by the United States Air Force Research Laboratory and the United States Air Force Artificial Intelligence Accelerator and was accomplished under Cooperative Agreement Number FA8750-19-2-1000. 
The views and conclusions contained in this document are those of the authors and should not be interpreted as representing the official policies, either expressed or implied, of the United States Air Force or the U.S. Government. 
The U.S. Government is authorized to reproduce and distribute reprints for Government purposes notwithstanding any copyright notation herein. 
Fsch, AP, and CR  would like to thank DARPA for funding this work through the Automating Scientific Knowledge Extraction and Modeling (ASKEM) program, Agreement No. HR0011262087. The views, opinions and/or findings ex-pressed are those of the authors and should not be interpreted as representing the official views or policies of the Department of Defense or the U.S. Government.


\newpage
\appendix
\section*{Appendices}
\addcontentsline{toc}{section}{Appendices}
\renewcommand{\thesubsection}{\Alph{subsection}}

\subsection{Supplementary code}
This is a list of the code provided along with the current manuscript.
All the following scripts can be found in the GitHub repository \href{https://github.com/ODINN-SciML/DiffEqSensitivity-Review}{\texttt{DiffEqSensitivity-Review}}. 
 
\begin{itemize}
    \item[$\clubsuit_\text{\ref{code:figure-comparison}}$] \textbf{Comparison of direct methods.} The script \url{https://github.com/ODINN-SciML/DiffEqSensitivity-Review/blob/main/code/DirectMethods/Comparison/direct-comparision.jl} reproduces Figure \ref{fig:direct-methods}.
    
    \item[$\clubsuit_\text{\ref{code:dual-number}}$] \textbf{Dual numbers definition.} The script \url{https://github.com/ODINN-SciML/DiffEqSensitivity-Review/blob/main/code/DirectMethods/DualNumbers/dualnumber_definition.jl} includes a very simple example of how to define a dual number using \texttt{struct} in Julia and how to extend simple unary and binary operations to implement the chain rule usign multiple distpatch. 
    
    \item[$\clubsuit_\text{\ref{code:AD-wrong}}$] \textbf{When AD is algorithmically correct but numerically wrong.} The script \url{https://github.com/ODINN-SciML/DiffEqSensitivity-Review/blob/main/code/SensitivityForwardAD/example-AD-tolerances.jl} includes the example shown in Section \ref{section:software-Forward-AD} and further elaborated in Section \ref{section:AD-incorrect} where forward AD gives the wrong answer when tolerances in the gradient are not computed taking into account both numerical errors in the numerical solution and the sensitivity matrix. Further examples of this phenomena can be found in the  and the Julia \url{https://github.com/ODINN-SciML/DiffEqSensitivity-Review/blob/main/code/SensitivityForwardAD/testgradient_julia.jl}.
    
    \item[$\clubsuit_\text{\ref{code:AD-wrong-JAX}}$] \textbf{When AD is algorithmically correct but numerically wrong (JAX)}. Python script \url{https://github.com/ODINN-SciML/DiffEqSensitivity-Review/blob/main/code/SensitivityForwardAD/testgradient_python.py}.
    \item[$\clubsuit_\text{\ref{code:complex-step}}$] \textbf{Complex step in numerical solver.} The script \url{https://github.com/ODINN-SciML/DiffEqSensitivity-Review/blob/main/code/DirectMethods/ComplexStep/complex_solver.jl} shows how to define the dynamics of the ODE to support complex variables and then compute the complex step derivative. 
    
    \item[$\clubsuit_\text{\ref{code:sensitivity-equation}}$] \textbf{Forward sensitivity equation. } The scrip \url{https://github.com/ODINN-SciML/DiffEqSensitivity-Review/blob/main/code/SolverMethods/Harmonic/forward_sensitivity_equations.jl} includes a manual implementation of the  forward sensitivity equations. This also includes how to compute the same sensitivity using \texttt{ForwardSensitivity} in Julia. 
    \item[$\clubsuit_\text{\ref{code:discrete-adjoint}}$] \textbf{Discrete adjoint method. }The script \url{https://github.com/ODINN-SciML/DiffEqSensitivity-Review/blob/main/code/SolverMethods/Harmonic/adjoint_discrete.jl} includes a manual implementation of the discrete adjoint method for the simple harmonic oscillator. 
    
    \item[$\clubsuit_\text{\ref{code:continuous-adjoint}}$] \textbf{Continuous adjoint method. }The script \url{https://github.com/ODINN-SciML/DiffEqSensitivity-Review/blob/main/code/SolverMethods/Harmonic/adjoint_continuous.jl} includes a manual implementation of the continuous adjoint method for the simple harmonic oscillator. 
\end{itemize}
\label{appedix:code}


\newpage
\printbibliography[heading=bibintoc, title={References}]

\end{document}